\documentclass[a4paper]{article}

\usepackage{MyLaTeX,MyArticle}
\usepackage{multirow}

\RequirePackage{pgf,pgffor}
\usepackage{tikz}
\newcommand{\cf}{\mathrm{cf}}
\newcommand{\Irr}{\mathrm{Irr}}
\newcommand{\AGL}{\mathrm{AGL}}
\newcommand{\vmin}{V_{\min}}
\newcommand{\Vmin}{V_{\min}}
\renewcommand{\top}{\mathrm{top}}
\newcommand{\minitab}[2][l]{\begin{tabular}{#1}#2\end{tabular}}

\numberwithin{equation}{section}
\title{Alternating Subgroups of Exceptional Groups of Lie Type}
\author{David A.\ Craven}
\date{\today}

\tikzstyle{every node}=[circle, fill=black!0,
                        inner sep=0pt, minimum width=4pt]

\renewcommand{\thedefn}{\arabic{defn}}
\begin{document}
\maketitle

\begin{abstract}
In this paper we examine embeddings of alternating groups and symmetric groups into almost simple groups of exceptional type. In particular, we prove that unless the alternating or symmetric group has degree 6 or 7, there is no maximal subgroup of any almost simple group with socle an exceptional group of Lie type that is an alternating or symmetric group. Furthermore, in the remaining open cases we give considerable information about the possible embeddings. Note that no maximal alternating or symmetric subgroups are known in the remaining cases.

This is the first in a sequence of papers aiming to substantially improve the state of knowledge about the maximal subgroups of exceptional groups of Lie type.
\end{abstract}

The maximal subgroups of the finite simple groups have been the object of intense study over the last few decades. For the exceptional groups of Lie type it seems plausible -- unlike for example in the classical group case -- to have a complete list of all maximal subgroups. If $G=G(q)$ is an exceptional group of Lie type, then all maximal subgroups of $G$ are known, except for an explicit, finite, list of candidates, all of which are almost simple.

Techniques for finding maximal subgroups of the larger-rank exceptional groups have until now been largely geometric (for example, \cite{magaardphd} and \cite{aschbacherE6Vun} for $F_4$ and $E_6$) or relate to algebraic groups (for example, \cite{liebecksaxl1987}, \cite{lst1996} and \cite{liebeckseitz1999}). Independently the author, in a series of papers of which this is the first, and Alastair Litterick \cite{litterickmemoir}, introduced methods from modular representation theory to attack the problem. These seem much more suited to the remaining, difficult, open cases of small almost simple groups being embedded in large exceptional groups and offer hope that a complete solution can be obtained. Representation theory is predominantly used to prove that maximal subgroups do not exist, and other methods will need to be used to prove existence and uniqueness of maximal subgroups that do, in fact, exist.

The author's rough estimation is that, from the list of possible isomorphism types of maximal subgroups that have not yet been eliminated or found, these methods can be used to solve 90\% of them or more, leaving a small number of possibilities. These methods really apply when the characteristic of the field over which the group of Lie type is defined, say $p$, divides the order of the potential maximal subgroup $H$. As an example of their power, one can use them to eliminate $\PSL_2(7)$ as a maximal subgroup of $E_8(3^a)$ for all $a$ in one line, using the tables in \cite{litterickmemoir} as a starting point.

In keeping with this theme of vastly reducing -- but not eliminating entirely -- the possibilities, we do not prove that there can be no maximal alternating or symmetric subgroups in exceptional groups, but we do prove this in almost all cases. Our strongest theorems are slightly more technical, so we begin with a weaker theorem that serves for many purposes. We should state that no maximal alternating or symmetric subgroup of an exceptional group of Lie type is known at the moment, and so the next theorem merely lists the possibilities not yet excluded.

\begin{thm}\label{thm:summary} Let $G$ be an almost simple exceptional group of Lie type over a field of characteristic $p$, and let $H$ be a maximal subgroup of $G$ such that $F^*(H)=\Alt(n)$ for some $n\geq 5$. Then $n=6$ or $n=7$, and one of the following holds:
\begin{enumerate}
\item $n=6$ and $(G,p)$ is one of $(F_4,3)$, $(F_4,\geq 7)$, $(E_6,\geq 7)$, $(E_7,5)$, $(E_7,\geq 7)$, $(E_8,2)$, $(E_8,3)$, $(E_8,\geq 7)$
\item $n=7$ and $(G,p)$ is one of $(E_7,5)$, $(E_8,3)$, $(E_8,5)$, $(E_8,7)$, $(E_8,\geq 11)$
\end{enumerate}
\end{thm}

In this theorem, $E_6$ means either $E_6$ or ${}^2\!E_6$, and $\geq p$ means all primes at least $p$, i.e., the case where $p\nmid |H|$, which includes characteristic $0$. This theorem is a summary of more complicated, and stronger, theorems about embeddings of alternating groups into simple exceptional groups of Lie type. In addition, because we will always prove that there exists a positive-dimensional subgroup of the algebraic group fixing certain specific subspaces that the potential maximal subgroup also fixes (normally lines on the minimal or adjoint module), using results from Section \ref{sec:intro} the above restrictions also apply to any almost simple group whose socle is an exceptional group, giving the full statement of the theorem above.

Let $\mathbb G$ be an exceptional algebraic group, let $\ms X$ be the set of all maximal positive-dimensional subgroups of $\mathbb G$, and let $\ms X^\sigma$ denote the corresponding subgroups of the finite group $G=\mathbb{G}^\sigma$ (i.e., the fixed points under the Frobenius endomorphism $\sigma$ defining $\mathbb G^\sigma$). We prove that, with the exceptions above, $H$ is always contained inside a member of $\ms X^\sigma$. This is stronger than simply stating that $H$ is not maximal, since there are other potential maximal subgroups of $G$ that could contain $H$. For example, inside $E_6$ in characteristic $0$, there is a copy of ${}^2\!F_4(2)'$ acting irreducibly on the minimal module, and this contains a copy of $\PSL_3(3)$. This subgroup is clearly not maximal, but is not contained in any member of $\ms X$. In addition, since we understand the set $\ms X^\sigma$, we can get some handle on the possible embeddings of $H$ into $G$.

\begin{thm}\label{thm:alt5}
Let $G=\mathbb{G}^\sigma$ be a simple group of type $F_4$, $E_6$, ${}^2\!E_6$, $E_7$ or $E_8$. If $H$ is a subgroup of $G$ with $F^*(H)\cong\Alt(5)$ then $H$ lies inside a member of $\ms X^\sigma$.
\end{thm}

For the alternating groups $\Alt(6)$ and $\Alt(7)$ much less is known. We summarize what we do know now. In what follows we denote simple modules by their dimension, distinguishing between modules of the same dimension with an index, so $4_2$ is a module of dimension $4$. $P(-)$ denotes the projective cover and we delineate between socle layers with a `$/$' character.

\begin{thm}\label{thm:alt6}
Let $G=\mathbb{G}^\sigma$ be a simple group of type $F_4$, $E_6$, ${}^2\!E_6$, $E_7$ or $E_8$. Let $H$ be a subgroup of $G$ with $F^*(H)\cong\Alt(6)$, and suppose that $H$ does not lie inside a member of $\ms X^\sigma$. Then one of the following holds:
\begin{enumerate}
\item $p=2$ and $G=E_8$, with $F^*(H)$ acting on $L(G)$ with composition factors $8_1^6,8_2^6,4_1^{16},4_2^{16},1^{24}$;
\item $p=3$ and $G=F_4,E_8$;
\item $p=5$ and $G=E_7$, with $F^*(H)$ acting on $V_{\min}$ and on $L(G)$ as
\[ 10^{\oplus 4}\oplus 8^{\oplus 2}\quad\text{and}\quad 10^{\oplus 2}\oplus 5_1^{\oplus 3}\oplus 5_2^{\oplus 3}\oplus P(8)^{\oplus 3}\oplus 8\]
respectively, or $2\cdot F^*(H)$ lying in $E_7$, acting on $V_{\min}$ and $L(G)$ as
\[10_2^{\oplus 3}\oplus 10_3\oplus 4_1/4_2\oplus 4_2/4_1\quad\text{and}\quad 10^{\oplus 5}\oplus P(8)^{\oplus 3}\oplus 8\]
respectively;
\item $p>5$ and $G=F_4,E_6,E_7,E_8$.
\end{enumerate}
\end{thm}

\begin{thm}\label{thm:alt7}
Let $G=\mathbb{G}^\sigma$ be a simple group of type $F_4$, $E_6$, ${}^2\!E_6$, $E_7$ or $E_8$. Let $H$ be a subgroup of $G$ with $F^*(H)\cong\Alt(7)$, and suppose that $H$ does not lie inside a member of $\ms X^\sigma$. Then one of the following holds:
\begin{enumerate}
\item $p=3$ and $G=E_8$, with $F^*(H)$ acting on $L(G)$ with composition factors $15^4,13^6,(10,10^*)^5,1^{10}$;
\item $p=5$ and $G=E_7$, with $F^*(H)$ acting on $V_{\min}$ as either $(10\oplus 10^*)^{\oplus 2}\oplus 8^{\oplus 2}$ or $20\oplus 4/14\oplus 14/4^*$, or $G=E_8$, with $F^*(H)$ acting on $L(G)$ as
\[ 35^{\oplus 4}\oplus 15^{\oplus 4}\oplus 10\oplus 10^*\oplus 8/6\oplus 6/8;\]
\item $p=7$ and $G=E_8$, with $F^*(H)$ acting on $L(G)$ as
\[ 35^{\oplus 3}\oplus 21\oplus 14_1\oplus 14_2^{\oplus 2}\oplus P(10)^{\oplus 2}\oplus 10\quad\text{or}\quad 35^{\oplus 4}\oplus 14_1^{\oplus 2}\oplus 10^{\oplus 6}\oplus 5^{\oplus 4};\]
\item $p>7$ and $G=E_8$, acting on $L(G)$ with composition factors $35^4,15^4,14_1^2,10,10^*$.
\end{enumerate}
\end{thm}

Once we move past $\Alt(7)$, however, we have complete theorems. For $\Alt(8)$ and above we can prove that every copy lies inside a positive-dimensional subgroup. This includes the difficult case of $\Alt(8)=\GL_4(2)$ inside $E_8(2^n)$, which starts to extend to $E_8$ the results on defining-characteristic embeddings of groups of rank at least $2$ considered for $F_4$, $E_6$ and $E_7$ in \cite{cmp2015un}. Extending this theorem to include $E_8$ is ongoing work of the author.

\begin{thm}\label{thm:alt8+}
Let $G=\mathbb{G}^\sigma$ be a simple group of type $F_4$, $E_6$, ${}^2\!E_6$, $E_7$ or $E_8$. If $H$ is a subgroup with $F^*(H)=\Alt(n)$ for $n\geq 8$, then $H$ lies inside a member of $\ms X^\sigma$.
\end{thm}

In \cite[Section 3.3]{litterickmemoir}, Litterick proves that, for $n\geq 10$, if $F^*(H)\cong \Alt(n)$ then $H$ is never maximal in an almost simple exceptional group. We include an alternative proof of the case $n=10$, $p=2$, $G=E_8$ in the last section, as we show that $H$ fixes a line on $L(G)$ more easily with a theoretical argument involving Frobenius reciprocity.

In addition, Litterick has proved many other cases for alternating groups, in his more wide-ranging project that included all simple groups other than Lie type in defining characteristic. We refer to \cite{litterickmemoir}, as yet not in publication but available on the arXiv, for full details as to the cases proved there, which are  more than in his PhD thesis. We have however mostly maintained our own proofs here both for completeness and as our methods are slightly different.

These results include and extend known results on maximal subgroups of exceptional groups. For $p\geq 5$ and $G=F_4$ we recover results of Magaard in \cite{magaardphd}. For $G=E_6$ we extend results of Aschbacher in \cite{aschbacherE6Vun}, where we remove the following possibilities for maximal subgroups: 
\begin{itemize}
\item $\Alt(6)$, $p=2$;
\item $\Alt(6)$, $p=5$; 
\item $\Alt(7)$, $p=5$;
\item $\Alt(8)$, $p=2$.
\end{itemize}
We leave $\Alt(6)$, $p\geq 7$ (or $p=0$) unresolved.

\bigskip

The method of proof here is as follows: we firstly use traces of semisimple elements to restrict the possible sets of composition factors for the action of $H$ on both the minimal and adjoint modules. This task was already accomplished in \cite{litterickmemoir} when $H$ is not isomorphic to a Lie type group in the same characteristic, so we only have to consider those cases.

Once we have these data, we then use the information about the action of unipotent elements on the minimal and adjoint modules in \cite{lawther1995} to build up a picture of the possible restrictions of the minimal and adjoint modules to the subgroup. Often this is enough to prove something like the subgroup fixes a line on one of these two modules, or that the set of composition factors cannot yield an embedding of the subgroup.

We also use a result from \cite{lst1996} that states that for certain unipotent classes, which we call `generic' classes in the next section, merely containing an element from them and not acting irreducibly on the minimal or adjoint module is enough to guarantee that the subgroup lies in the set $\ms X$ above. This is where many embeddings do lie, and so it is useful to exclude large numbers of embeddings from consideration early on.

If the group has a cyclic Sylow $p$-subgroup then we have complete information about the module category for the group, and so we can get very good information on the possible embeddings. Potential embeddings that seem not to be attackable with using the techniques here are listed in Theorems \ref{thm:alt6} and \ref{thm:alt7}. Other ideas such as using structure constants, which have been used before to success, might be one avenue for these, and of course the characteristic $0$ possibilities, as they have been successful before, for example for $\Alt(5)$ \cite{lusztig2002}.

Many of the arguments for $p=2,3$ from later sections require use of a computer, however. This is firstly to prove the existence or non-existence of certain modules, and secondly to determine the actions of unipotent elements on these modules. While it might be possible in many cases to perform this analysis without a computer, doing so would substantially increase both the size of the paper and the likelihood of mistakes appearing. The author believes that the probability of there being bugs in Magma that invalidate the results here is lower than the chances of making an error in hundreds of pages of module computations.

Having said that, despite the author's use of a computer to produce various modules, many claims can be checked using the structures of projective indecomposable modules for the small alternating groups in characteristics $2$ and $3$, that have appeared in the literature.

\bigskip

The structure of this paper is as follows: we begin with a section collating the notation we need, and the preliminary results, including a result on passing from maximal subgroups of the simple group to those of an almost simple group. We then briefly summarize those aspects of the theory of blocks with cyclic defect groups that we need, particularly how the structure of the Brauer tree informs the structure of the projective indecomposable modules. In Secton \ref{sec:modulesalt} we give information about simple modules for alternating groups, in particular their dimensions and $\Ext^1$, the Brauer trees of $\Alt(n)$ for $p=5,7$ and $n\leq 8$, and some useful lemmas for specific groups. The succeeding section classifies the composition factors of large-degree alternating groups on the adjoint module for $E_8$ in characteristic $2$.

After these preliminary sections, Section $n$ considers $\Alt(n)$, for $5\leq n\leq 10$.

\renewcommand{\thedefn}{\thesection.\arabic{defn}}

\section{Notation and preliminary results}
\label{sec:intro}

Throughout this paper, to avoid confusion with algebraic groups, we denote by $\Alt(n)$ the alternating group on $n$ letters, and so to remain consistent use $\Sym(n)$ for the symmetric group on $n$ letters.

Although we will remind the reader often, we establish a consistent naming convention for the various permutations we will need here:

\begin{center}
\begin{tabular}{cc}
\hline Label & Permutation
\\ \hline
$t$ & $(1,2)(3,4)$
\\ $u$ & $(1,2,3,4,5)$
\\ $v$ & $(1,2,3,4)(5,6)$
\\ $w$ & $(1,2,3,4,5,6,7)$
\\ $x$ & $(1,2,3)$
\\ $y$ & $(1,2,3)(4,5,6)$
\\ $z$ & $(1,2,3,4,5,6,7,8,9)$
\\ \hline
\end{tabular}
\end{center}

Let $\mathbb G$ be a simple algebraic group, and let $\sigma$ be a Frobenius endomorphism of $\mathbb G$. The finite group $G=\mathbb G^\sigma$ is a finite group of Lie type. We can also think of $G$ as $G(q^\delta)$, a finite group of Lie type defined over the field $\F_{q^{\delta}}$. In our case $G=\mathbb G^\sigma=G(q^\delta)$ is one of ${}^2\!B_2(q^2)$, ${}^3\!D_4(q^3)$, $E_6(q)$, ${}^2\!E_6(q^2)$, $E_7(q)$, $E_8(q)$, $F_4(q)$, ${}^2\!F_4(q^2)$, $G_2(q)$ and ${}^2\!G_2(q^2)$, or their universal central extensions. Since we are examining maximal subgroups, and the maximal subgroups are known for all of the above groups except for those of twisted rank at least four, we assume that $G$ is one of $F_4(q)$, $E_6(q)$, ${}^2\!E_6(q^2)$, $E_7(q)$ and  $E_8(q)$. (The maximal subgroups of the remaining groups are given in \cite{wilsonrob}, which includes references to the original papers.) We will assume that by $E_6$ and $E_7$ we mean the simply connected version, so $Z(E_6)$ and $Z(E_7)$ have orders $3$ and $2$ respectively. (We will remind the reader of this regularly.)

Let $p\mid q$ be a prime, and let $k$ be a field of characteristic $p$. We specifically do not take $k$ to be algebraically closed, because embeddings of (say) $H$ into $F_4(2)$ produce $\F_2H$-modules, not modules over the algebraically closed field. This is not important when proving that a subgroup fixes a line on a module, since this is true regardless of the field over which the module is considered, but is important when asking whether, for example, $\SL_2(4)$ stabilizes a $2$-space on a module.

If $M$ is a module for a group, $\soc(M)$, the \emph{socle} of $M$, is the largest semisimple submodule of $M$, and $\top(M)$ is the largest semisimple quotient of $M$. Write $\uparrow$ and $\downarrow$ for induction and restriction, and write $k_H$ or $1_H$ for the trivial module for the group $H$. If the group $H$ is obvious we will simply write $k$ or $1$. The kernel of the map $M\to \top(M)$ is denoted by $\rad(M)$, the \emph{Jacobson radical} of $M$. Write $P(M)$ for the projective cover of $M$, the smallest projective module that has $M$ as a quotient. If $I$ is a set of simple modules, the \emph{$I$-radical} of $M$ is the largest submodule whose composition factors lie in $I$, and the \emph{$I$-residual} is the smallest submodule for which the quotient only has composition factors lying in $I$. By $M^*$ we refer to the dual of $M$. The $I$-radical and $I$-residual are related by the following: the $I$-radical of $M^*$ is the dual of the quotient by the $I$-residual of $M$. In analogy with the notation of a $p'$-subgroup, $I'$-radical and $I'$-residual mean $(\Irr(H)\setminus I)$-radical and $(\Irr(H)\setminus I)$-residual.

Write $V_{\min}$ for one of the minimal modules for $G$, namely $L(\lambda_1)$ for $F_4$, either $L(\lambda_1)$ or $L(\lambda_6)$ for $E_6$ and ${}^2\!E_6$, $L(\lambda_1)$ for $E_7$ and $L(\lambda_1)$ for $E_8$. We write $L(G)$ for the simple, non-trivial constituent of the Lie algebra module, which is $L(\lambda_4)$, $L(\lambda_5)$, $L(\lambda_7)$ and $L(\lambda_1)$ respectively. These have the following dimensions:
\begin{center}
\begin{tabular}{ccc}
\hline Group & $\dim(V_{\min})$ & $\dim(L(G))$
\\\hline $F_4$ & $26-\delta_{q,3}$ & $52$
\\ $E_6$ & $27$ & $78-\delta_{q,3}$
\\ $E_7$ & $56$ & $133-\delta_{q,2}$
\\ $E_8$ & $248$ & $248$
\\ \hline
\end{tabular}
\end{center}
We will remind the reader when $L(G)$ does not have the usual dimension.

The actions of certain reductive subgroups, and all Levi subgroups, on these modules is helpfully tabulated in \cite{thomas2016}. There are many such subgroups and it is not necessary to reproduce the full list of tables, but we describe some of the more commonly used subgroups in Table \ref{tab:commonsubgroups}.

\begin{table}\begin{center}\begin{tabular}{cccc}
\hline Group & Subgroup & Factors on $\vmin$ & Factors on $L(G)$
\\ \hline $F_4$ & $B_4$ & 1000,0001,0000 & 0100,0001
\\ & $A_2\tilde A_2$ & (10,10),(01,01),(00,11) & (11,00),(00,11),(10,02),(01,20)
\\\hline $E_6$ & $F_4$ & 0001,0000 & 1000,0001
\\ & $A_1A_5$ & $(1,\lambda_1)$,$(0,\lambda_4)$ & $(2,0)$,$(0,\lambda_1+\lambda_5)$,$(1,\lambda_3)$
\\ \hline $E_7$ & $E_6$ & $\lambda_1$,$\lambda_6$,$0^2$ & $\lambda_2$,$\lambda_1$,$\lambda_6$,$0$
\\ & $A_7$ & $\lambda_2$,$\lambda_6$&$\lambda_1+\lambda_7$,$\lambda_4$
\\ & $D_6A_1$ & $(\lambda_1,1)$,$(\lambda_5,0)$ & $(\lambda_6,1)$,$(\lambda_2,0)$,$(0,2)$
\\ \hline $E_8$ & $A_1E_7$ & N/A & $(2,0)$,$(0,\lambda_1)$,$(1,\lambda_7)$
\\ & $D_8$ & N/A & $\lambda_2$, $\lambda_7$
\\ & $A_8$ & N/A & $\lambda_1+\lambda_8$, $\lambda_3$, $\lambda_5$

\\ \hline
\end{tabular}\end{center}
\caption{Actions of some common subgroups on $\vmin$ and $L(G)$}\label{tab:commonsubgroups}
\end{table}

The actions of unipotent elements on these modules are given in \cite{lawther1995}. From that we can extract much information. For example, we have Tables \ref{tab:unipe8p3} and \ref{tab:unipe8p5} of certain actions of elements of orders $3$ and $5$ on $L(G)$ for $E_8$. In these we only see the \emph{non-generic} classes: we should define what this means now. We use the Bala--Carter--Pommerening notation, in particular the precise notation used in \cite{lawther1995} which will be our main reference for unipotent actions, and this allows us to compare the actions of a unipotent class in different characteristics.

\begin{defn} Let $p$ be a prime, let $G$ be an algebraic group in characteristic $p$, and let $C$ be a unipotent class in $G$. Let $L(\lambda)$ be a highest-weight module for $G$. If the Jordan blocks of the action of $u\in C$ on $L(\lambda)$ are the same size as for the same class $C$ in $G$ for all primes $q>>0$, then $C$ is said to be \emph{generic} at the prime $p$. Otherwise, $C$ is \emph{non-generic}.
\end{defn}

As an example, as we see in \cite[Table 7]{lawther1995}, the class $E_7(a_2)$ in $E_7$ acts on $V_{\min}$ with blocks $18,16,10,8,4$ for all $p\neq 2,3,13,17$, but for example with blocks $13^4,4$ for $p=13$. This means that this class is generic for $V_{\min}$ whenever $p\neq 2,3,13,17$, and non-generic for these primes. Notice that, on $L(G)$, \cite[Table 8]{lawther1995} shows that $E_7(a_2)$ is non-generic for $p=2,3,5,13,17,19$, so being non-generic does depend on the module under consideration.

Our next two results mean that, if $H$ is a subgroup of an exceptional algebraic group $G$, and $H$ contains a unipotent element of order $p$ from a generic class, then there exists a positive-dimensional subgroup $X$ of $G$ containing $H$ and stabilizing the same subspaces of $V_{\min}$ as $H$. In particular, if $H$ does not act irreducibly on $V_{\min}$ then $H$ is contained in a known maximal subgroup of $G$, since the set of maximal positive-dimensional subgroups $\ms X$ of $G$ are known \cite{liebeckseitz2004}, and since $X\neq G$.

This result appears in \cite[Lemma 1.14]{lst1996}.

\begin{lem}\label{lem:unipotentreduction} Let $H$ be a subgroup of an algebraic group $\mathbb G$, and let $u\in H$ be a unipotent element of order $p$. If $u$ is contained in an $A_1$ subgroup that acts with $p$-restricted composition factors on a module $M$ for $G$, then there exists a positive-dimensional subgroup $\mathbb X$ of $\mathbb G$ such that $H\leq \mathbb X$, and $H$ and $\mathbb X$ stabilize the same subspaces of $M$. In particular, if $\mathbb G$ acts irreducibly on $M$ and $H$ acts reducibly on $M$, then $H$ is contained in a member of $\ms X$, and if $H$ is contained in $G=\mathbb G^\sigma$ then it is contained in a member of $\ms X^\sigma$.
\end{lem}

(To get from \cite[Lemma 1.14]{lst1996} to this result, let $X$ be the subgroup generated by $H$ and the unipotent subgroup $U$ containing $u\in H$ that is constructed in \cite[Lemma 1.14]{lst1996}.)

Now we need to know that such $A_1$s actually exist for the minimal and adjoint modules. Together with the previous lemma, this means that we can discount generic unipotent classes from consideration.

\begin{lem}\label{lem:genericgood}  Let $p$ be an odd prime and let $G$ be one of $F_4$, $E_6$, $E_7$ and $E_8$. Let $u$ be a unipotent element in $G$ of order $p$. If $u$ belongs to a generic class for either the minimal or adjoint module, then $u$ lies inside an $A_1$ subgroup whose composition factors on that module are $p$-restricted.
\end{lem}
\begin{pf}
In the case of the adjoint module, we simply use the tables in \cite{lawthertesterman1999}, which construct $A_1$ subgroups above unipotent classes, and whenever $u$ comes from a generic class the composition factors are $p$-restricted. It remains to consider the minimal module (except for $E_8$, where the two coincide).

We must construct, for each unipotent class, an algebraic $A_1$ with $p$-restricted composition factors, whenever $u$ is in the generic case. Using the embeddings $F_4<E_6<E_7$ we can reduce our work by embedding an $A_1$ inside $F_4$ into both $E_6$ and $E_7$.

We use the tables in \cite{lawther1995} both to determine for which primes we are in the generic case, and to give a list of the unipotent classes that we will use here. We also consult \cite{thomas2016} for the actions of Levi subgroups and irreducible $A_1$s.

Let $G=F_4$. For classes $A_1$ and $\tilde A_1$ we use the $A_1$ Levi subgroups, which act as needed. For $A_1+\tilde A_1$ we use the diagonally embedded $A_1$ inside the $A_1\tilde A_1$ Levi subgroup. For $A_2$ we use the irreducible $A_1$ inside the $A_2$ Levi subgroup.

For the rest of the classes we need $p\geq 5$. For $\tilde A_2$ we can use the irreducible $A_1$ inside the $\tilde A_2$ Levi subgroup. For $A_2+\tilde A_1$ we use the $A_1$ embedded in the Levi subgroup $A_2\tilde A_1$ acting as $L(2)$ on the $A_2$ and as $L(1)$ on the $A_1$ factor. The same subgroup of $\tilde A_2A_1$ deals with $u$ lying in $\tilde A_2+A_1$.

For $C_3(a_1)$, we consider the $C_3$ Levi subgroup, which acts on the minimal module as the sum of two copies of the natural module $L(100)$, and one copy of the adjoint $L(010)$, which is the exterior square of the natural minus a trivial. If $A_1$ is embedded as $L(1)\oplus L(3)$ then the exterior square of this is 
\[\Lambda^2(L(1))\oplus \Lambda^2(L(3))\oplus L(1)\otimes L(3)=L(0)\oplus (L(0)\oplus L(5))\oplus (L(2)\oplus L(4)),\]
yielding the correct action. (Remember that $\Lambda^2(A\oplus B)=\Lambda^2(A)\oplus \Lambda^2(B)\oplus A\otimes B$.)

For $u$ lying in $F_4(a_3)$, \cite[Table 2]{lawthertesterman1999} suggests we can look inside $C_3A_1$: this acts as the sum of $010\otimes 0$ (the exterior square of the natural minus a trivial, tensored by the trivial) and $100\otimes 1$ (the natural tensored by the natural). Inside here, we take a diagonal $A_1$ which projects non-trivially along the $C_3$ and acts as $L(1)\oplus L(3)$ on the natural, as in the previous case of $C_3(a_1)$, but also as $L(1)$ on the $A_1$ factor. This then acts as $L(1)^{\otimes 2}\oplus L(3)\otimes L(1)=L(0)\oplus L(2)^{\oplus 2}\oplus L(4)$ on $100\otimes 1$ and as $\Lambda^2(L(1)\oplus L(3))=L(0)^{\oplus 2}\oplus L(2)\oplus L(4)^{\oplus 2}$ minus a copy of $L(0)$ on $010\otimes 0$, as needed.

For $B_3$ we need $p\geq 7$. Inside the $B_3$ Levi subgroup we take an $A_1$ acting as $L(6)$, i.e., irreducibly, on the natural, and hence as $L(0)\oplus L(6)$ on the $8$-dimensional module, giving the correct action.

We now need $p\geq 11$. For $C_3$ we take the $A_1$ lying in the $C_3$ Levi subgroup acting as $L(5)$: since $\Lambda^2(L(5))=L(0)\oplus L(4)\oplus L(8)$ we get the correct action on the minimal module.

For $F_4(a_2)$, \cite[Table 2]{lawthertesterman1999} suggests that as for $F_4(a_3)$ we can look inside $C_3A_1$: the diagonal $A_1$ acting as $L(5)$ on the natural for $C_3$ and the natural on the $A_1$ factor works, because the $100\otimes 1$ restricts to $L(4)\oplus L(6)$, and the $010\otimes 0$ restricts to $L(4)\oplus L(8)$, as needed.

For $F_4(a_1)$ and $F_4$, we see in \cite[Table 10]{thomas2016} that there are two irreducible $A_1$s, namely subgroups 7 and 10, that cover these two classes, completing the proof for $F_4$.

\bigskip

For $E_6$ we can take the $A_1$s from $F_4$ and add a trivial, dealing with many classes. We run through those that are left.

The first such class is $A_4$, which requires $p\geq 7$ from now on. The $A_4$ Levi subgroup acts as (up to duality, which is not important for $A_1$s) three copies of the natural, one of its exterior square, and two trivials. If we embed as $A_1$ as $L(4)$ then $\Lambda^2(L(4))=L(2)\oplus L(6)$, and so we get $L(6),L(4)^3,L(2),L(0)^2$ as the factors of the $A_1$ on the minimal module, as needed.

The next class is $A_4+A_1$, for which we embed an $A_1$ diagonally into the $A_4A_1$ Levi as $L(4)\otimes L(1)$. This acts on $1000\otimes 1$ as $L(5)\oplus L(3)$, and on $0001\otimes 0$ as $L(4)$. It acts on $0010\otimes 0$ as $L(6)\oplus L(2)$, as in the previous case, and so we get the correct factors for covering an element from class $A_4+A_1$.

There are only two classes left: $D_5(a_1)$ needs $p\geq 11$, and this $A_1$ can be found inside $D_5$, acting as $L(7)\oplus L(5)\oplus L(1)$ on the $16$-dimensional spin module and $L(6)\oplus L(2)$ on the $10$-dimensional natural module. The last one is $E_6(a_1)$, which needs $p\geq 13$, and is covered by subgroup 6 from \cite[Table 11]{thomas2016}. This completes the proof for $G=E_6$.

\bigskip

As with $E_6$ and $F_4$, we can use $E_6$ to exclude many classes for $E_7$, and deal with only those that are left. However, this still leaves twenty-four classes, and so we want to cut a few of those down before constructing individual groups.

We use the subgroup $F_4A_1$ for this, which acts with factors $0000\otimes 3$ and $1000\otimes 1$. Notice that every generic class of order $p$ in $F_4$ for the module $1000$ is covered by an $A_1$ acting with $p$-restricted factors, as we have just proved it, and obviously the same is true for the $A_1$ factor; therefore at least for $V_{\min}$ we know that every unipotent class of $F_4A_1$ that is generic is covered by an $A_1$ subgroup acting with $p$-restricted composition factors also. (If one is worried that there is a composition factor $L(p-1)\otimes L(1)$ in $V_{\min}$, note that the class is not generic for $p$, and only for primes larger than $p$.)

Consulting the first and third columns of Table 38 in \cite{lawther2009}, we get the unipotent classes contained in $F_4A_1$: the third column contains the classes $(3A_1)''$, $4A_1$, $A_2+3A_1$, $(A_3+A_1)''$, $A_3+2A_1$, $D_4(a_1)+A_1$, $A_3+A_2$, $A_3+A_2+A_1$, $D_5(a_1)+A_1$, $D_6(a_2)$, $E_7(a_5)$, $E_7(a_4)$ and $E_7(a_2)$.

This leaves only eleven unipotent classes to cover with $A_1$s with $p$-restricted composition factors, which is enough that we can construct the subgroups explicitly. The list of classes left to cover is
\begin{center}
$(A_5)''$, $D_4+A_1$, $A_4+A_2$, $A_5+A_1$, $D_6(a_2)$, $A_6$, $D_5+A_1$, $D_6$, $E_7(a_3)$, $E_7(a_2)$, $E_7(a_1)$, $E_7$.
\end{center}

In the $A_5$ Levi subgroup one sees $(A_5)''$ as the regular class, and the regular $A_1$ covering it acts as $L(3)\oplus L(5)\oplus L(9)$ on the exterior cube of the natural, hence has the right properties for $p\geq 11$. Taking this class and tensoring it by an $A_1$ factor inside the $A_5A_2$ Levi subgroup provides a module acting with composition factors of dimension $6$ (from $L(\lambda_1)\otimes L(0)$ and its dual), $10,6,4$ (from $L(\lambda_3)\otimes L(0)$) and $7,5$ (from $L(\lambda_1)\otimes L(1)$ and its dual), hence covers $A_5+A_1$.

Inside $A_7$ a multiplicity-free module for $\PGL_2$ will have a fixed-point-free exterior square, and so we can take $A_1$s acting as $L(0)\oplus L(6)$ and $L(2)\oplus L(4)$, and these act on $V_{\min}$ as $L(10)\oplus L(6)^{\oplus 2}\oplus L(2)$ and $L(6)^{\oplus 2}\oplus L(4)\oplus L(2)^{\oplus 3}$, containing representatives from classes $A_6$ and $A_4+A_2$.

For $E_7(a_1)$ and $E_7$ we need $p\geq 23$ and $p\geq 29$ respectively, and the maximal $A_1$s inside $E_7$ provide the subgroups.

\medskip

For the rest we need to examine the $D_6A_1$ maximal-rank subgroup. Consider firstly a copy of $\SL_2$ inside $D_6$ that acts on the natural as $L(10)\oplus L(0)$. This acts on the spin module as $L(15)\oplus L(9)\oplus L(5)$, at least as long as the characteristic is at least $17$. This will cover the unipotent class $D_6$.

To cover $E_7(a_3)$ we simply take the diagonal subgroup between this $A_1$ and the $A_1$ factor. This leaves the $32$ alone but tensors the $L(10)\oplus L(0)$ by $L(1)$ to get $L(11)\oplus L(9)\oplus L(1)$, finishing this case.

This also shows how to construct $A_1$s above $D_4+A_1$ and $D_5+A_1$: we take the $A_1$s lying above the $D_4$ and $D_5$ classes, which are known to exist and act as $p$-restricted composition factors whenever we are in the generic case, since they lie inside the $E_6$ Levi subgroup, and as with $F_4A_1$ we tensor with the $A_1$ in the second factor of $D_6A_1$: as this $A_1$ acts as $L(8)\oplus L(0)^{\oplus 3}$ on the natural, we get $p$-restricted composition factors, as needed. Similarly, the $A_1$ covering the $D_4$ class acts as $L(6)\oplus L(0)^{\oplus 5}$ on the natural, and so we get the right composition factors for the diagonal subgroup with the $A_1$ factor.

The last unipotent class left to cover is $D_6(a_1)$. Now we consider an $\SL_2$ acting on the natural as $L(8)\oplus L(2)$, which lifts to the spin module with factors $L(11)$, $L(9)$, $L(5)$ and $L(3)$, at least for $p\geq 13$. This yields the final class of $A_1$s needed to prove the result.
\end{pf}

\begin{table}
\begin{center}
\begin{tabular}{cc}
\hline Class & Action on $L(G)$
\\\hline $3A_1$ & $3^{31},2^{50},1^{55}$
\\ $A_2$ & $3^{57},1^{77}$
\\ $4A_1$ & $3^{44},2^{40},1^{36}$
\\ $A_2+A_1$ & $3^{58},2^{20},1^{54}$
\\ $A_2+2A_1$ & $3^{65},2^{16},1^{21}$
\\ $A_2+3A_1$ & $3^{70},2^{14},1^{10}$
\\ $2A_2$ & $3^{78},1^{14}$
\\ $2A_2+A_1$ & $3^{79},2^2,1^7$
\\ $2A_2+2A_1$ & $3^{80},2^4$
\\ \hline
\end{tabular}\caption{Non-generic unipotent classes for $E_8$ and order $3$}\label{tab:unipe8p3}
\end{center}
\end{table}

\begin{table}
\begin{center}
\begin{tabular}{cc}
\hline Class & Action on $L(G)$
\\\hline $A_2$ & $4^{2},3^{54},1^{78}$ 
\\$A_2+A_1$ & $4^{40},3^{12},2^{18},1^{16}$
\\$A_2+2A_1$ & $4^{14},3^{30},2^{34},1^{34}$
\\$A_2+3A_1$ & $4^{26},3^{6},2^{62},1^{2}$
\\$2A_2$ & $4^{28},3^{36},1^{28}$
\\$2A_2+A_1$ & $4^{40},3^{12},2^{18},1^{16}$
\\$2A_2+2A_1$ & $4^{44},3^4,2^{28},1^4$
\\$A_3$ & $4^{46},2^{10},1^{44}$
\\$A_3+A_1$ & $4^{46},2^{24},1^{16}$
\\$A_3+2A_1$ & $4^{46},2^{30},1^4$
\\$A_3+A_2$ & $4^{50},3^{10},2^{6},1^6$
\\$D_4(a_1)$ & $4^{54},3^2,1^{26}$
\\$D_4(a_1)+A_1$ & $4^{54},3^2,2^{10},1^6$
\\$A_3+A_2^{(2)}$ & $4^{54},3^2,2^{10},1^6$
\\$A_3+A_2+A_1$ & $4^{54},3^2,2^{12},1^2$
\\$D_4(a_1)+A_2$ & $4^{56},3^8$
\\$2A_3$ & $4^{60},2^4$
\\ \hline
\end{tabular}\caption{Unipotent classes for $E_8$ and order $4$}\label{tab:unipe8p4}
\end{center}
\end{table}

\begin{table}
\begin{center}
\begin{tabular}{cc}
\hline Class & Action on $L(G)$
\\\hline $A_3$ & $5^{13},4^{32},1^{55}$
\\ $2A_2+A_1$ & $5^{14},4^{14},3^{23},2^{18},1^{17}$
\\ $2A_2+2A_1$ & $5^{18},4^{12},3^{20},2^{20},1^{10}$
\\ $A_3+A_1$ & $5^{21},4^{16},3^9,2^{14},1^{24}$
\\ $A_3+2A_1$ & $5^{25},4^{10},3^{14},2^{14},1^{13}$
\\ $D_4(a_1)$ & $5^{29},3^{25},1^{28}$
\\ $D_4(a_1)+A_1$ & $5^{29},4^6,3^{14},2^{14},1^9$
\\ $A_3+A_2$ & $5^{30},4^8,3^{13},2^8,1^{11}$
\\ $A_3+A_2+A_1$ & $5^{32},4^8,3^{10},2^{10},1^6$
\\ $D_4(a_1)+A_2$ & $5^{36},3^{20},1^8$
\\ $2A_3$ & $5^{38},4^{12},1^{10}$
\\ $A_4$ & $5^{45},1^{23}$
\\ $A_4+A_1$ & $5^{45},3,2^6,1^8$
\\ $A_4+2A_1$ & $5^{45},3^4,2^4,1^3$
\\ $A_4+A_2$ & $5^{46},3^5,1^3$
\\ $A_4+A_2+A_1$ & $5^{46},4^2,3^2,2^2$
\\ $A_4+A_3$ & $5^{48},4^2$
\\ \hline
\end{tabular}\caption{Non-generic unipotent classes for $E_8$ and order $5$}\label{tab:unipe8p5}
\end{center}
\end{table}

\begin{lem}\label{lem:smallspacestabilizers} Suppose that $H$ is a subgroup of $G=G(q)$, and suppose that
\begin{itemize}
\item $H$ stabilizes a $1$-space on either $V_{\min}$ or $L(G)$,
\item $H$ stabilizes a $2$-space on $V_{\min}$ for $G$ of type $F_4$, $E_6$ or $E_7$, or
\item $H$ stabilizes a $3$-space on $V_{\min}$ for $G$ of type $E_6$.
\end{itemize} Then $H$ is contained in a member of $\ms X^\sigma$.
\end{lem}
\begin{pf} If $H$ fixes a $1$-space on $L(G)$ then $H$ is contained in either a parabolic or semisimple subgroup by Seitz \cite[(1.3)]{seitz1991}, and thus is in a member of $\ms X^\sigma$. If $H$ stabilizes a line on $V_{\min}$ then $H$ is contained in a member of $\ms X^\sigma$, for example see \cite[Lemma 2.2]{liebeckseitz2005}. (For $G=F_4$, $H$ is contained in $B_4$ or a maximal parabolic, for $G=E_6$ $H$ is contained in $F_4$ or a $D_5$ parabolic, and for $G=E_7$ $H$ is contained inside an $E_6$ or $D_6$ parabolic, or inside an $E_6$ or ${}^2\!E_6$ subgroup with an automorphism on top.)

For the others, we simply compute the dimension of the subspace stabilizer inside the algebraic group, and note that it is non-zero. For example, the dimension of the stabilizer of a $3$-space in $\vmin$ for $G=E_6$ is at least
\[ 78-(27+26+25-3-2-1)=6>0.\]
Having done this, the stabilizer of a subspace $W$ of the appropriate dimension above is therefore positive dimensional. To see $\sigma$-stability, we note that $W$ is fixed by $H\leq G$, so $W$ is defined over $\F_q$ and in $\vmin^\sigma$.
\end{pf}

The proves non-maximality of simple subgroups of simple groups, but when we move to almost simple exceptional groups we need to be more careful, as $\Vmin$ is not always stable under the outer automorphisms.

The next lemma deals with so-called novelty maximal subgroups. Suppose that $G$ is an almost simple group with socle $X$, and let $M$ be a maximal subgroup of $G$ not containing $X$. If $M\cap X$ is a maximal subgroup of $X$ then $M$ is said to be an \emph{ordinary} maximal subgroup, and if $X\cap M$ is not maximal in $X$ then $M$ is said to be a \emph{novelty} maximal subgroup. Notice that, while ordinary maximal subgroups can be easily found from a list of maximal subgroups of $X$, novelties are harder, and arise when a subgroup $H$ is stabilized by a group of outer automorphisms of $X$ (hence extending to a larger subgroup of $G$) but any subgroup of $X$ containing $H$ is not stabilized by the outer automorphisms, so it becomes maximal.

The alternating subgroups of exceptional groups considered here never become novelty maximal subgroups.

\begin{lem}\label{lem:novelty} Let $G$ be an almost simple group with socle an exceptional group of Lie type. Suppose that $H$ is an almost simple subgroup of $G$ such that $F^*(H)$ has abelian outer automorphism group.
\begin{enumerate}
\item If the preimage of $F^*(H)$ in the simply connected form of $F^*(G)$ fixes a line on either the minimal or adjoint module for $F^*(G)$ then $H$ is contained in a member of $\ms X^\sigma$. In particular, $H$ is not a novelty maximal subgroup of $G$.
\item If there exists a positive-dimensional subgroup $X$ of $F^*(G)$ containing $F^*(H)$ and stabilizing the same subspaces of a simple module $V$ for $F^*(G)$ that is $G$-stable and such that $H$ acts reducibly on $V$, then $H$ is not a novelty maximal subgroup of $G$.
\item Suppose that $G$ induces a graph automorphism $\theta$ on $F^*(G)$, and let $\bar H$ denote the preimage of $F^*(H)$ in the simply connected form $\bar G$ of $F^*(G)$. If $W$ is an $\bar H$-submodule of the minimal module $\Vmin$ for $\bar G$ such that $W^\phi$ is not an $\bar H$-submodule of $\Vmin^\theta$ for any $\phi\in\Out(\bar H)$, then the $G$-conjugacy class of $F^*(H)$ splits upon restriction to $F^*(G)$. In particular, $H$ is not maximal in $G$.
\end{enumerate}
\end{lem}
\begin{pf}
The first part follows from the orbit-stabilizer theorem: if $G$ is an almost simple group that does not induce a graph automorphism on $F^*(G)$ then the automorphisms of $G/F^*(G)$ stabilize both the minimal and adjoint modules, and so the dimension of a line stabilizer is positive.

If $G$ does induce a graph automorphism then $L(G)$ is still stabilized, so that is fine, but $\Vmin$ is swapped with another module of the same dimension: however, $\Vmin\oplus \Vmin^\theta$ has dimension at most that of $G$ (in the case $F_4$ and $p=2$) and so the line stabilizer on that module is still at least $1$-dimensional. As $F^*(H)$ fixes a line, and $\Out(F^*(H))$ is abelian, $H$ still fixes a line on $\Vmin\oplus \Vmin^\theta$, so lies inside a positive-dimensional subgroup of $G$, hence $H$ lies inside a member of $\ms X^\sigma$.

(Notice that if $H\leq G^\sigma$ then while not every line is $\sigma$-stable, there is a $\sigma$-stable line, so $H$ is contained in a member of $\ms X^\sigma$.)

The second part is easy: the subgroup $\langle X,H\rangle$ is proper in $G$ since it acts reducibly on $V$ and hence $H$ cannot be a novelty maximal subgroup.

For the third part, we note that, in order for $F^*(H)$ to be stabilized by $\theta$, we must have an automorphism $\phi\in\Out(F^*(H))$ such that $\Vmin^\theta\downarrow_H=\Vmin\downarrow_H^\phi$. By assumption this is not true, so $\theta$ does not normalize $F^*(H)$ and $N_G(F^*(H))=N_{G'}(F^*(H))$ for some proper subgroup $G'$ of $G$. In particular, $H\leq G'$, and so is not maximal.
\end{pf}

In order to examine embeddings of alternating groups into exceptional groups, we need to know the potential composition factors of the restrictions of $V_{\min}$ and $L(G)$ to the subgroup. One way to do this is to know the traces of semisimple elements of the algebraic group, i.e., the possible Brauer characters of the restriction. We end up using traces of elements of order up to $15$ (as there are elements of order $15$ in $\Alt(8)$) but it would be impractical to list these here. We simply list the integral traces of elements of orders $2$, $3$ and $5$, since this is enough for many purposes. The others can be downloaded from the author's webpage, or computed with Litterick's program in \cite{litterick}.

\begin{lem}\label{lem:traces} The traces of semisimple elements of orders $2$, $3$ and $5$, whose entries are integers, on $V_{\min}$ and $L(G)$ respectively, are as follows.
\begin{center}
\begin{tabular}{cccc}
\hline Group & Order & Trace on $V_{\min}$ & Trace on $L(G)$
\\\hline $F_4$ & $2$ & $2$, $-6$ & $-4$, $20$
\\ & $3$ & $-1$, $8$ & $-2$ or $7$, $7$
\\ & $5$ & $1$ & $2$
\\\hline $E_6$ & $2$ & $3$, $-5$ & $-2$, $14$
\\ & $3$ & $0$, $9$ & $-3$ or $6$, $15$
\\ & $5$ & $2$ & $3$
\\\hline $E_7$ & $2$ & $8$, $-8$ & $5$, $5$
\\ & $3$ & $-25$, $-7$, $2$, $20$ & $52$, $7$, $-2$ or $7$, $34$
\\ & $5$ & $6$ & $8$
\\ \hline $E_8$ & $2$ & N/A & $-8$, $24$
\\ & $3$ & N/A & $-4$, $5$, $14$, $77$
\\ & $5$ & N/A & $-2$, $3$, $23$
\\ \hline
\end{tabular}
\end{center}
\end{lem}

Since we are aiming to prove that a module has a trivial submodule or quotient, we need a result that will guarantee such an outcome. A slightly weaker version of this has been considered before, for example \cite[Lemma 1.2]{lst1996} or \cite[Proposition 3.6]{litterick}. We start with a definition.

\begin{defn}
Let $G$ be a finite group, $k$ be a field, and let $M$ be a $kG$-module. Suppose that, for every composition factor $V$ of $M$, $H^1(G,V)=H^1(G,V^*)$. (This is true, for example, if there exists an automorphism swapping simple modules and their duals.) Suppose further that $H^1(G,k)=0$, i.e., $O^p(G)=G$. Write $\cf(M)$ for the multiset of composition factors of $M$. The \emph{pressure} of $M$ is the quantity
\[ \sum_{V\in \cf(M)} (H^1(G,V)-\delta_{V,k}).\]
\end{defn}

The pressure of a module is an indicator of how easy it is to hide trivial composition factors in the middle of it. Before we prove that we need the following lemma, which gives some indication that the pressure of a module is an interesting invariant.

\begin{lem}\label{lem:negativepressure}
Let $G=O^p(G)$ be a finite group, $k$ be a field, and let $M$ be a $kG$-module. Suppose that, for every composition factor $V$ of $M$, $H^1(G,V)=H^1(G,V^*)$.
\begin{enumerate}
\item If $M$ is a module with negative pressure, then $M$ has a trivial submodule and a trivial quotient.

\item Suppose that $M$ contains trivial composition factors, but no trivial submodules. If $M$ has pressure $n$ then $\dim(H^1(G,M))\leq n$.
\end{enumerate}
\end{lem}
\begin{pf}(i) is known, for example \cite[Proposition 3.6]{litterick}.

For (ii), let $M'$ be the module with no trivial submodules obtained from $M$ by extending by $k^m$, where $m=\dim(H^1(G,M))$, and suppose that $M$ has pressure less than $m$, so that $M'$ has negative pressure. The module $(M')^*$ still has negative pressure but no trivial quotient, a contradiction.
\end{pf}

Even better than modules of negative pressure, we have the following result on pressures, which helps in eliminating a variety of cases from consideration.

\begin{prop}\label{prop:lowpressuremodules}
Let $G=O^p(G)$ be a finite group, $k$ be a field, and let $M$ be a $kG$-module. Suppose that, for every composition factor $V$ of $M$, $H^1(G,V)=H^1(G,V^*)$. Suppose that $M$ contains trivial composition factors, but no trivial submodule or quotient.

If $M$ has pressure $n$ then there is no subquotient $A$ of $M$ with pressure greater than $n$ or less than $-n$. In particular, if $M$ has a trivial composition factor but has non-positive pressure, then $M$ has a trivial submodule or quotient, and if $M$ has a composition factor $V$ whose $1$-cohomology has dimension greater than the pressure of $M$, then $M$ has a trivial submodule or quotient.
\end{prop}
\begin{pf}
Let $A$ be a subquotient of $M$, and let $B$ and $C$ be the other two factors involved in $M$, i.e., $B$ is a submodule of $M$ and $C$ a quotient of $M$ whose kernel contains $B$ and such that the surjective map $M/B\to C$ has kernel isomorphic to $A$. If $A$ has pressure greater than $n$ then at least one of $B$ and $C$ must have negative pressure, without loss of generality (take duals) $B$. Since $B$ has negative pressure, it has a trivial submodule, a contradiction. Similarly, if $A$ has pressure less than $-n$, then one of $B$ and $C$ has pressure at most $n$, say $B$ again. Since the pressure of $C$ is less than $-n$, it has a trivial submodule of dimension at least $n+1$ by repeated application of Lemma \ref{lem:negativepressure}. However, $H^1(G,B)$ has dimension at most $n$, again by Lemma \ref{lem:negativepressure}, and so we cannot fit this trivial submodule from $C$ on top of $B$ without getting a trivial submodule in $M$.
\end{pf}

Modules $M$ of pressure $1$ are fairly easy to characterize: ignoring non-trivial composition factors with no $1$-cohomology, the structure of $M$ is uniserial, alternating between non-trivial and trivial modules. (Of course, $M$ need not be uniserial, but ignoring non-trivial composition factors with no $1$-cohomology it is. What this means is that there is no subquotient $k\oplus k$ and no subquotient $V\oplus W$ with $H^1(G,V),H^1(G,W)>0$.)

The reason for the name pressure is that modules of low pressure will have socle structures that look long and thin (close to uniserial) and as the pressure increases they tend to get squashed down.

\bigskip

Suppose that $H$ is a subgroup of a finite group $G$, and let $M$ be a $kG$-module. Frobenius reciprocity gives us information about how to relate $M$ and $M\downarrow_H$. Recall that Frobenius reciprocity states that if $N$ is a $kH$-module, then
\[ \Hom(N,M\downarrow_H)=\Hom(N\uparrow^G,M).\]
The next result follows immediately from this.

\begin{lem}\label{lem:frobrep}
Let $G$ be a finite group and let $H$ be a subgroup of $G$. Let $k$ be a field and let $M$ be a $kG$-module. Write $N$ for the permutation module of $G$ on the cosets of $H$. If $\alpha\in M$ is fixed by $H$ then there is a map $\phi:N\to M$ whose image contains $\alpha$.
\end{lem}

The next result finds singular subspaces of orthogonal and symplectic groups preserved by a particular subgroup. It comes in handy when analysing copies of $\Alt(5)$ in the $D_6$ parabolic of $E_7$ in characteristic $3$, at the end of Section \ref{sec:alt5}.

\begin{lem}\label{lem:singularsubspace} Let $\F$ be a field, let $H$ be a finite group, and let $H$ act on an $\F$-vector space $V$ of dimension $n$ preserving a bilinear form $f$. If $-1$ is a sum of $m$ squares in $\F$ then $H$ stabilizes a totally isotropic $n$-dimensional subspace of $V^{\oplus (m+1)}$.

Consequently, if $H$ is a subgroup of an orthogonal or symplectic group in positive characteristic, and the action of $H$ on the natural module has a submodule isomorphic with three copies of a given module of dimension $n$, then $H$ stabilizes a totally isotropic subspace of the natural module of dimension $n$, and hence lies inside a parabolic subgroup of the orthogonal or symplectic group.
\end{lem}
\begin{pf}
If $\lambda_1,\dots,\lambda_m\in \F$, then let $H$ act on the diagonal subspace $\{(v,\lambda_1 v,\dots,\lambda_m v):v\in V\}$: the restriction of $f$ to this subspace has norm $(1+\lambda_1^2+\cdots +\lambda_m^2)(v,w)$, for $v,w\in V$. If $\lambda_1^2+\cdots+\lambda_m^2=-1$ then this means the form vanishes, and we get the result. 

If the characteristic of $\F$ is $p$ then we claim that $-1$ is always the sum of two squares in $\F_p$: to see this, note that if non-squares are not the sum of two squares, then the squares form a subfield of $\F_p$.

The consequence arises because the stabilizers of totally isotropic subspaces are parabolic subgroups.
\end{pf}

We often know that a subgroup $H$ of $G$ is contained in a member of $\ms X^\sigma$, and we want to find all possible embeddings. We of course can run through all members of $\ms X^\sigma$, but some of them, like the normalizer of a torus, are maximal while the quotient by the connected component of the identity could contain a copy of our simple group, i.e., $X/X^0$ is not soluble. The next result states that, if an alternating group $H$ is contained inside a member of $\ms X^\sigma$ then $H$ is contained inside a \emph{connected} parabolic or reductive subgroup, rather than inside a subgroup like the normalizer of a torus, where $H$ lies inside the Weyl group. In this lemma we leave out the case where $G=E_7$, and $H=\PSU_3(3)$ in all characteristics and $\SL_2(8)$ in characteristic $2$, but this does not mean that these are genuine counterexamples, merely that our methods do not immediately cover it, and the author's later results do not need these specific instances of the lemma.

\begin{lem}\label{lem:insideconnected}
Let $H$ be a simple subgroup of a member $X$ of $\ms X$. Either there exists some member $Y$ of $\ms X$ such that $H\leq Y^0$, or $G=E_7$ and $H=\PSU_3(3)$, or $G=E_7$, $p=2$, and $H=\SL_2(8)$.
\end{lem}
\begin{pf}
In \cite{liebeckseitz2004} all maximal subgroups of positive dimension are found, and so it is easy to see which members of $\ms X$ have $X/X^0$ insoluble. If $G=G_2$ or $G=F_4$ then all subgroups $X/X^0$, which are subgroups of the Weyl group, are soluble, so the result holds.

If $G=E_6$ then the only maximal subgroup $X$ with $X/X^0$ insoluble is the normalizer of the torus, and as the Weyl group has $\PSU_4(2)$ of index $2$, we have that $H$ is $\PSU_4(2)$. If $p=2$ then $H$ is a subgroup of $N_G(T)$, and an easy calculation shows that $H$ acts on $\vmin$ as $1\oplus 6/14/6$, so $H$ is contained in a line stabilizer. Outside of characteristic $2$, $N_G(T)$ is a non-split extension and there is no such subgroup $H$ of it. (This is easy to check by computer for some prime $p$, and the extension to all primes is formal.)

\medskip

For $G=E_7$, we can have that $X=N_G(T)$ or $X=A_1^7.\PSL_2(7)$. In the first case, the Weyl group is $2\times \Sp_6(2)$, by \cite[Lemma 3.10]{litterick} we reduce to the case where $H=\PSL_2(8),\PSU_3(3),\Sp_6(2)$.

If $H=\Sp_6(2)$, then again there is copy of $H$ inside $N_G(T)$ unless $p=2$, in which case $H$ is a subgroup of $N_G(T)$, and acts on $\vmin$ as
\[ 1/6,1,14/6,6/1,14/6/1,\]
so it fixes a line on $\vmin$ and therefore lies inside an $E_6$ parabolic or $D_6$ parabolic, as needed. (Since the $1$-cohomology of the reflection representation is $0$ for $p=3,5,7$ by an easy computer calculation, there is only one class of complements to the torus in $N_G(T)$.) The same is not true for $\PSU_3(3)$ and $\SL_2(8)$ though, which have $1$-dimensional $1$-cohomology on this module in characteristic $3$.

If $H=\SL_2(8)$ then the non-split extension of $T$ by $\Sp_6(2)$ restricts to a non-split extension of $T$ by $H$, so $H$ is not a subgroup of $N_G(T)$ for $p$ odd.

We leave $\PSU_3(3)$, and consider the subgroup $X=A_1^7.\PSL_2(7)$, which can be thought of as $A_1\wr \PSL_2(7)$ (with $\PSL_2(7)$ acting on seven points). The conjugacy classes of complements in a permutational wreath product is in bijection with the homomorphisms from the point stabilizer to the base group, in this case maps from $\Sym(4)$ to $\SL_2(q)$. Of course, there are exactly two of these, and both maps restrict to maps from $\Sym(4)$ to a maximal torus of $\SL_2(q)$, say $T_1$. By the description in \cite{houghton1975} -- but see \cite[p.208]{hassanabadi1978} for a clearer description -- the two conjugacy classes of complements in $A_1\wr \PSL_2(7)$ intersect the two classes of complements in $T_1\wr \PSL_2(7)$, and hence if $H=\PSL_2(7)$ then $H$ can be chosen to lie in $N_G(T)$, which has already been dealt with. This completes the proof for $E_7$.

\medskip

Finally, let $G=E_8$. If $H$ is contained in $N_G(T)$ then we are done by \cite[Lemma 3.10]{litterick}, so we may assume that $H$ is contained in $A_1\wr \AGL_3(2)$. There are two classes of simple subgroup of this, both $\PSL_2(7)$ of course, one transitive and one fixing a point. The one fixing a point centralizes an $A_1$ subgroup and lies in the reductive subgroup $A_1E_7$. The transitive one has point stabilizer a Frobenius group of order $21$, so again there are exactly two homomorphisms from this to $A_1$, both of which restrict to homomorphisms to $T_1$, and so we again lie inside $N_G(T)$, so inside $X^0$ for some $X\in \ms X$. This completes the proof.
\end{pf}

We end this section with one obvious, but useful, result, which can be found in, for example \cite[Lemma 4.12]{freyun}. It is simply because $H$ lies inside the centralizer of an involution, and the two classes of involutions in $E_8$ in odd characteristic have centralizers $A_1E_7$ and $D_8$.

\begin{lem}\label{lem:centinv}
Let $H$ be a perfect subgroup of $G=E_8$ in characteristic $p\neq 2$, and suppose that $C_G(H)$ contains an involution. Either $H\leq D_8$ or $H\leq A_1E_7$. If $H$ has no representation of dimension $2$ over $\bar{\F}_p$ then $H$ is contained in either $E_7$ or $D_8$.
\end{lem}

\section{Blocks with cyclic defect group}

In this section we will summarize the theory of blocks with cyclic defect groups, at least as much as is needed for our purposes. More details can be found in, for example \cite{alperin} and \cite{feit}. This includes the case where the Sylow $p$-subgroup of a finite group $G$ is cyclic.

Throughout this section, $G$ is a finite group and $k$ is an algebraically closed field of characteristic $p$, $B$ is a block of $kG$ with cyclic defect group $D$. The number of simple $B$-modules is $e$, and $(|D|-1)/e$, an integer, is the \emph{exceptionality} of $B$. If $D$ is a Sylow $p$-subgroup then $e=|N_G(D)|/|C_G(D)|$. \textbf{We will assume that $e>1$}.

To $B$ we may associate a planar-embedded tree, the \emph{Brauer tree} of $B$ (technically it is merely a ribbon graph, but a planar-embedded tree is easier to envisage), whose edges are labelled by the simple $B$-modules, and whose vertices are labelled by the ordinary characters of $B$. If the exceptionality of $B$ is $1$, then any two different ordinary irreducible characters in $B$ have different $p$-modular reductions, but if the exceptionality is $n>1$ then there is a unique $p$-modular reduction that occurs $n$ times amongst all $p$-modular reductions of irreducible characters, with all others occurring with multiplicity $1$. The $n$ vertices of the Brauer tree that correspond to characters with the same $p$-modular reduction are identified and referred to as the \emph{exceptional vertex}.

If $\chi$ is an irreducible ordinary character of $B$ then the Brauer characters that are constiuents of the $p$-modular reduction of $\chi$ occur with multiplicity $1$. The edges incident to a given vertex of the Brauer tree label the composition factors of the $p$-modular reduction of the corresponding ordinary character.

We need to describe the embedding of the tree now. We order the vertices so that $M$ and $N$ appear next to each other in clockwise order around some mutually incident vertex if and only if $\Ext^1(M,N)=k$, and in all other cases $\Ext^1(M,N)=0$.

\medskip

In our case of alternating groups all Brauer trees are lines, and so the issue of ordering of the vertices is not important. However, for arbitrary simple groups in non-defining characteristic the Brauer tree can be very complicated, and it is for this reason that we describe the general situation.

\medskip

Having described the Brauer tree, we now explain how to use the tree to construct the projective cover of each simple module in $B$. If $S$ is a simple $B$-module, corresponding to an edge in the tree with incident vertices labelled by $\chi$ and $\phi$, then $\chi+\phi$ is the projective character of $S$. Furthermore, the projective cover $P(S)$ has the following structure:

The top and socle of $P(S)$ are both $S$. Removing both of these, the module splits as the direct sum of two (possibly one is the zero module) uniserial modules $M$ and $N$, corresponding to the two incident vertices $\chi$ and $\phi$ respectively.

Write $n$ for the exceptionality of $\chi$ if it is exceptional, and $n=1$ otherwise. Starting from $S$, write the edges appearing in clockwise order in $n$ complete revolutions around the vertex $\chi$. For example, if there are four edges, $S=S_1$, $S_2$, $S_3$ and $S_4$ around $\chi$, and $n=2$, we get the list
\[ S_1,S_2,S_3,S_4,S_1,S_2,S_3,S_4.\]
Delete the first copy of $S$ from the start: these are the radical layers of the uniserial module $M$. The same process around $\phi$ produces $N$, and $P(S)$ has structure
\[ S/(M\oplus N)/S.\]

Notice that if the Brauer tree is a star (i.e., all but one vertex has valency $1$) and the central vertex is exceptional (if there is such a vertex), then all projective modules are uniserial. If $G$ is $p$-soluble then the Brauer tree of $B$ is always a star with exceptional node in the centre.

\medskip

The quotients of the projective modules give many indecomposable modules in $B$, and in general, using the Green correspondence and the structure of the Brauer tree of blocks of $p$-soluble groups, it is easy to see that there are $(|D|-1)e$ non-isomorphic indecomposable modules in $B$. Hence it is possible to list them all, and in particular know when we have all of them.

\medskip

If $B$ is a block with Brauer tree a line, and with simple modules $S_1$, $S_2$,\dots, $S_e$, listed in order along the line, and with exceptionality $1$, then the structure of the projectives is
\[ P(S_i)=S_i/S_{i-1},S_{i+1}/S_i,\]
for $1<i<e$, and
\[ P(S_1)=S_1/S_2/S_1,\qquad P(S_e)=S_e/S_{e-1}/S_e.\]
This yields obvious indecomposable modules $S_i$, $P(S_i)$, $S_i/S_{i-1}$ (and its dual) and $S_i/S_{i-1},S_{i+1}$ (and its dual) for all appropriate $i$. The remaining simple modules (up to duality) have two socle layers, with socle
\[ S_i,S_{i+2},S_{i+4},\dots,S_{i+2m}\]
and top
\[ S_{i+1},S_{i+3},\dots,S_{i+2m\pm 1},\]
for some $m$ so that all indices are at most $e$. To construct these modules, we generalize in the obvious way the construction of the module $S_1,S_3/S_2,S_4$ obtained by taking the sum of $S_1,S_3/S_2$ and $S_3/S_4$ and taking the kernel of the appropriate diagonal homomorphism onto $S_3$. This process constructs all non-projective indecomposable modules, up to duality.

We quickly count the number of modules so obtained. We have $e$ projective modules, $e$ simple modules, and for each $1\leq i\leq e$, exactly $2(e-i)$ non-simple indecomposable modules, the above modules and their duals. This yields
\[ 2e+2\sum_{i=1}^e (e-i)=2e+e(e-1)=e(e+1)=(|D|-1)e,\]
as needed.

If the Brauer tree of $B$ is a line, but this time the exceptionality of $B$ is $2$ and the exceptional vertex is at the end of the Brauer tree of $B$, then we have exactly the same indecomposable modules as above, except that $S_i=S_{e-i+1}$ for all $i$.

The Brauer tree of any block of any alternating group is always a line, with exceptionality $1$ or $2$, and if it is $2$ then the exceptional vertex lies at the end of the tree, so this covers all cases.

\medskip

One final point to make is the following lemma.

\begin{lem}\label{lem:covertrivialBrauer} Let $B$ be a block with cyclic defect group, whose Brauer tree is a line, Suppose that the exceptionality of $B$ is $1$, or it is $2$ and the exceptional vertex has valency $1$. If $S$ is a simple $B$-module and $M$ is an indecomposable module with $S$ as a composition factor but with no submodule or quotient isomorphic to $S$, then $M=P(T)$ for some $T$, and $T$ and $S$ label edges that share a vertex of the Brauer tree. In particular, if $S$ is the trivial module then $M$ is uniquely determined, and has the structure $T/1,U/T$ for some simple $B$-module $U$.
\end{lem}

Since we often want to prove that $V_{\min}$ has a trivial submodule or quotient, such a lemma is of tremendous importance.

\section{Modules for alternating groups}
\label{sec:modulesalt}
In this section we summarize some information about modules for alternating groups. We begin by giving in Table \ref{tab:altmodules} all simple modules in all characteristics for all alternating groups $\Alt(n)$ for $5\leq n\leq 9$, and for $p=2,3$. The first row is the principal block, the second row is modules in the second block, and so on. The only exception to this is when we write `$a$ and $b$', where these are two projective simple modules (hence appearing in different blocks) and we have done so simply to save space.
\begin{table}\begin{center}
\begin{tabular}{ccc}
\hline Group & $p=2$ & $p=3$
\\\hline $\Alt(5)$ & $1$, $2_1$, $2_2$ & $1$, $4$
\\ & $4$ & $3_1$ and $3_2$
\\\hline $\Alt(6)$ & $1$, $4_1$, $4_2$ & $1$, $3_1$, $3_2$, $4$
\\ & $8_1$ and $8_2$ & $9$
\\\hline $\Alt(7)$ &$1$, $14$, $20$&$1$, $10$, $10^*$, $13$
\\ &$4_1$, $4_2$, $6$& $6$, $15$
\\ \hline $\Alt(8)$ & $1$, $4$, $4^*$, $6$, $14$, $20$, $20^*$ & $1$, $7$, $13$, $28$, $35$
\\ &$64$ & $21$
\\ && $45$ and $45^*$
\\\hline $\Alt(9)$ & $1$, $8_1$, $8_2$, $20$, $20^*$, $26$, $78$ & $1$, $7$, $21$, $35$, $41$
\\ & $8_3$, $48$, $160$ & $27$, $189$
\\ && $162$
\\ \hline 
\end{tabular}
\end{center}\caption{Simple modules in characteristics $2$ and $3$ for alternating groups}\label{tab:altmodules}
\end{table}
Here, the outer automorphism of $\Alt(n)$ induced by the transposition $(1,2)$ in the symmetric group swaps all simple modules that are not self dual with their duals, and swaps $i_1$ and $i_2$ for all $i$ apart from $4_1$ and $4_2$ for $\Alt(6)$ in characteristic $2$.

We also need to describe $1$-cohomology for these simple modules if we want to apply Proposition \ref{prop:lowpressuremodules}.

\begin{prop}\label{prop:ext1simples}
Let $5\leq n\leq 9$, and $p=2,3$. If $M$ is a simple module for $\Alt(n)$ in characteristic $p$ then $H^1(\Alt(n),M)=0$ unless $M$ appears in Table \ref{tab:exts}. Here, all $1$-cohomologies are $1$-dimensional unless the module is in bold, in which case it is $2$-dimensional.
\end{prop}

\begin{table}
\begin{center}\begin{tabular}{ccc}
\hline Group & $p=2$ & $p=3$
\\\hline $\Alt(5)$ & $2_1$, $2_2$ & $4$
\\\hline $\Alt(6)$ & $4_1$, $4_2$ & $\b{4}$
\\\hline $\Alt(7)$ & $14$, $20$& $10$, $10^*$, $\b{13}$
\\ \hline $\Alt(8)$ & $6$, $14$, $20$, $20^*$ & $13$, $35$
\\\hline $\Alt(9)$ & $20$, $20^*$, $\b{26}$, $78$ & $7$, $35$, $41$
\\ \hline 
\end{tabular}\end{center}\caption{Simple modules with non-trivial $1$-cohomology (bold indicates $2$-dimensional)}\label{tab:exts}
\end{table}

For $p=5$ and $p=7$, the Sylow $p$-subgroups are cyclic and so we can simply describe the Brauer trees. If there is an exceptional node we represent it by filling in the node. By $\chi_i$ we mean an ordinary character of degree $i$. We start with $p=5$.

For $\Alt(5)$, we have a principal block and a single projective simple module of dimension $5$.
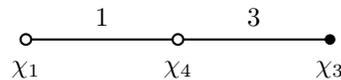
\begin{center}\begin{tikzpicture}[thick,scale=2]
\draw (0,0.8) node {$\chi_1$};
\draw (1,0.8) node {$\chi_{4}$};
\draw (2,0.8) node {$\chi_{3}$};

\draw (0,1) -- (2,1);
\draw (0,1) node [draw] (l2) {};
\draw (1,1) node [draw] (l2) {};
\draw (2,1) node [fill=black!100] (l2) {};
\draw (0.5,1.15) node{$1$};
\draw (1.5,1.15) node{$3$};
\end{tikzpicture}\end{center}

For $\Alt(6)$ there is the principal block with the following Brauer tree, and three projective simple modules, $5_1$, $5_2$ (not permuted by the $\Sym(6)$ outer automorphism, but are by the other two) and $10$.
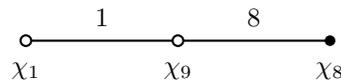
\begin{center}\begin{tikzpicture}[thick,scale=2]
\draw (0,0.8) node {$\chi_1$};
\draw (1,0.8) node {$\chi_{9}$};
\draw (2,0.8) node {$\chi_{8}$};

\draw (0,1) -- (2,1);
\draw (0,1) node [draw] (l2) {};
\draw (1,1) node [draw] (l2) {};
\draw (2,1) node [fill=black!100] (l2) {};
\draw (0.5,1.15) node{$1$};
\draw (1.5,1.15) node{$8$};
\end{tikzpicture}\end{center}

For $\Alt(7)$ there is the principal block with the following Brauer tree, and four projective simple modules, $10$, $10^*$, $15$ and $35$.
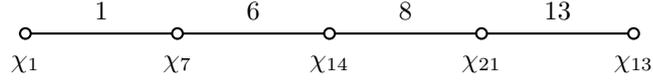
\begin{center}\begin{tikzpicture}[thick,scale=2]
\draw (0,0.8) node {$\chi_1$};
\draw (1,0.8) node {$\chi_{7}$};
\draw (2,0.8) node {$\chi_{14}$};
\draw (3,0.8) node {$\chi_{21}$};
\draw (4,0.8) node {$\chi_{13}$};

\draw (0,1) -- (4,1);
\draw (1,1) node [draw] (l2) {};
\draw (2,1) node [draw] (l2) {};
\draw (3,1) node [draw] (l2) {};
\draw (4,1) node [draw] (l2) {};
\draw (0,1) node [draw] (ld) {};
\draw (0.5,1.15) node{$1$};
\draw (1.5,1.15) node{$6$};
\draw (2.5,1.15) node{$8$};
\draw (3.5,1.15) node{$13$};
\end{tikzpicture}\end{center}
We also need the two dual faithful blocks of $3\cdot \Alt(7)$, which have the same structure.
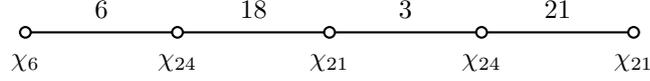
\begin{center}\begin{tikzpicture}[thick,scale=2]
\draw (0,0.8) node {$\chi_6$};
\draw (1,0.8) node {$\chi_{24}$};
\draw (2,0.8) node {$\chi_{21}$};
\draw (3,0.8) node {$\chi_{24}$};
\draw (4,0.8) node {$\chi_{21}$};

\draw (0,1) -- (4,1);
\draw (1,1) node [draw] (l2) {};
\draw (2,1) node [draw] (l2) {};
\draw (3,1) node [draw] (l2) {};
\draw (4,1) node [draw] (l2) {};
\draw (0,1) node [draw] (ld) {};
\draw (0.5,1.15) node{$6$};
\draw (1.5,1.15) node{$18$};
\draw (2.5,1.15) node{$3$};
\draw (3.5,1.15) node{$21$};
\end{tikzpicture}\end{center}

For $\Alt(8)$, we have the principal block and a non-principal block with the two trees below, plus five projective simple modules, $20$, $35$, $45$, $45^*$ and $70$.

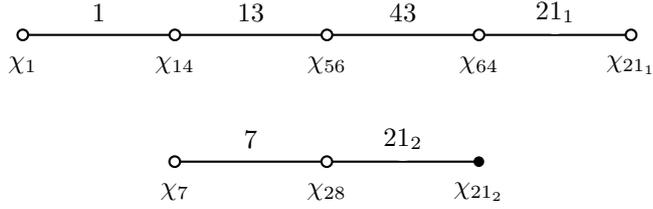
\begin{center}\begin{tikzpicture}[thick,scale=2]
\draw (0,0.8) node {$\chi_1$};
\draw (1,0.8) node {$\chi_{14}$};
\draw (2,0.8) node {$\chi_{56}$};
\draw (3,0.8) node {$\chi_{64}$};
\draw (4,0.8) node {$\chi_{21_1}$};

\draw (0,1) -- (4,1);
\draw (1,1) node [draw] (l2) {};
\draw (2,1) node [draw] (l2) {};
\draw (3,1) node [draw] (l2) {};
\draw (4,1) node [draw] (l2) {};
\draw (0,1) node [draw] (ld) {};
\draw (0.5,1.15) node{$1$};
\draw (1.5,1.15) node{$13$};
\draw (2.5,1.15) node{$43$};
\draw (3.5,1.15) node{$21_1$};
\end{tikzpicture}\end{center}

\begin{center}\begin{tikzpicture}[thick,scale=2]
\draw (0,0.8) node {$\chi_7$};
\draw (1,0.8) node {$\chi_{28}$};
\draw (2,0.8) node {$\chi_{21_2}$};

\draw (0,1) -- (2,1);
\draw (0,1) node [draw] (l2) {};
\draw (1,1) node [draw] (l2) {};
\draw (2,1) node [fill=black!100] (l2) {};
\draw (0.5,1.15) node{$7$};
\draw (1.5,1.15) node{$21_2$};
\end{tikzpicture}\end{center}

We now turn to $p=7$. Of course, we only need to consider $\Alt(7)$ and $\Alt(8)$. For $\Alt(7)$, apart from the principal block we have four projective simple modules, $14_1$, $14_2$ (not swapped by the outer automorphism), $21$ and $35$.
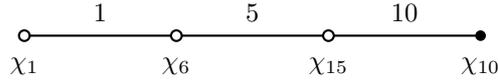
\begin{center}\begin{tikzpicture}[thick,scale=2]
\draw (0,0.8) node {$\chi_1$};
\draw (1,0.8) node {$\chi_{6}$};
\draw (2,0.8) node {$\chi_{15}$};
\draw (3,0.8) node {$\chi_{10}$};

\draw (0,1) -- (3,1);
\draw (0,1) node [draw] (l2) {};
\draw (1,1) node [draw] (l2) {};
\draw (2,1) node [draw] (l2) {};
\draw (3,1) node [fill=black!100] (l2) {};
\draw (0.5,1.15) node{$1$};
\draw (1.5,1.15) node{$5$};
\draw (2.5,1.15) node{$10$};
\end{tikzpicture}\end{center}

For $\Alt(8)$, apart from the principal block we have nine projective simple modules, namely $7$, $14$, $21_1$, $21_2$, $21_2^*$, $28$, $35$, $56$ and $70$, with the outer automorphism swapping $21_2$ and $21_2^*$.
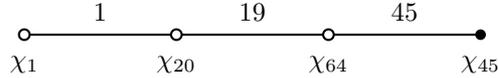
\begin{center}\begin{tikzpicture}[thick,scale=2]
\draw (0,0.8) node {$\chi_1$};
\draw (1,0.8) node {$\chi_{20}$};
\draw (2,0.8) node {$\chi_{64}$};
\draw (3,0.8) node {$\chi_{45}$};

\draw (0,1) -- (3,1);
\draw (0,1) node [draw] (l2) {};
\draw (1,1) node [draw] (l2) {};
\draw (2,1) node [draw] (l2) {};
\draw (3,1) node [fill=black!100] (l2) {};
\draw (0.5,1.15) node{$1$};
\draw (1.5,1.15) node{$19$};
\draw (2.5,1.15) node{$45$};
\end{tikzpicture}\end{center}

To end this section, we give a couple of facts about specific modules for $\Alt(5)$ and $\Alt(8)$ in characteristic $2$.

\begin{lem}\label{lem:numberof2sneeded} Let $M$ be a module for $H\cong\Alt(5)$ in characteristic $2$, and suppose that $M$ has no trivial submodules or quotients. If $u\in H$ is an involution, suppose that $u$ acts on $M$ with exactly $a$ Jordan blocks of size $1$, and there are exactly $a+2b$ trivial composition factors in $M$. There are at least $2a+3b$ composition factors of dimension $2$ in $M$.
\end{lem}
\begin{pf}
Firstly, the projective cover $P(2_i)$ has the form $2_i/1/2_{3-i}/1/2_i$, and since we may assume that $H$ does not fix a line or hyperplane on $M$, it is a submodule of copies of $P(2_1)$ and $P(2_2)$, together with copies of $4$, which break off as they are projective. If $M$ has five socle layers then there must be a projective summand, and so we write $M$ as $U\oplus W$, where $U$ has no projective summands and $W$ is projective. Notice that $U$ has at most three socle layers, since it cannot have five by above, and if it had four then it would have a trivial quotient. 

For $U$ not to fix a line or hyperplance, $U$ must have at least twice as many $2$-dimensional factors as trivials, so $2a$. For $W$, $P(1)$ is not a summand, and so we have copies of $P(2_1)$ and $P(2_2)$, which have three $2_i$s for each trivial composition factor. This proves the result.
\end{pf}

\begin{lem}\label{lem:alt8no20s}
Let $G=\Alt(8)$ and $k$ be a field of characteristic $2$. Suppose that $M$ is a $kG$-module with no composition factors of dimension $20$. If $M$ has no trivial submodules or quotients, and $M$ has $n$ trivial composition factors, then $M$ also has composition factors of either $14^n,6^{n+1}$ or $14^{n+1},6^n$.
\end{lem}
\begin{pf}
Since $\Ext^1(1,4)=\Ext^1(1,4^*)=0$, we can remove any $4$s and $4^*$s from the top and bottom of $M$, so that $M$ is a submodule of a sum of $P(6)$ and $P(14)$.

Of course, we can perform the same reductions to the top as well, so we can assume that this submodule has top consisting entirely of $6$s and $14$s. We therefore examine the largest submodules of $P(14)$ and $P(6)$ containing no $20$s or $20^*$s. For $P(14)$ this has a single trivial composition factor, and the smallest submodule containing it with no trivial quotient is
\[ 14/1,6/14,\]
which satisfies the statement, and for $P(6)$ the module has exactly two trivial composition factors, and so we need to examine this. The smallest submodule containing both trivials and with no trivial quotient is given by the following diagram.

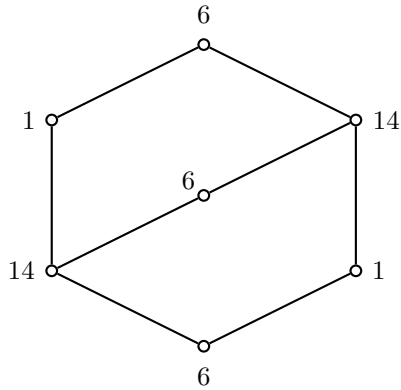
\begin{center}\begin{tikzpicture}[thick,scale=2]
\draw (0,-0.2) node {$6$};
\draw (0,2.2) node {$6$};
\draw (-1.15,1.5) node {$1$};
\draw (1.15,0.5) node {$1$};
\draw (-0.1,1.1) node {$6$};
\draw (1.2,1.5) node {$14$};
\draw (-1.2,0.5) node {$14$};

\draw (0,0) node [draw] (l1) {};
\draw (1,0.5) node [draw] (l2) {};
\draw (-1,0.5) node [draw] (l3) {};
\draw (0,1) node [draw] (l4) {};
\draw (-1,1.5) node [draw] (l5) {};
\draw (1,1.5) node [draw] (l6) {};
\draw (0,2) node [draw] (l7) {};
\draw (l1) -- (l2);
\draw (l1) -- (l3);
\draw (l3) -- (l4);
\draw (l4) -- (l6);
\draw (l5) -- (l7);
\draw (l6) -- (l7);
\draw (l3) -- (l5);
\draw (l2) -- (l6);
\end{tikzpicture}\end{center}
Since there are no uniserial modules $6/1/6$ or $14/1/6$, the statement clearly holds for submodules with a single trivial factor. As $M$ contains a diagonal submodule of one of these as a submodule, the statement holds for $M$.
\end{pf}

\section{Alternating subgroups in positive-dimensional subgroups}

In this section we classify the possible composition factors of $H\cong\Alt(n)$ on $\vmin$ and $L(G)$ for $G=F_4,E_6,E_7,E_8$ in characteristic $2$, for $n=7,8,9$. This is particularly useful when attacking $\Alt(8)$ and $\Alt(9)$, because there are many possible sets of composition factors (i.e., consistent with the traces of semisimple elements) that do not actually occur. We establish the base case of the induction in Theorem \ref{thm:alt7} by showing that every $\Alt(7)$ occurs in a positive-dimensional subgroup, hence easily classifiable. We then restrict the composition factors for $\Alt(8)$ and $\Alt(9)$ to $\Alt(7)$ to cut down substantially on the amount of work we have to do in later sections to prove Theorem \ref{thm:alt8+}.

\begin{prop}\label{prop:alt7compfactors} Let $G$ be one of $F_4$, $E_6$, $E_7$ and $E_8$ in characteristic $2$, and let $H\cong\Alt(7)$.
\begin{enumerate}
\item If $G=F_4$ then the possible composition factors of $H$ on $V_{\min}$ and $L(G)$ are as follows:
\begin{center}
\begin{tabular}{cc}
\hline $V_{\min}$ & $L(G)$
\\\hline $6,(4,4^*)^2,1^4$ & $14,6^3,(4,4^*)^2,1^4$
\\ $14,6^2$ & $14,6^3,(4,4^*)^2,1^4$
\\\hline \end{tabular}
\end{center}
\item If $G=E_6$ then the possible composition factors of $H$ on $V_{\min}$ and $L(G)$ are as follows:
\begin{center}
\begin{tabular}{cc}
\hline $V_{\min}$ & $L(G)$
\\\hline $6,(4,4^*)^2,1^5$ & $14,6^4,(4,4^*)^4,1^8$
\\ $14,6^2,1$ & $14^2,6^5,(4,4^*)^2,1^4$
\\\hline
$15,6',6'$ & $20,14,6^4,(4,4^*)^2,1^4$
\\\hline \end{tabular}
\end{center}
(The first two are embeddings of $\Alt(7)$ into $E_6$, the last is $3\cdot \Alt(7)\leq E_6$.)
\item If $G=E_7$ then the possible composition factors of $H$ on $V_{\min}$ and $L(G)$ are as follows:
\begin{center}
\begin{tabular}{cc}
\hline $V_{\min}$ & $L(G)$
\\\hline $6^2,(4,4^*)^4,1^{12}$ & $14,6^6,(4,4^*)^8,1^{19}$
\\ $14^2,6^4,1^4$ & $14^4,6^9,(4,4^*)^2,1^7$
\\ $6^8,4,4^*$ & $14^8,6,1^{15}$
\\\hline \end{tabular}
\end{center}
\item If $G=E_8$ then the possible composition factors of $H$ on $L(G)$ are as follows:
\begin{center}
\begin{tabular}{c}
\hline $L(G)$
\\\hline $14,6^{10},(4,4^*)^{16},1^{46}$
\\ $14^8,6^{17},(4,4^*)^2,1^{18}$
\\ $20^4,14^4,6^8,(4,4^*)^7,1^8$
\\ $20^4,14^{10},6^2,4,4^*,1^8$
\\\hline \end{tabular}
\end{center}
\end{enumerate}
\end{prop}
\begin{pf} There is no embedding of $H\cong\Alt(7)$ into subgroups of type $A_1$, $A_2$, $B_2$ and $G_2$, and so if we see these factors in an algebraic group we can remove them. Since we only considering the sets of composition factors we may replace parabolic subgroups by their Levis, and so may assume that $X$, a positive-dimensional subgroup of $G$, is a product of simple groups not of type $A_1$, $A_2$, $B_2$ and $G_2$.

\bigskip

Inside $F_4$, we know that there are five potential sets of composition factors by \cite{litterickmemoir}, each of which lies inside a positive-dimensional subgroup by Theorem \ref{thm:alt7}. Let $X$ be a maximal such subgroup: if $X=B_4$ then we can embed $H$ as $4,4^*,1$ or $6,1^3$ on the natural, both of which yield $6,(4,4^*)^2,1^4$ on $V_{\min}$ and $14,6^3,(4,4^*)^2,1^4$ on $L(G)$.

If $X=C_4$ then there are two sets of composition factors on the natural, namely $4,4^*$ and $6,1^2$, both of which yield $14,6^2$ on $L(\lambda_2)$ (which is the restriction of $V_{\min}$ to $X$), and these factors determine the action on $L(G)$. Since $H\not\leq A_2\tilde A_2$ this case need not be considered.

The only Levi subgroups that can be considered are $B_3$ and $C_3$, which lie in $B_4$ and $C_4$ respectively, and so we are done.

\bigskip

For $G=E_6$, since $H$ must fix a line on $V_{\min}$ we have that $H\leq F_4$ or $H\leq D_5$. If $H\leq F_4$ then we know the answer, and if $H\leq D_5$ then the factors on the natural are $6,1^4$ or $4,4^*,1^2$; in either case, $H$ fixes a line on the natural, so that $H\leq B_4$, but this is contained in $F_4$ and we are done.

We could also have $H\cong3\cdot \Alt(7)$ contained in $E_6$ with centres coinciding. By the appropriate table in \cite{litterickmemoir} we see that the composition factors on $V_{\min}$ have dimensions $15,6,6$, and this means that $H$ must lie in the $A_5$ Levi subgroup, acting irreducibly on the natural module. The composition factors of this embedding on $L(G)$ are easily calculable (as they are two copies of $\Lambda^3(M)$ and one of $M\otimes M^*$, where $M$ is the natural module for the $A_5$), and are as appear in the table above.

\bigskip

For $G=E_7$, $H\cong\Alt(7)$ either fixes a line on $V_{\min}$ or has composition factors $6^8,4,4^*$, by \cite{litterickmemoir}. The line stabilizers are contained in an $E_6$ parabolic or are subgroups $q^{1+32}(q-1)B_5(q)$, which acts on $V_{\min}$ with composition factors $32,11^2,1^2$. Since the non-trivial irreducibles for $H$ have even dimension, we get in the latter case that $H$ fixes a line on the $11$ so lies inside $D_5\leq E_6$ anyway. We therefore get those inherited from $E_6$, which are the first two rows of the table. For the final row, this is the case stated above: if $H$ is embedded as $4\oplus 4$ inside the $A_7$ maximal-rank subgroup then it acts on $V_{\min}$ as the sum of $\Lambda^2(4)^{\oplus 2}\oplus 4\otimes 4$ and its dual, but this is $6^{\oplus 2}\oplus 6/4/6$, as needed.

\bigskip

For $G=E_8$, we know that $H$ lies inside a positive-dimensional subgroup, so we consider the maximal such subgroups $X$, which are given in \cite{liebeckseitz2004}. Using Lemma \ref{lem:insideconnected} we can exclude those subgroups $X\in \ms X$ with $X^0\leq Y$ for some other $Y\in \ms X$, for example the normalizer of a torus.

We can remove any quotient from $X$ that does not contain a copy of $\Alt(7)$, so for example replace $(D_4^2).(\Sym(3)\times 2)$ by $D_4^2$.

If $H$ lies in $X=E_7$ then we get the first two cases, as the second and third case for $E_7$ yield the same factors for $E_8$. Thus we may assume that $H$ does not lie inside $E_7$. This eliminates the $E_7A_1$, $E_6A_2$ and $F_4G_2$ reductive subgroups.

For $X=A_8$, we must have $H\leq A_8$ since the faithful irreducible modules for $3\cdot \Alt(7)$ have dimensions $6$, $15$ and $24$. We can embed $\Alt(7)$ inside $A_8$ with factors
\[ 4^2,1\qquad 4,4^*,1\qquad 6,1^3,\qquad 4,1^5.\]
These result in composition factors on $L(G)$ of
\[ 20^4,14^4,6^8,(4,4^*)^7,1^8,\qquad 20^4,14^4,6^8,(4,4^*)^7,1^8,\qquad 14^8,6^{17},(4,4^*)^2,1^{18},\qquad 14,6^{10},(4,4^*)^{16},1^{46},\]
respectively.

For $X=D_8$, we run through the possible actions of $H$ on the natural $16$, noting that it cannot have more than three trivial composition factors on this module as else it would lie inside $D_6\leq E_7$. There are only three possibilities, namely
\[ (4,4^*)^2,\qquad 6,4,4^*,1^2,\qquad 14,1^2.\]

Using the traces of elements of orders $3$, $5$ and $7$, and comparing them with the nine classes of elements of order $3$, the $53$ classes of elements of order $5$, and the $209$ classes of elements of order $7$ from $D_8$, we find that the composition factors on $L(\lambda_7)$ are determined uniquely for the second and third embeddings, but there are two possibilities for the first embedding.

\begin{center}
\begin{tabular}{ccc}
\hline Factors on $L(\lambda_1)$ & Factors on $L(\lambda_2)$ & Factors on $L(\lambda_7)$
\\\hline $(4,4^*)^2$ & $14^4,6^8,4,4^*,1^8$ & $20^4,(4,4^*)^6$ or $14^4,6^9,4,4^*,1^{10}$
\\ $6,4,4^*,1^2$ & $20^2,14^2,6^4,(4,4^*)^3,1^4$ & $20^2,14^2,6^4,(4,4^*)^4,1^4$
\\ $14,1^2$ & $20^2,14^4,6^2,4,4^*,1^4$ & $20^2,14^6,1^4$
\\\hline\end{tabular}
\end{center}

The sum of the second and third columns give the action on $L(G)$, which contributes the last line to the table in the result.

For $X=A_4A_4$, the only action of $H$ on a $5$-space is as $4\oplus 1$ (or $4^*\oplus 1$), and so $H\leq A_3A_3\leq A_8$, so we have already done this case.

For the maximal Levi subgroups, they are
\[ D_7,\quad A_7,\quad A_1A_6,\quad A_1A_2A_4,\quad A_4A_3,\quad D_5A_2,\quad E_6A_1,\quad E_7,\]
and after removing $A_1$s and $A_2$s, these are contained within $D_8$, $A_8$, $A_8$, $A_8$, $A_4A_4$, $E_7$, $E_7$ and $E_7$ respectively, so we have gone through the complete list of maximal positive-dimensional subgroups of $G$.
\end{pf}

\begin{prop}\label{prop:alt8compfactors} Let $G$ be one of $F_4$, $E_6$, $E_7$ and $E_8$ in characteristic $2$, and let $H\cong\Alt(8)$.
\begin{enumerate}
\item If $G=F_4$ then the possible composition factors of $H$ on $V_{\min}$ and $L(G)$ are as follows:
\begin{center}
\begin{tabular}{cc}
\hline $V_{\min}$ & $L(G)$
\\\hline $6,(4,4^*)^2,1^4$ & $14,6^3,(4,4^*)^2,1^4$
\\ $14,6^2$ & $14,6^3,(4,4^*)^2,1^4$
\\\hline \end{tabular}
\end{center}
\item If $G=E_6$ then the possible composition factors of $H$ on $V_{\min}$ and $L(G)$ are as follows:
\begin{center}
\begin{tabular}{cc}
\hline $V_{\min}$ & $L(G)$
\\\hline $6,(4,4^*)^2,1^5$ & $14,6^4,(4,4^*)^4,1^8$
\\ $14,6^2,1$ & $14^2,6^5,(4,4^*)^2,1^4$
\\\hline \end{tabular}
\end{center}
\item If $G=E_7$ then the possible composition factors of $H$ on $V_{\min}$ and $L(G)$ are as follows:
\begin{center}
\begin{tabular}{cc}
\hline $V_{\min}$ & $L(G)$
\\\hline $6^2,(4,4^*)^4,1^{12}$ & $14,6^6,(4,4^*)^8,1^{19}$
\\ $14^2,6^4,1^4$ & $14^4,6^9,(4,4^*)^2,1^7$
\\ $6^8,4,4^*$ & $14^8,6,1^{15}$
\\\hline \end{tabular}
\end{center}
\item If $G=E_8$ then the possible composition factors of $H$ on $L(G)$ are as follows:
\begin{center}
\begin{tabular}{c}
\hline $L(G)$
\\\hline $14,6^{10},(4,4^*)^{16},1^{46}$
\\ $14^8,6^{17},(4,4^*)^2,1^{18}$
\\ $(20,20^*)^2,14^4,6^8,(4,4^*)^7,1^8$
\\ $64^2,20,20^*,14^4,6^2,4,4^*,1^2$
\\\hline \end{tabular}
\end{center}
\end{enumerate}
\end{prop}
\begin{pf}
By the construction of the $\Alt(7)$s inside $G=F_4$, $E_6$, $E_7$ and $E_8$, for each set of composition factors there is an $\Alt(7)$ representing these factors that extends to an $\Alt(8)$ inside $G$. We need to determine the composition factors of that extension.

Notice that the restrictions to $\Alt(7)$ of all but one of the simple $\Alt(8)$-modules are simple, and are unique except that $20$ and $20^*$ both restrict to the same module. The last module is $64$, and the restriction to $\Alt(7)$ has factors $20,14^3,1^2$.

Thus if $\Alt(7)$ extends to $H\cong\Alt(8)$ and does not have $20$ as a composition factor on a module $M$, then the composition factors of the $\Alt(8)$ on $M$ are uniquely determined. This completes the list for $G\neq E_8$, and for two entries for $G=E_8$.

Notice that, since $L(G)$ is self dual, we cannot have an odd number of $64$s as then we would be left with an odd number of $20$s in $L(G)\downarrow_H$. This deals with the third possibility for the action of $\Alt(7)$ on $L(G)$ given in Proposition \ref{prop:alt7compfactors}. We are left with $\Alt(7)$ acting as $20^4,14^{10},6^2,4,4^*,1^8$, which can act as either
\[ (20,20^*)^2,14^{10},6^2,4,4^*,1^8,\qquad 64^2,20,20^*,14^4,6^2,4,4^*,1^2.\]

Using the traces of semisimple elements of order up to $15$, the potential set of composition factors for $\Alt(8)$s inside $E_8$ is
\[ 14^8,6^{17},(4,4^*)^2,1^{18},\quad 20,20^*,14^6,6^{10},(4,4^*)^7,1^8,\quad (20,20^*)^2,14^4,6^8,(4,4^*)^7,1^8,\]
\[ 64,20,20^*,14^5,6^6,(4,4^*)^4,1^6,\quad 64,(20,20^*)^2,14^3,6^4,(4,4^*)^4,1^6,\quad 64^2,20,20^*,14^4,6^2,4,4^*,1^4,\]
\[ 64^2,(20,20^*)^2,14^2,4,4^*,1^4,\quad 14,6^{10},(4,4^*)^{16},1^{46}.\]
Only one of the two possibilities we generated above is on this list, and so we are done.
\end{pf}

This proposition \emph{does not require} the case $n=8$ from Theorem \ref{thm:alt8+} to prove it and so it can be used in the proof of Theorem \ref{thm:alt8+}.

Similarly, the next proposition does not require the $n=9$ case of Theorem \ref{thm:alt8+}.

\begin{prop}\label{prop:alt9compfactors} Let $G$ be one of $F_4$, $E_6$, $E_7$ and $E_8$ in characteristic $2$, and let $H\cong\Alt(9)$.
\begin{enumerate}
\item If $G=F_4$ then the possible composition factors of $H$ on $V_{\min}$ and $L(G)$ are as follows:
\begin{center}
\begin{tabular}{cc}
\hline $V_{\min}$ & $L(G)$
\\\hline $8_1,8_2,8_3,1^2$ & $26,8_1,8_2,8_3,1^2$
\\ $26$ & $26,8_1,8_2,8_3,1^2$
\\\hline \end{tabular}
\end{center}
\item If $G=E_6$ then the possible composition factors of $H$ on $V_{\min}$ and $L(G)$ are as follows:
\begin{center}
\begin{tabular}{cc}
\hline $V_{\min}$ & $L(G)$
\\\hline $8_1,8_2,8_3,1^3$ & $26,8_1^2,8_2^2,8_3^2,1^4$
\\ $26,1$ & $26^2,8_1,8_2,8_3,1^2$
\\\hline \end{tabular}
\end{center}
\item If $G=E_7$ then the possible composition factors of $H$ on $V_{\min}$ and $L(G)$ are as follows:
\begin{center}
\begin{tabular}{cc}
\hline $V_{\min}$ & $L(G)$
\\\hline $8_1^2,8_2^2,8_3^2,1^8$ & $26,8_1^4,8_2^4,8_3^4,1^{11}$
\\ $26^2,1^4$ & $26^4,8_1,8_2,8_3,1^5$
\\\hline \end{tabular}
\end{center}
\item If $G=E_8$ then the possible composition factors of $H$ on $L(G)$ are as follows:
\begin{center}
\begin{tabular}{c}
\hline $L(G)$
\\\hline $26,8_1^8,8_2^8,8_3^8,1^{30}$
\\ $26^8,8_1,8_2,8_3,1^{16}$
\\ $26^4,(20,20^*)^2,8_1^5,8_2^2,1^8$
\\ $48^2,26^4,8_3^5,1^8$
\\ $48,26^2,(20,20^*)^2,8_1^3,8_2^3,8_3^2,1^4$
\\\hline \end{tabular}
\end{center}
\end{enumerate}
(Here, $8_1$ and $8_2$ are swapped under the outer automorphism of $\Alt(9)$.)
\end{prop}
\begin{pf} Since $\Alt(9)\leq F_4(2)$ by, for example, \cite{liebeckseitz1999} (see also \cite{litterickmemoir}), we consult the tables in \cite{litterickmemoir} to find its composition factors, and see that there are two classes swapped by the graph automorphism. These classes propagate through $E_6$ and $E_7$, and using \cite{litterickmemoir} we see that there are no more possible sets of composition factors, completing the proof for these groups. Hence $G=E_8$.

If $H\leq X=E_7$ then we get the first two lines of the tables in the proposition, so we may assume that $H$ is not contained in $E_7$.

If $H\leq A_8$ then as $8_1$ and $8_2$ are swapped by the outer automorphism of $H$, so we get that $H$ acts on the natural module with factors either $8_1,1$ or $8_3,1$: these yield
\[ 26^4,(20,20^*)^2,8_1^5,8_2^2,1^8, \qquad 48^2,26^4,8_3^5,1^8\]
on $L(G)$, yielding the next two rows.

The only other positive-dimensional subgroup in which $H$ can lie is $D_8$. If $H$ acts on the $16$ with eight trivial factors then $H\leq D_6\leq E_7$, so we may assume that the factors on the $16$ are $8_i,8_j$. 

Examining the table in \cite{litterickmemoir}, and using Proposition \ref{prop:alt8compfactors}, we can determine the potential sets of composition factors that can occur, given their restrictions to $\Alt(8)$ must exist. As well as those already found, there are two more:
\[ 78^2,48,20,20^*,1^4,\qquad 48,26^2,(20,20^*)^2,8_1^3,8_2^3,8_3^2,1^4.\]
The exterior square of $8_i^{\oplus 2}$ is $26^4,8_i,1^8$, and so cannot be either of the two remaining sets of composition factors. The other cases are $8_1\oplus 8_2$ and $8_1\oplus 8_3$: their exterior squares are
\[ \Lambda^2(8_1\oplus 8_2)=48,26^2,8_3^2,1^4,\qquad \Lambda^2(8_1\oplus 8_3)=26^2,20,20^*,8_1,8_2^2,1^4.\]
This means the first case cannot occur in any positive-dimensional subgroup. To see whether the second can occur we simply compute its Brauer character and compare it to that produced by $8_1\oplus 8_2$ and $8_1\oplus 8_3$, computing traces using the program for semisimple elements in \cite{litterick}.

There are $650$ classes of semisimple elements of order $9$ and $9375$ of order $15$ in $D_8$. Using the traces of these elements it is possible to construct the Brauer character of $H$ on $L(G)$. If $H$ acts as $8_1\oplus 8_2$ on the natural then there are two options for the Brauer character on $L(G)$, arising from the modules $48^2,26^4,8_3^5,1^8$ and $48,26^2,(20,20^*)^2,8_1^3,8_2^3,8_3^2,1^4$. For $8_1\oplus 8_3$ we have again two possible characters on $L(G)$, arising from $48,26^2,(20,20^*)^2,8_1^3,8_2^3,8_3^2,1^4$ again and also $26^4,(20,20^*)^2,8_1^5,8_2^2,1^8$ this time.

Now suppose that $H$ does not lie in a positive-dimensional subgroup of $G$. We are still yet to find the single case from \cite{litterickmemoir} allowed by Proposition \ref{prop:alt8compfactors} but not in any positive-dimensional subgroup. Since $H^1(H,M)$ has dimension $2$ for $M=26$, dimension $1$ for $M=20,20^*,78$, and $0$ otherwise, the first of these has a trivial submodule Proposition \ref{prop:lowpressuremodules}, and so lies inside a positive-dimensional subgroup of $G$, a contradiction.
\end{pf}

\section{$\Alt(5)$}
\label{sec:alt5}
There are four possibilities for primes $p$ when $H$ is $\Alt(5)$: $p=2,3,5$ and primes larger than $5$. Recall our assumption that $E_6$ and $E_7$ are used to denote the simply connected forms, so that $Z(E_6)$ and $Z(E_7)$ have orders $3$ and $2$ respectively.

\begin{prop} Suppose that $p\neq 2,3,5$.
\begin{enumerate}
\item If $G=F_4$ then there is a unique conjugacy class of subgroups isomorphic to $\Alt(5)$ that fixes a line on neither $\Vmin$ nor $L(G)$. This is contained in the $A_2\tilde A_2$ subgroup, and hence is not maximal.
\item If $G=E_6,E_7$, and $H\cong \Alt(5)$, then $H$ fixes a line on either $\Vmin$ or $L(G)$, and hence is not maximal.
\item If $G=E_7$ and $H\cong 2\cdot \Alt(5)$ with $Z(H)=Z(G)$, then either $H$ fixes a line on $L(G)$ or then $H$ fixes a $2$-space on $\Vmin$ and $\Vmin\downarrow_H$ is not $\Out(H)$-stable, and hence $H$ is not maximal.
\item If $G=E_8$ then there is a unique conjugacy class of subgroups isomorphic to $\Alt(5)$ that does not fix a line on $L(G)$. It is contained in a $D_8$ subgroup, and hence is not maximal.
\end{enumerate}
\end{prop}
\begin{pf}
In \cite{litterickmemoir} we find tables of all possible sets of composition factors for embeddings of $H\cong\Alt(5)$ into $G=F_4$, $E_6$, $E_7$ and $E_8$ for $p>5$, which immediately proves (ii) by Lemma \ref{lem:smallspacestabilizers} (as $H$ must act semisimply on $L(G)$ and $V_{\min}$). If we embed $2\cdot \Alt(5)$ into the simply connected version of $E_7$ with centres coinciding, then there is a single possibility that has no trivial composition factors on $L(G)$ from \cite[Table 6.137]{litterickmemoir}, and that has two isomorphic $2$-dimensional composition factors, but not their Frobenius twists, so that $H$ stabilizes a $2$-space inside $V_{\min}$ and so is not maximal by Lemma \ref{lem:smallspacestabilizers}; this proves (iii). (The condition of $\Out(H)$-stability is there to prove that in this case there is no almost simple subgroup with socle $H$ being maximal in an almost simple group of type $E_7$, via Lemma \ref{lem:novelty}(iii).)

For (i) and (iv), so $G$ is one of $F_4$ and $E_8$, we need to consider those possible embeddings of $H$ into $G$ with no trivial composition factors on either $V_{\min}$ or $L(G)$. In this case, there are only two cases to consider: one in $E_8$ and one in $F_4$.

For $G=E_8$, the fixed-point-free embedding was proved to be unique by Lusztig \cite{lusztig2002} and contained in $D_8$ by \cite[Table 7.6]{frey1998} (it is pattern 844). For $G=F_4$, this is proved to not be maximal by Lemma 16.4 of \cite{magaardphd}, and indeed is unique by \cite[Table 4.18]{frey1998}. It can be found inside the $A_2\tilde A_2$ subgroup, acting irreducibly along both factors as different $3$-dimensional modules. This completes the proof.
\end{pf}

\begin{prop}
Suppose that $p=5$. If $G$ is one of $F_4$, $E_6$, $E_7$ and $E_8$, and $H$ is $\Alt(5)$, or $H$ is $2\cdot \Alt(5)$ for $G=E_7$ with $Z(H)=Z(G)$, then $H$ fixes a line on either $\Vmin$ or $L(G)$, and hence is not maximal.
\end{prop}
\begin{pf}
For $p=5$ we have the Brauer tree in Section \ref{sec:modulesalt}, and Lemma \ref{lem:covertrivialBrauer} states that the only indecomposable module for $H$ that has trivial composition factors but no trivial submodules or quotients is $P(3)$, which has composition factors $3^3,1$. Thus if $M$ is any module for $H$ in characteristic $5$ with fewer than three times as many $3$s as $1$s, then $M$ has a trivial submodule or quotient. We now prove that this always occurs for $V_{\min}$ or $L(G)$. (Since this is a defining-characteristic embedding of $\Alt(5)=\PSL_2(5)$ the possible composition factors are not tabulated in \cite{litterickmemoir}.)

The traces of $t=(1,2)(3,4)$ and $x=(1,2,3)$ on the three irreducible modules are as follows:

\begin{center}
\begin{tabular}{ccc}
\hline Dimension of $M$ & Trace on $t$ & Trace on $x$
\\\hline $1$&$1$&$1$
\\ $3$&$-1$&$0$
\\ $5$&$1$&$-1$
\\ \hline
\end{tabular}
\end{center}

For $F_4$ the traces of $t$ and $x$ are given in Lemma \ref{lem:traces}, and this yields the following possible sets of composition factors for $H$ acting on $V_{\min}$:
\[ 3^6,1^8,\qquad 5,3^7,\qquad 5^3,3^3,1^2.\]
Since we need at least three times as many $3$s as $1$s, else we fix a line, we see that the first and third case clearly fix a line on $V_{\min}$, leaving us with the second case. In this case the trace of $t$ on $V_{\min}$ is $-6$, so the trace of $t$ on $L(G)$ is $20$, whence $H$ has at least twelve trivial composition factors, so we need at least thirty-six $3$-dimensional factors to avoid fixing a line, but this is more than $52=\dim(L(G))$, so that $H$ fixes a line on $L(G)$.

For $E_6$ the traces of the elements on $V_{\min}$ are $1$ more than for $F_4$, so we end up with the composition factors above, or no trivial composition factors at all, but this is impossible. The first and third case fall as before, and in the second case we switch to $L(G)$, on which $t$ has trace $14$. Since the trace of $x$ on $V_{\min}$ is $0$, on $L(G)$ it is $-3$ or $6$. This gives composition factors of $H$ on $L(G)$ as either $5^{11},3^5,1^8$ or $5^8,3^8,1^{14}$, and we again see that $H$ fixes a line on $L(G)$.

For $G=E_7$ the trace of an involution is $\pm8$ and the trace of an element of order $3$ is one of $-25$, $-7$, $2$ and $20$. This yields four possible sets of composition factors for $V_{\min}\downarrow_H$, namely
\[ 3^{12},1^{20},\qquad 5^2,3^{14},1^4,\qquad 5^9,3^3,1^2,\qquad 5^6,3^6,1^8.\]
Of these, only the second could potentially not fix a line on $V_{\min}$, and in this case the trace of an involution on $L(G)$ must be $5$, and the trace of an element of order $3$ on $L(G)$ is either $-2$ or $7$, whence the composition factors of $L(G)\downarrow_H$ are
\[ 5^{13},3^{19},1^{11}\qquad \text{or}\qquad 5^{10},3^{22},1^{17}.\]
In either case, $H$ clearly fixes a line on $L(G)$.

The remaining possibility is that $H\cong 2\cdot \Alt(5)$ embeds inside $G=E_7$ with $Z(H)=Z(G)$. The faithful simple modules for $H$ are of dimension $2$ and $4$, on which $x$ acts with trace $-1$ and $1$ respectively. Since the trace of $x$ on $V_{\min}$ is one of $-25,-7,2,20$, we cannot have the last case, but the rest yield
\[ 4,2^{26},\qquad 4^7,2^{14},\qquad 4^{10},2^8.\]
As with the case of $\Alt(5)$, only a projective module can have a $2$-dimensional composition factor but no $2$-dimensional submodule or quotient, and since $P(4)$ has structure $4/2/4$, we need at least twice as many $4$s as $2$s not to fix a $2$-space on $V_{\min}$. This clearly never happens, and so we are never maximal for $p=5$ and $G$ of type $E_7$.

For $E_8$ the trace of an involution on $L(G)$ is $24$ or $-8$ and the trace of an element of order $3$ on $L(G)$ is one of $-4$, $5$, $14$ and $77$. This leads to seven possible sets of composition factors (since we cannot simultaneously have traces of $-8$ and $77$) and we have
\[ 5^{28},3^{28},1^{24},\quad 5^{25},3^{31},1^{30},\quad 5^{22},3^{34},1^{36},\quad 5^{20},3^{44},1^{16},\quad 5^{17},3^{47},1^{22},\quad 5^{14},3^{50},1^{28},\quad 5,3^{55},1^{78}.\]
Although the fourth possibility gets close to the limit, it still must fix a line, and hence $H$ always fixes a line on $L(G)$, as needed.
\end{pf}

\begin{prop} Suppose that $p=3$. If $G=F_4,E_6,E_7,E_8$ and $H\cong \Alt(5)$ or $G=E_7$ and $H\cong 2\cdot \Alt(5)$ with $Z(H)=Z(G)$, then $H$ fixes a line on the simple module $L(G)$.
\end{prop}
\begin{pf} Proposition \ref{prop:ext1simples} and Lemma \ref{lem:covertrivialBrauer} can be used (together with the fact that the Sylow $3$-subgroup is cyclic) to see that the only indecomposable module with a trivial composition factor but no trivial submodule or quotient is the projective module $P(4)$, which has composition factors $4^2,1$. Thus if $M$ is any module for $H$ in characteristic $3$ with fewer than twice as many $4$s as $1$s, then $M$ has a trivial submodule or quotient. We can now consult the tables in \cite{litterickmemoir} to see whether, for $G=F_4,E_6,E_7,E_8$, there are any such sets of composition factors. We see that $H$ always fixes a line on $L(G)$ for $F_4$, $E_6$ (remember to remove a single trivial as $L(G)$ is not simple) and $E_7$ (and in most cases fixes a line on $V_{\min}$ as well), and always fixes a line on $L(G)$ for $G=E_8$. We must also consider $2\cdot \Alt(5)$ inside $E_7$ with $Z(H)=Z(G)$, and here we fix a line on $L(G)$ also. 
\end{pf}

\begin{prop} Suppose that $p=2$.
\begin{enumerate}
\item If $G=F_4,E_6,E_7$ and $H\cong \Alt(5)$ then $H$ fixes either a line or a $2$-space on $\Vmin$, and hence is not maximal.

\item If $G=F_4,E_6,E_7$ and $H\cong \Alt(5)$ with $\Vmin\downarrow_H^\phi=\Vmin\downarrow_H^*$ for $\phi$ a generator of $\Out(H)$, then $H$ fixes a line on $\Vmin$.

\item If $G=E_8$ and $H\cong \Alt(5)$ then $H$ fixes a line on $L(G)$.
\end{enumerate}

\end{prop}
\begin{pf} The modules for $H$ have dimension $1$, $2$, $2$ and $4$, and there are semisimple elements $x$ of order $3$ and $y$ of order $5$. The two modules of dimension $2$ are swapped by the outer automorphism of $H$.

If $V$ is a module for $H$ whose Brauer character is $\Out(G)$-stable, then $V$ has factors $1$, $4$ and $2_1\oplus 2_2=4_2$. Each of these has a rational trace for an element of order $5$, and there is a unique rational class of semisimple elements of order $5$ in each of $F_4$, $E_6$ and $E_7$. These have traces $1$, $2$ and $6$ respectively on $\Vmin$, and together with the possible traces of elements of order $3$, yields very few possible sets of composition factors for $\Vmin\downarrow_H$, namely
\[ 4^4,4_2,1^7,\qquad\text{and}\qquad 4,4_2^4,1^7\]
for $E_6$ and these with one trivial removed for $F_4$, and
\[ 4^8,4_2^2,1^{16},\qquad\text{and}\qquad4_2^8,4^2,1^{16}\]
for $E_7$. If $G=E_7$ then these have non-positive pressure and so $H$ always fixes a line on $\Vmin$. If $G=F_4,E_6$ then we have six and seven trivial modules respectively, and from \cite{lawther1995} we see that an involution acts with at least two and three blocks of size $1$ respectively, so from Lemma \ref{lem:numberof2sneeded} we see that we need at least ten and twelve composition factors of dimension $2$ not to fix a line on $\Vmin$, and we have at most eight. This completes the proof of (ii).

For (i), this is easy: as $4$ is projective it breaks off, and so the only way that $H$ does not stabilize a $1$- or $2$-space on a module $V$ is if $V\cong 4^{\oplus m}$ for some $m$. But then the trace of an element of order $5$ on $V$ is $-m<0$, so this is not true for $\Vmin$ as we stated above the trace of $u=(1,2,3,4,5)$ if it is rational.

\medskip

We thus move to $G=E_8$. There are four and fourteen classes of elements of order $3$ and $5$ respectively, and this yields many possibilities for the composition factors of $H$ on $L(G)$. Using a computer we can list them all, and by Lemma \ref{lem:numberof2sneeded} we firstly need that there are more $2$-dimensional factors than trivials, else we certainly fix a line.

By \cite[Table 9]{lawther1995}, we see that $t$ in $H$ acts with at least eight blocks of size $1$, and so if there are $8+2b$ trivial composition factors, we use Lemma \ref{lem:numberof2sneeded} to see that we need at least $16+3b$ composition factors of dimension $2$. We now check in each case that there are fewer $2$-dimensional factors than this, proving that $H$ always fixes a line on $L(G)$, as needed.
\end{pf}

For the next section, we want to understand copies of $H\cong\Alt(5)$ inside $G=E_7$ in characteristic $3$ with specific composition factors on the minimal module.

\begin{lem}\label{lem:specificalt5} Let $p=3$. If $H\cong \Alt(5)$ is a subgroup of $G=E_7$ with composition factors $4^6,3_1^6,3_2^4,1^2$ on $\Vmin$, then $H$ is contained in an $A_6$-parabolic subgroup of $G$ and and the action of $H$ on $L(G)$ has $P(1)^{\oplus 3}$ as a summand.
\end{lem}
\begin{pf} Let $H$ be embedded in $G$ with those factors. The action of $x=(1,2,3)$ on $\Vmin$ must have at least sixteen blocks of size $3$, whence it is $3^{18},1^2$ by \cite[Table 7]{lawther1995}, and the only way this can work is if $H$ embeds as
\[ (4/1/4)^{\oplus 2}\oplus 3_1^{\oplus 6}\oplus 3_2^{\oplus 4}\oplus 4^{\oplus 2}.\]

Since we know that $H$ lies in a positive-dimensional subgroup, we run through the list of maximal connected subgroups from \cite{liebeckseitz2004} and see how many embeddings we can find.

Suppose that $H$ is contained in the $A_6$-parabolic subgroup $X$ of $G$: this acts on $\Vmin$ with composition factors $7,7^*,21,21^*$, where $7$ is the natural module and $21$ its exterior square, and $7$ is a submodule. Since $7$ is a submodule of $\Vmin\downarrow_X$, $H$ cannot fix a line on the natural module, and so the only embedding can be as $4\oplus 3_1$ (or $4\oplus 3_2$). The action of $X$ on $L(G)$ has $7\otimes 7^*$, $\Lambda^3(7)$ and $\Lambda^3(7^*)$ as factors, and it is easy to check that all three of these have a copy of $1/4/1=P(1)$ as a summand, so the result holds if $H$ is contained in an $A_6$-parabolic.

\medskip

We now check the other parabolics: since $H$ does not fix a line on $\Vmin$, it does not lie in the $E_6$ parabolic. The $D_5A_1$ parabolic has two $2$-dimensional composition factors on $\Vmin$, so not compatible with the composition factors of $H$. For the $A_1A_5$ parabolic $X$, notice that both $(0,\lambda_1)$ and $(1,\lambda_1)$ are composition factors of $\Vmin\downarrow_X$, so in fact $H$ is contained in the $A_5$ parabolic, which is in the $A_6$ parabolic, so done. If $H$ lies inside the $A_2A_4$ parabolic then, as it cannot act irreducibly on a $5$-dimensional module it must fix a line or hyperplane, hence lie inside $A_2A_3$, which is contained inside the $A_6$ parabolic again. A similar argument kills off the $A_1A_2A_3$ parabolic.

If $X$ is the $D_6$ parabolic, then $H$ cannot fix a line on the natural module, since the action of $X$ on $\vmin$ is $12/32/12$. Computing with the involutions and the thirty-four classes of elements of order $5$, one sees that $3_1^{\oplus 3}\oplus 3_2$ is the only possibility for the restriction of the $12$ to $H$. We can therefore apply Lemma \ref{lem:singularsubspace} to find that the image of $H$ in $D_6$ lies inside a parabolic subgroup, hence $H$ lies inside a different parabolic subgroup of $G$, but all other parabolics have been dealt with.

\medskip

Having proved the result for parabolics, we move on to reductive maximal subgroups, starting with the maximal-rank subgroups.

If $X$ is the $A_7$ maximal-rank subgroup then $H$ embeds as either $\Alt(5)$ or $2\cdot\Alt(5)$ in the natural module for $\SL_8$. In the former case, if $H$ fixes a line or hyperplane if lies inside an $A_6$-parabolic, considered before, and the only other case if $4\oplus 4$, but $\Lambda^2(4\oplus 4)$ has a trivial submodule, and this is a summand of $\Vmin\downarrow_H$. If $2\cdot \Alt(5)$ embeds in $\SL_8$ then, as the faithful simple modules for $2\cdot\Alt(5)$ have dimension $2$ and $6$, it either has a $2$-dimensional submodule or quotient, whence its exterior square has a trivial submodule or quotient, not permitted.

For $X=A_2A_5$ it must act irreducibly as $3_1$ along the $A_2$ factor and cannot fix a line along the $A_5$ factor, else it would be inside the $A_2A_4$ parabolic. This means it must act as $3_i\oplus 3_j$ along the $A_4$ factor: this however means that either $i=j=2$ or $H$ fixes a line on $(10,1000)$, which is a submodule of $\vmin\downarrow_X$. If $i=j=2$ though, we have that $4/1/4$ appears four times in $\vmin\downarrow_H$, not allowed either.

For $X=A_1D_6$, if $H$ embeds in the $D_6$ factor we are done, whence $H$ acts irreducibly along the $A_1$ factor, say as $2_1$ on its natural module. Since the summands of $\Vmin\downarrow_X$ are the spin module for the $D_6$ factor and the product of the two naturals, this means that $H$ embeds as $2\cdot \Alt(5)$ when acting on the natural for the $D_6$, with composition factors $2_1$, $2_2$ and $6$. Since $x=(1,2,3)$ acts as $3^{18},1^2$ on $\Vmin$, it must act projectively on the $24$-dimensional summand of $\Vmin\downarrow_X$, hence its projection along $D_6$ acts projectively on the $12$. Thus the projection of $H$ along $D_6$ acts with summands from the list \[2_1/2_2/2_1,\qquad 2_2/2_1/2_2,\qquad 6,\]
i.e., the faithful projectives for $2\cdot \Alt(5)$. The first is not possible as $H$ would then fix a line on $\Vmin$ as $2_1^{\otimes 2}$ would be a submodule. If $H$ embeds as the sum of two isomorphic $6$-spaces then since $\F_9$ must be a subfield of our field of definition, else $\Vmin\downarrow_H$ is not definable, $-1$ is a square in our field, so Lemma \ref{lem:singularsubspace} applies and the projection of $H$ stabilizes a totally isotropic $6$-space, thus lies inside one of the two $A_5$-parabolics of $D_6$; this means that $H$ is contained inside $A_7$ or $A_2A_5$, examined before. Thus the projection of $H$ along $D_6$ acts on the natural module as $2_2/2_1/2_2\oplus 6$, but this does not fix an orthogonal form, compleing the result for this group.

Now for the other maximal connected reductive subgroups: for $A_1F_4$, since $Z(F_4)=1$ $H$ cannot act along the $A_1$ factor and so must lie inside the $F_4$ factor, but this is contained in $E_6$, not allowed. For $A_1G_2$ we again get that $H$ is contained in the $G_2$, and this time $G_2\leq A_6$, as was seen in \cite[pp.33--35]{seitz1991}, so this case is also done.

Finally, for $H\leq X=G_2C_3$, then the summands of $\Vmin\downarrow_X$ are $(10,100)$ and $(00,001)$, and one has that the exterior cube of $L(100)$ is the sum of $L(100)\oplus L(001)$. To embed $H$ into $C_3$ we need to choose a $6$-dimensional module: $3_1\oplus 3_2$, $1/4/1$, $4\oplus 1^{\oplus 2}$ and $3_i\oplus 1^{\oplus 3}$ do not support symplectic forms (the quickest way to see this is to note that they are not summands of their exterior cubes) and $3_1^{\oplus 2}$ has four trivial submodules on its exterior cube, so not correct. Thus $H$ does not embed into $X$, completing the proof.
\end{pf}

\section{$\Alt(6)$}

The cases under consideration are $p=2,3,5$ and primes larger than $5$. We continue our assumption that $E_6$ and $E_7$ are used to denote the simply connected forms, so that $Z(E_6)$ and $Z(E_7)$ have orders $3$ and $2$ respectively.

For $p\neq 2,3,5$, in Theorem \ref{thm:alt6} we say nothing, so we will not give a formal proposition here, and just list those entries from \cite{litterickmemoir} that are fixed-point free on both $\Vmin$ and $L(G)$.

For $F_4$ there is one possibility that fixes no line on either $V_{\min}$ or $L(G)$, and that acts as $9^2,8_1$ and $10^2,8_1,8_2^3$ on $V_{\min}$ and $L(G)$ respectively. 

For $G=E_6$, there is a single set of factors for $\Alt(6)$s that fix lines on neither module, and these have factors $9,8_1,5_1,5_2$ and $10^2,9^2,8_1^2,8_2^3$ on $V_{\min}$ and $L(G)$ respectively.

For $3\cdot \Alt(6)$ embedding into $E_6$ with centres coinciding, there is only one possibility whose image in $L(G)$ does not fix a line, and this acts as $9^2,6,3$ on $\Vmin$, so stabilizes a $3$-space on $\Vmin$, hence not maximal in $G$ by Lemma \ref{lem:smallspacestabilizers}, although we don't consider almost simple groups with socle $\Alt(6)$ inside almost simple groups with socle $E_6$.

For $E_7$ there is a single example, acting as $10^4,8_1^2$ on $V_{\min}$ and as $10^2,9^3,8_1^4,8_2^3,5_1^3,5_2^3$ on $L(G)$, which exists inside the $A_7$ subgroup, but there might be other classes.

For $2\cdot \Alt(6)$ inside $E_7$ with the centres coinciding, there is (up to outer automorphism) a single possible action on $V_{\min}$ of $8_3^2,10_2^3,10_3$, acting on $L(G)$ as $10^5,9^3,8_1^4,8_2^3$.

For $\Alt(6)$ inside $E_8$, there are four potential collections of composition factors for $L(G)\downarrow_H$ with no fixed points, namely (up to an outer automorphism that swaps $5_1$ and $5_2$)
\[10^{10},9^6,8_1^4,8_2^4,5_1^6\qquad 10^{10},9^6,8_1^4,8_2^4,5_1^3,5_2^3,\qquad 10^8,9^8,8_1^6,8_2^6,\qquad 10^8,9^4,8_1^7,8_2^7,5_1^2,5_2^2.\]

Since we prove specific things about other primes, we now switch to formal propositions, as with the previous and later sections.

\begin{prop} Suppose that $p=5$.
\begin{enumerate}
\item If $G=F_4,E_6$ and $H\cong \Alt(6)$ then $H$ fixes a line on $\Vmin$, and hence is not maximal.
\item If $G=E_6$ and $H\cong 3\cdot\Alt(6)$ with $Z(H)=Z(G)$ then either $H$ fixes a line on $L(G)$ or $\Vmin\downarrow_H$ has exactly one of a semisimple submodule of dimension $12$ and a semisimple quotient of dimension $12$, and in particular satisfies the conditions for non-maximality of Lemma \ref{lem:novelty}(iii) in an almost simple group with socle $E_6$.
\item If $G=E_7$ and $H\cong\Alt(6)$ then either $H$ fixes a line on $L(G)$ or $\Vmin\downarrow_H$ and $L(G)\downarrow_H$ are
\[ 10^{\oplus 4}\oplus 8^{\oplus 2}\qquad\text{and}\qquad 10^{\oplus 2}\oplus 5_1^{\oplus 3}\oplus 5_2^{\oplus 3}\oplus P(8)^{\oplus 3}\oplus 8.\]
\item If $G=E_7$ and $H\cong 2\cdot\Alt(6)$ with $Z(H)=Z(G)$ then either $H$ fixes a line on $L(G)$ or $\Vmin\downarrow_H$ and $L(G)\downarrow_H$ are
\[ 10_2^{\oplus 3}\oplus 10_3\oplus 4_1/4_2\oplus 4_2/4_1\qquad\text{and}\qquad 10_1^{\oplus 5}\oplus P(8)^{\oplus 3}\oplus 8.\]
\item If $G=E_8$ and $H\cong \Alt(6)$ then $H$ fixes a line on $L(G)$.
\end{enumerate}
\end{prop}
\begin{pf} From the Brauer tree in Section \ref{sec:modulesalt}, and Lemma \ref{lem:covertrivialBrauer}, the appropriate module is projective cover of $8$, which is $8/1,8/8$, and so we need at least three times as many $8$s as $1$s. We can go through the tables in \cite{litterickmemoir} to see if this is the case. For $F_4$ we are therefore done using $V_{\min}$ or $L(G)$, and for $E_6$ we are done using $V_{\min}$. If $H\cong3\cdot \Alt(6)$ however, we can have the action of $H$ on $L(G)$ having factors $10^2,8^7,1^2$, with all other cases fixing lines on $L(G)$. In this case we need $P(8)^{\oplus 2}$ in $L(G)\downarrow_H$, and the $10^2$ are projective, leaving $8$ remaining. This yields an action of $H$ on $L(G)$ of $10^{\oplus 2}\oplus P(8)^{\oplus 2}\oplus 8$, on which the element $y=(1,2,3,4,5)$ acts with blocks $5^{15},3$, so class $A_4+A_1$. This acts on $V_{\min}$ as $5^5,2$. The composition factors of $V_{\min}\downarrow_H$ are one of
\[ 15,6^2,\qquad 6^3,3^3,\]
with neither of these modules definable over $\F_5$, so that $H\leq E_6(5^{2a})$ for some $a$. In this case, if there is a $3$ in the socle of $\Vmin\downarrow_H$ then $H$ is contained in a member of $\ms X^\sigma$ by Lemma \ref{lem:smallspacestabilizers}. The Brauer tree of the block containing $3$ has two edges with exceptional vertex at the end labelled by $3$, so the projectives are
\[ 6/3/6,\qquad 3/3,6/3,\]
and therefore we need twice as many $6$s as $3$s in order not to stabilize a $3$-space; thus $6^3,3^3$ does indeed stabilize a $3$-space. Since the $15$ is projective it splits off from $V_{\min}\downarrow_H$, and so we must have $15\oplus 6\oplus 6$, and so $u$ cannot act on $V_{\min}$ as $5^5,2$, a contradiction.

To deal with almost simple groups of type $E_6$, we need to examine the $6^3,3^3$ case more closely. Since the action of $u$ is $5^5,2$, there is a unique non-projective summand of $\Vmin\downarrow_H$, necessarily of dimension $12$. This must be one of $3/3,6$ or $3,6/3$, with the remaining summand being $P(6)$. Therefore $\Vmin\downarrow_H$ has $6^{\oplus 2}$ as a submodule or quotient, but not both, and so since the graph automorphism interchanges $\Vmin$ and its dual, we may apply Lemma \ref{lem:novelty}(iii) to get no maximal subgroups in the almost simple group case.

For $E_7$ there is a single possibility for the composition factors of $H$ on $L(G)$ and $V_{\min}$ that has at least three times as many $8$s as $1$s on $L(G)$, and this is $H$ acting with factors $10^4,8^2$ for $V_{\min}$ and $10^2,8^{10},5_1^3,5_2^3,1^3$ for $L(G)$.

The only way this could work without fixing a line or hyperplane on $L(G)$ is
\[ 10^{\oplus 2}\oplus 5_1^{\oplus 3}\oplus 5_2^{\oplus 3}\oplus P(8)^{\oplus 3}\oplus 8,\]
on which $u=(1,2,3,4,5)$ acts as $5^{26},3$. This is consistent with coming from class $A_4+A_2$.

For $H\cong 2\cdot \Alt(6)$ in $G=E_7$ with the centres coinciding, the unique action on $L(G)$ with at least three times as many $8$s as $1$s is similar to above, namely $10^5,8^{10},1^3$, so again has $u$ coming from class $A_4+A_2$. This must act on $V_{\min}$ as 
\[10_2^{\oplus 3}\oplus 10_3\oplus 4_1/4_2\oplus 4_2/4_1,\]
and everything is again consistent.

In $E_8$, we are left with two sets of composition factors on $L(G)$, after examining the table in \cite{litterickmemoir}. These are
\[ 10^8,8^{18},5_1^2,5_2^2,1^4,\qquad 10^9,8^{16},5_1^4,5_2,1^5.\]
For the first case, we need $P(8)^{\oplus 4}$ to hide the four trivial modules, and then have six $8$-dimensional modules left over. Taking projectives into account, if the $8$s form simply an $8^{\oplus 6}$, we would get $46$ blocks of size $5$ in the action of $u$ on $L(G)$. Since the only modules constructible using only $8$s are $8$ and $8/8$, on which $u$ acts as $5,3$ and $5^3,1$ respectively, we get the action of $u$ to be one of 
\[ 5^{46},3^6,\qquad 5^{47},3^4,1,\qquad 5^{48},3^2,1^2,\qquad 5^{49},1^3.\]
None of these appears on Table \ref{tab:unipe8p5}, and so $H$ must fix a line on $L(G)$.

For the second one, the only way it could not have a trivial submodule is $P(8)^{\oplus 5}\oplus 8$, which means that $u$ acts on $L(G)$ with Jordan blocks $5^{49},3$, which does not appear on Table \ref{tab:unipe8p5}. This proves that $H$ always fixes a line.
\end{pf}

For $p=3$, Theorem \ref{thm:alt6} gave no information about $F_4$ and $E_8$: for $F_4$, this will be considered in a later paper of the author on $\SL_2(q)$ subgroups of exceptional groups not of type $E_8$. For $E_8$ there are many possible sets of composition factors for the action of $H$ on $L(G)$ and do not consider this here.

\begin{prop}Suppose that $p=3$.
\begin{enumerate}
\item If $G=E_6$ and $H\cong \Alt(6)$ then $H$ fixes a line on either $\Vmin$ or $L(G)$.
\item If $G=E_7$ and $H\cong \Alt(6)$ then $H$ fixes a line on $\Vmin$ or $L(G)$.
\item If $G=E_7$ and $H\cong 2\cdot \Alt(6)$ with $Z(H)=Z(G)$ then $H$ always stabilizes a $\sigma$-stable $2$-space, and if $\phi$ is any automorphism of $H$ that stabilizes $\Vmin\downarrow_H$ then $\langle H,\phi\rangle$ also stabilizes a $2$-space on $\Vmin$.
\end{enumerate}
\end{prop}
\begin{pf} From Section \ref{sec:modulesalt}, we know that there are five simple modules: $1$, $3_1$, $3_2$, $4$ and the projective $9$. Only the $4$ has any $1$-cohomology, and it is $2$-dimensional. The structure of $P(4)$ is
\[ 4/1,1,3_1,3_2/4,4,4/1,1,3_1,3_2/4,\]
and so has five $4$s for every four trivials. If an indecomposable module is not projective and has no trivial submodules or quotients, then by removing the $3_i$s from top and bottom we may assume that it has three socle layers, so in particular needs as many $4$s as $1$s. Thus if $M$ is a module with more trivial factors then $4$-dimensional factors then $H$ fixes a line or hyperplane on $M$.

Let $G=E_6$. If $H\cong\Alt(6)$ acts on $\vmin$ with more trivial composition factors than $4$-dimensional ones then it fixes a line by the previous paragraph, and so we may assume the contrary. Using the traces of semisimple elements, one finds only three possible sets of composition factors on $\vmin$, up to swapping $3_1$ and $3_2$, and these are
\[ 9,4^3,3_1,1^3,\qquad 9,3_1^3,3_2^3,\qquad 4^3,3_1^3,3_2^2.\]
The restriction of the last case to the point stabilizer $\Alt(5)$ must be semisimple as the $3_i$s become projective and the $4$s have no self extensions, so that $x=(1,2,3)$ acts on $\vmin$ with Jordan blocks $3^{17},1^5$, not present in \cite[Table 5]{lawther1995}, and so this case cannot occur.

For the middle case we turn to the Lie algebra, where the traces of the semisimple elements on $\vmin$ means that $t=(1,2)(3,4)$ and $v=(1,2,3,4)(5,6)$ act with trace $13$, and $u=(1,2,3,4,5)$ acts with trace $2$ on the $77$-dimensional module $L(G)$. This is enough to uniquely determine the composition factors of $L(G)\downarrow_H$, and they are
\[ 9^6,4^2,3_1,3_2,1^9,\]
so that $H$ centralizes at least a $7$-space on $L(G)$.

Finally we have the first case. However, we only have three $4$s here, so we cannot have two $4$s both above and below the three trivials, so this case must also fix a line. (If there are an odd number $2n-1$ of trivial composition factors then we need at least $2n$ $4$-dimensional factors to avoid fixing a line, in a slight improvement to our earlier result.) Thus $H$ always fixes a line on $\vmin$ or $L(G)$, as needed.

If $G=E_7$ then, as we said before, we are only interested in possible sets of composition factors on $V_{\min}$ that have at least as many $4$s as $1$s. Using the traces of all semisimple elements, it turns out that, up to application of an outer automorphism swapping $3_1$ and $3_2$, there are five such sets of composition factors, namely
\[ 4^{10},3_1^2,1^{10},\qquad 9,4^9,3_1,1^8,\qquad 4^6,3_1^9,3_2,1^2,\qquad 4^6,3_1^6,3_2^4,1^2,\qquad 9,4^5,3_1^5,3_2^4.\]
The first three of these yield unique sets of composition factors on $L(G)$, namely
\[ 9^6,4^9,3_1^6,3_2^5,1^{10},\qquad 9^5,4^{10},3_1^7,3_2^5,1^{12},\qquad 4^{16},3_1^{10},3_2^6,1^{21}.\]
Notice that each of these has more trivials than $4$s, so the first three cases fix a line on $L(G)$. 
For the final case, we claim it cannot occur. To see this, restrict to the point stabilizer $\Alt(5)$ inside $H$: the $3_i$s become projective and the $4$s have no self-extensions, with the $9$s restricting to the projective $P(4)$, and so this module becomes
\[ P(4)\oplus 3_1^{\oplus 5}\oplus 3_2^{\oplus 4}\oplus 4^{\oplus 5},\]
on which $x=(1,2,3)$ acts with Jordan blocks $3^{17},1^5$. A quick check of \cite[Table 7]{lawther1995} proves that this is not a valid unipotent class.

We are left with the fourth case: restricting to an $\Alt(5)$ subgroup $L$, we notice that it has composition factors $4^6,3_1^6,3_2^4,1^2$, and so Lemma \ref{lem:specificalt5} applies. Thus $L(G)\downarrow_L$ has three copies of $P(1)$ as summands. We claim that, if $V$ is a module for $H$ with no trivial submodules, and its restriction to $L$ has $P(1)^{\oplus m}$ as a summand if and only if $V$ has $P(4)^{\oplus m}$ as a summand. To see this, it suffices to prove the case $m=1$ by induction, and notice that the restriction of $P(4)$ for $H$ to $L$ is $P(4)^{\oplus 2}\oplus P(1)$, so one direction is true. For the other, since $V$ has no trivial submodules, in order to restrict to $1/4/1$ there must be five socle layers to $V$, but $P(4)$ only has five socle layers, so the whole projective must be a submodule of $V$, hence a summand.

Thus $L(G)\downarrow_H$ has $P(4)^{\oplus 3}$ as a summand, so since $P(4)$ has structure
\[ 4/1,1,3_1,3_2/4,4,4/1,1,3_1,3_2/4,\]
there are at least six copies of $3_2$ in $L(G)\downarrow_H$. However, using the Brauer character of $\Vmin\downarrow_H$, we get two options for the composition factors of $L(G)\downarrow_H$, and these are
\[ 9^3,4^{15},3_1^6,3_2^5,1^{13}\qquad\text{and}\qquad 9^5,4^{11},3_1^6,3_2^5,1^{11},\]
neither of which has six $3_2$s. Thus $H$ fixes a line on $L(G)$, as required.

\medskip

Of course, we also have to deal with $H\cong2\cdot \Alt(6)=\SL_2(9)$ embedding in $G=E_7$ with the centres coinciding. We firstly assume that $G=G(q)$ for $q$ an odd power of $2$. In this case the irreducible modules are $4=2_1\oplus 2_2$ and $12=6_1\oplus 6_2$. An element of order $8$ acts with trace $0$ on both of these, and so this element acts with trace $0$ on all of $\vmin$. However, no such element exists in $E_7$, and so $H$ does not embed with $Z(H)=Z(G)$ unless $q$ is an even power of $2$. In particular, this means that we may assume that there is no $2$-dimensional submodule of $\vmin\downarrow_H$. Hence $\vmin\downarrow_H$ is a submodule of copies of $P(6_1)$ and $P(6_2)$. The structure of $P(6_1)$ is
\[ 6_1/2_1/2_2/2_1/6_1,\]
and so we actually see that $\vmin\downarrow_H$ either has a $2$-dimensional submodule or quotient, or is actually projective. Since $56$ is not a multiple of $3$, this is impossible, so $H$ always fixes a $2$-dimensional subspace of $\vmin$, as needed.
\end{pf}

For $p=2$ we get a result for everything but $E_8$, where there is one set of composition factors that might yield an example that does not fix a line on $L(G)$.

\begin{prop} Suppose that $p=2$.
\begin{enumerate}
\item If $G=F_4$ and $H\cong \Alt(6)$ then either $H$ or its image under the graph automorphism fixes a line on $\Vmin$, and hence is not maximal.
\item If $G=E_6$ and $H\cong \Alt(6)$ then $H$ fixes a line on $\Vmin$, and hence is not maximal.
\item If $G=E_6$ and $H\cong 3\cdot \Alt(6)$ then $H$ fixes a line on $L(G)$ and hence is not maximal.
\item If $G=E_7$ and $H\cong \Alt(6)$ then $H$ fixes a line on $\Vmin$.
\item If $G=E_8$ and $H\cong \Alt(6)$ then either $H$ fixes a line on $L(G)$ or the composition factors of $L(G)\downarrow_H$ are
\[8_1^6,8_2^6,4_1^{16},4_2^{16},1^{24}.\]
\end{enumerate}
\end{prop}
\begin{pf} From Section \ref{sec:modulesalt} we see that there are five simple modules: $1$, $4_1$, $4_2$, and two projectives $8_1$ and $8_2$. The projective cover of $4_i$ is given by
\[ P(4_i)=4_i/1/4_{3-i}/1/4_i/1/4_{3-i}/1/4_i,\]
so we see that for $H$ not to fix a line on a module we must have at least five modules of dimension $4$ for each four of dimension $1$.

If, in a module $V$, there are the same number of $8_1$s as $8_2$s then $u=(1,2,3,4,5)$ acts with rational trace, namely $+1$ on $1$ and $16=8_1\oplus 8_2$, and $-1$ on the $4_i$. There is a single rational class of elements of order $5$ in $F_4$, $E_6$ and $E_7$, with trace on $\Vmin$ given by $1$, $2$ and $6$ respectively. The last thing we need to note is that the minimal module for $F_4$ and $E_6$ restricts to $H$ with at least two and three trivial composition factors, since there are at least two and three Jordan blocks of size $1$ in the action of $t=(1,2)(3,4)$ on $\Vmin$ by \cite{lawther1995}, and $t$ has no blocks of size $1$ on all non-trivial simple modules. These facts prove that $\Vmin\downarrow_H$ always has at least as many trivials as $4$-dimensionals, hence has non-positive pressure and fix a line or hyperplane on $\Vmin$, as needed.

Thus we may assume that there are more $8_1$s than $8_2$s in $\Vmin\downarrow_H$, and in particular that there is at least one $8_1$. If $G=F_4$ then the $8_1$, together with the two trivials, and three $4$s needed to avoid fixing a line on $\Vmin$, gives us $22$ of $26$ dimensions, so our composition factors must be $8,4^4,1^2$. The trace of $x=(1,2,3)$ on this module is one of $-7,-4,-1,2,5$ depending on the numbers of $4_1$s and $4_2$s, and $-1$ is the only one of these that is a trace of a class in $F_4$, so the composition factors are $8_1,4_1^2,4_2^2,1^2$. The Brauer character of its image under the graph automorphism is easy to compute, and is the character of $8_2^3,1^2$, which clearly fixes a line on $\Vmin$, as needed.

\medskip

For $E_6$, we have three trivial composition factors, and hence four $4$s else we fix a line, and an $8_1$, so we have one more trivial, but the same composition factors as before, namely
\[ 8_1,4_1^2,4_2^2,1^3.\]
Since this has pressure $1$, and the $8_1$ is projective so breaks off, by the discussion after Proposition \ref{prop:lowpressuremodules} if we ignore the $8_1$ the module must be uniserial. However, the element $t=(1,2)(3,4)$ acts projectively on a module of the form $1/4_i/1$, and so $t$ would have to act on $\Vmin$ as $2^{13},1$, not a valid action by \cite[Table 5]{lawther1995}. Thus $H$ fixes a line on $\Vmin$, as needed.

\medskip

If $H\cong 3\cdot \Alt(6)$ with $Z(H)=Z(G)$ then the composition factors for $\Vmin\downarrow_H$ consist of $3_1$, $3_2$ and $9$, with the $3_i$ being swapped by the outer automorphism of $H$ fixing $Z(H)$.

Using the traces of elements of orders $3$ and $5$, there are five possible sets of composition factors for $\Vmin\downarrow_H$, up to swapping the $3_i$, and these are
\[ 3_1^8,3_2,\quad 3_1^6,3_2^3,\quad 9,3_1^3,3_2^3,\quad 9^2,3_1^2,3_2,\]
and all of these fix $3$-spaces on $\Vmin$, so $H$ lies inside a positive-dimensional subgroup of $G$ by Lemma \ref{lem:smallspacestabilizers}. However, we need more because the graph automorphism does not stabilize $\Vmin$, so we switch to $L(G)$. Here the third and fourth sets of composition factors for $\Vmin\downarrow_H$ do not have corresponding factors on $L(G)$, so these do not exist, and the other two have factors
\[ 8_1^8,1^{14}\qquad\text{and}\qquad8,3_1^7,3_2^7,1^{14}\]
respectively, both of which have non-positive pressure and so $H$ fixes a line on $L(G)$.

\medskip

Now consider $G=E_7$. The element $v=(1,2,3,4)(5,6)$ acts projectively on all non-trivial simple modules, so we can get a lower bound on the number of trivial composition factors by examining the action of the unipotent classes on $\Vmin$ from \cite[Table 7]{lawther1995}. The maximum number of blocks of size $4$ in the action of $v$ is twelve, so there are at least eight trivials, hence either $H$ fixes a line on $\Vmin$ or there are nine $4$s, and also an $8_1$ from our assumption. This takes up $52$ of our $56$ dimensions, and so the composition factors of $\Vmin\downarrow_H$ must be $8_1,4^{10},1^8$, for some distribution of the $4$s among $4_1$ and $4_2$. Since $\Vmin\downarrow_H$ has pressure $2$, this means that it is a submodule of $8_1\oplus P(4_i)\oplus P(4_j)$ for some $i,j$. However, $P(4_i)$ is equal to
\[ 4_i/1/4_{3-i}/1/4_i/1/4_{3-i}/1/4_i,\]
which has dimension $24$, so $\Vmin$ is actually projective. But $v$ does not act projectively on $\Vmin$, a contradiction, so that $H$ fixes a line on $\Vmin$.

\medskip

We lastly turn to $E_8$. Using a computer and the traces of elements of orders $3$ and $5$, we are left with two possibilities (up to applying outer automorphisms), namely
\[ 8_1^6,8_2^6,4_1^{16},4_2^{16},1^{24},\qquad 8_1^5,8_2^5,4_1^{19},4_2^{16},1^{28}.\]
In these two cases there are at least fifty-five blocks of size $4$ in the action of $v=(1,2,3,4)(5,6)$ on $L(G)$. From Table \ref{tab:unipe8p4} we see that there are only two possible classes in $E_8$ in which $v$ can lie, namely $D_4(a_1)+A_2$ acting as $4^{56},3^8$, and $2A_3$ acting as $4^{60},2^4$. Notice that, in order to make a block of size $4$ in the action of $v$ we need a projective factor $P(4_i)$.

Notice that the second possibility for the composition factors of $L(G)\downarrow_H$ has factors $8^{10},4^{35},1^{28}$, so that there are exactly five $4$s for every four $1$s. Since this is the exact ratio as for $P(4_i)$, we see that $L(G)\downarrow_H$ consists solely of projective modules, not allowed by the action of $v$ above. Thus the second set of possible composition factors must fix a line.

For the first set of factors, if $v$ acted with eight blocks of size $3$ then we would need thirty-two $4$s to cover these, giving us no more $4$s for the rest of the $1$s, so this cannot occur. However, we are able to construct modules with action of $v$ given by $4^{60},2^4$, for example
\[ (8_1\oplus 8_2)^{\oplus 6}\oplus (P(4_1)\oplus P(4_2))^{\oplus 2}\oplus (4_1/1/4_2/1/4_1\oplus 4_2/1/4_1/1/4_2)^{\oplus 2}.\]

The restriction of this to $L=\Alt(5)$ has factors $1^{48},(2_1,2_2)^{34},4^{16}$, which cannot be distributed amongst the factors of $L(G)\downarrow_{E_7}$, so that $L$ does not lie in the $E_7$ parabolic. However, it can lie inside the $A_8$ maximal-rank subgroup and inside the $A_7$ parabolic, with composition factors $2_1,2_1,2_2,2_2$ on the natural module for the $A_7$.
\end{pf}

\bigskip

In order to be of use in the next section, just as in the last section, we include a lemma.  This is deducible from \cite{freyun}, but we produce a quick proof here.

\begin{lem}\label{lem:noalt6withspecificfactors} There is no copy of $H\cong\Alt(6)$ in $G=E_8$ with composition factors $10^{11},9^6,8_1^3,8_2^3,5_1^2,5_2^5,1$ on $L(G)$, in characteristic at least $7$ or $0$.
\end{lem}
\begin{pf} Firstly, $\Alt(6)$ does not embed, even projectively, in $A_1$, nor in $B_2$ or $G_2$. Secondly, as $H$ is a reductive subgroup of $G$, so is $C_G(H)$, and because $H$ has a trivial submodule on $L(G)$, $C_G(H)$ has dimension $1$, and in particular has an involution in it, hence we may apply Lemma \ref{lem:centinv} to get that we need only consider $E_7$ and $D_8$.

Since $E_7$ centralizes a $3$-space on $L(G)$, clearly $H$ cannot lie in it, so $H\leq D_8$.
\begin{center}
\begin{tabular}{cccc}
\hline Element&$\lambda_1$ & $\lambda_2$ & $\lambda_7$
\\\hline Possible& $7$ & $21$ & $-16$
\\traces for $x$&$-2$ & $3$ & $2$
\\\hline\multirow{3}{*}{\minitab[c]{Possible\\ traces for $y$}} & $-5$ & $15$ & $-1$
\\&$4$&$6$&$8$
\\ &$13$ & $78$ & $-64$
\\ \hline
\end{tabular}
\end{center}
We thus search for $16$-dimensional modules for $\Alt(6)$, or $2\cdot\Alt(6)$, with $x=(1,2,3)$ acting with trace $7$ or $-2$, and $y=(1,2,3)(4,5,6)$ acting with trace $-5$, $4$ or $13$.

Since $L(\lambda_2)=\Lambda^2(L(\lambda_1))$, it is easy to compute the traces of the elements of order $3$ in potential embeddings of $H$ into $D_8$.  Firstly, we cannot have more than two trivial factors on the $16$-dimensional because then we would get more than one trivial on $L(\lambda_2)$, not allowed. This leaves 
\[ 10,5_i,1,\qquad 9,5_i,1^2,\qquad 5_i^2,5_{3-i},1,\qquad 8_1^2,\qquad 8_1,8_2\]
The element $x$ acts on these with trace $4$ or $1$, $4$ or $1$, $4$ or $1$, $-2$ and $-2$ respectively, so only the one with $8$s works, but $y$ acts on these with trace $-2$ as well, not allowed.

For embedding $2\cdot \Alt(6)$ into $\Omega_{16}^+$, we have modules $4_1$, $4_2$, $8_1$ and $8_2$, on which $x$ acts with trace $1$, $-2$, $-1$ and $-1$ respectively, and $y$ acts with trace $-2$, $1$, $-1$ and $-1$ respectively. This means we need a different number of $4_1$s to $4_2$s, and which of the $8_i$ we use does not matter, so either $4_i^2,8_1$, $4_i^4$ or $4_i^3,4_{3-i}$, none of which works with the traces above. Hence $H$ cannot lie in $D_8$ either, and we are done.
\end{pf}

\section{$\Alt(7)$}

We consider the cases $p>7$, $p=7$, $p=5$, $p=3$ and $p=2$ in turn. We continue our assumption that $E_6$ and $E_7$ are used to denote the simply connected forms, so that $Z(E_6)$ and $Z(E_7)$ have orders $3$ and $2$ respectively.

\begin{prop} Suppose that $p\neq 2,3,5,7$.
\begin{enumerate}
\item If $G=F_4$ and $H\cong \Alt(7)$ then $H$ does not embed into $G$.
\item If $G=E_6$ and $H\cong \Alt(7)$, or $H\cong 3\cdot \Alt(7)$ with $Z(H)=Z(G)$, then $H$ fixes a line on $L(G)$ and hence is not maximal.
\item If $G=E_7$ and $H\cong \Alt(7)$ then $H$ fixes a line on $L(G)$ and hence is not maximal.
\item There is no embedding of $H\cong 2\cdot \Alt(7)$ into $G=E_7$ with $Z(H)=Z(G)$.
\item If $G=E_8$ and $H\cong \Alt(7)$ then either $H$ fixes a line on $L(G)$ and hence is not maximal, or $H$ acts on $L(G)$ with composition factors
\[ 35^4,15^4,14_1^2,10,10^*.\]
(Such a subgroup exists inside a $D_8$ maximal-rank subgroup acting as $15\oplus 1$ on the natural module.)
\end{enumerate}
\end{prop}
\begin{pf} This follows immediately from the tables in \cite{litterickmemoir} together with Lemma \ref{lem:smallspacestabilizers}.
\end{pf}

\begin{prop} Suppose that $p=7$.
\begin{enumerate}
\item If $G=F_4$ and $H\cong \Alt(7)$ then $H$ does not embed into $G$.
\item If $G=E_6$ and $H\cong \Alt(7)$, or $H\cong 3\cdot \Alt(7)$ with $Z(H)=Z(G)$, then $H$ fixes a line on $L(G)$ and hence is not maximal.
\item If $G=E_7$ and $H\cong \Alt(7)$ then $H$ fixes a line both on $\Vmin$ and on $L(G)$ and hence is not maximal.
\item There is no embedding of $H\cong 2\cdot \Alt(7)$ into $G=E_7$ with $Z(H)=Z(G)$.
\item If $G=E_8$ and $H\cong \Alt(7)$ then either $H$ fixes a line on $L(G)$ and hence is not maximal, or $H$ acts on $L(G)$ as
\[ 35^{\oplus 3}\oplus  21\oplus 14_1\oplus 14_2^2\oplus P(10)^{\oplus 2}\oplus 10\quad\text{or}\quad 35^{\oplus 4}\oplus 14_1^{\oplus 2}\oplus 10^{\oplus 6}\oplus 5^{\oplus 4}.\]
(The second of these is the reduction modulo $p$ of the case for $p=0$.)
\end{enumerate}
\end{prop}
\begin{pf} The Brauer tree in Section \ref{sec:modulesalt} implies that the projective modules are
\[ 1/5/1,\qquad 5/1,10/5,\qquad 10/5,10/10,\]
and the only indecomposable module with a trivial composition factor but no trivial submodule or quotient is $P(5)$ by Lemma \ref{lem:covertrivialBrauer}.

The action of $H$ on the adjoint module $L(G)$ has at least three trivials and at most two $5$-dimensionals, whether $H\cong\Alt(7)$ or $H\cong 3\cdot \Alt(7)$ inside $E_6$, and so $H$ always fixes a line on $L(G)$. For $G=E_7$, the unique set of composition factors has the same number of $5$s as $1$s, so fixes a line on $V_{\min}$ (it also does on $L(G)$).

For $G=E_8$ however, there are three potential sets of composition factors that have no trivial composition factor at all, and others with enough $5$s and $10$s to cover the trivials they do have.

Those without trivials are, up to projectives,
\[ 10^8,\qquad 10^7,5^2,\qquad 10^6,5^4.\]
In the case of $10^8$ we can only produce $10$s and $10/10$, on which $w=(1,2,3,4,5,6,7)$ acts as $7,3$ and $7^2,6$ respectively. Thus we have blocks of size $7$ and $6^i,3^{8-2i}$, not allowed by \cite[Table 9]{lawther1995}. Hence $H$ does not embed with these factors.

In the cases of $10^7,5^2$ and $10^6,5^4$, we have the following possible indecomposable modules:
\[ 10,\quad 5, \quad 10/10,\quad 10/5\oplus 5/10,\quad 10,5/10\oplus 5,10/10,\quad 5,10/5,10,\quad 10/5,10/10.\]
The actions of $w$ on these elements are (up to projectives) $3$, $5$, $6$, $1^2$, $4^2$, $2$ and $0$. For $10^7,5^2$ we are allowed only $P(10)^{\oplus 2}\oplus 10$, and for $10^6,5^4$ we are only allowed the semisimple case $10^{\oplus 6}\oplus 5^{\oplus 4}$, leading to the two cases listed in the proposition.

If the composition factors of $L(G)\downarrow_H$ have trivials but have twice as many $5$s and as many $10$s, then up to projective simple modules we have one of
\[ 10^9,5^2,1,\quad 10^8,5^4,1,\quad 10^{10},5^4,1^2.\]
Removing copies of $10,5^2,1$ (i.e., $P(5)$) yields $10^8$, $10^7,5^2$ and $10^8$, and so the first and last of these cannot occur without fixing a line, but the second can. However, the set of composition factors on $L(G)$ is $35^2,14_1,14_2^2,10^8,5^4,1$, (where $14_1$ is a factor of $5^{\otimes 2}$). This restricts to $\Alt(6)$ with composition factors $10^{11},9^6,8_1^3,8_2^3,5_1^2,5_2^5,1$, which happen to be those which were proved not to exist in Lemma \ref{lem:noalt6withspecificfactors}. Thus this embedding does not exist.
\end{pf}

\begin{prop} Suppose that $p=5$.
\begin{enumerate}
\item If $G=F_4$ and $H\cong \Alt(7)$ then $H$ fixes a line on $\Vmin$, and hence is not maximal.
\item If $G=E_6$ and $H\cong \Alt(7)$ then either $H$ fixes a line on $L(G)$, or is contained inside an algebraic $A_2$ subgroup with action of $\Vmin$ given by (up to duality) $8/19$, with two classes of $H$ being swapped by the graph automorphism. In particular, $H$ is not maximal.
\item If $G=E_6$ and $H\cong 3\cdot \Alt(7)$ with $Z(H)=Z(G)$, then $H$ fixes a line on $L(G)$, and hence is not maximal.
\item If $G=E_7$ and $H\cong \Alt(7)$, then $H$ fixes a line on either $\Vmin$ or $L(G)$, and hence is not maximal, or $H$ acts on $\Vmin$ as  \[(10\oplus 10^*)^{\oplus 2}\oplus 8^{\oplus 2}\qquad\text{or}\qquad 10\oplus 10^*\oplus P(8)^{\oplus 3}\oplus 8.\]
\item If $G=E_7$ and $H\cong 2\cdot \Alt(7)$ with $Z(H)=Z(G)$, then $H$ acts on $\Vmin$ as
\[20\oplus 4/14\oplus 14/4^*.\]
\item If $G=E_8$ and $H\cong \Alt(7)$, then either $H$ fixes a line on $L(G)$, hence is not maximal, or $H$ acts on $L(G)$ as
\[ 35^{\oplus 4}\oplus 15^{\oplus 4}\oplus 10\oplus 10^*\oplus 8/6\oplus 6/8.\]
(This is again the reduction modulo $p$ of the case for $p=0$.)

\end{enumerate}
\end{prop}
\begin{pf} Using the Brauer tree in Section \ref{sec:modulesalt}, the projective modules are
\[ 1/13/1,\quad 13/1,8/13,\quad 8/6,13/8,\quad 6/8/6.\]
For $G=F_4$ we consult the table in \cite{litterickmemoir} and see that $H$ has two trivial and three $8$-dimensional factors, so acts semisimply on $\Vmin$ and in particular fixes a line on $\Vmin$.

For $G=E_6$, there are three possible sets of composition factors on $\vmin$ according to \cite{litterickmemoir}. Two of these have trivial composition factors on $L(G)$, and no $13$s in either case, and hence these fix lines. The final case has composition factors $6$, $8$ and $13$ on $V_{\min}$ and $35^2,8$ on $L(G)$. This means that $u=(1,2,3,4,5)$ acts with factors $5^{15},3$ on $L(G)$ (since the $35$s are projective), and hence comes from class $A_4+A_1$, which acts with blocks $5^5,2$ on $V_{\min}$. This means that $H$ acts indecomposably on $V_{\min}$, hence as $8/6,13$ or its dual.

The restriction of $H$ to $L=\Alt(6)$ has structure $5_1\oplus 5_2\oplus 8/1,8$, which clearly fixes a line, and so lies in the $D_5$-parabolic (as the other line stabilizer is $F_4$ which has a trivial summand on $V_{\min}$). The $D_5$-parabolic acts with structure $10/16/1$, and obviously the embedding of $L$ into this acts as $8/8$ on the spin module and $5_1\oplus 5_2$ on the natural. In particular, this proves that the image of $L$ is uniquely determined up to conjugacy in the $D_5$ Levi.

To determine the number of classes in the $D_5$ parabolic we need to consider $1$-cohomology on the spin module, which is the action of the Levi on the unipotent radical of the parabolic. The $1$-cohomology of $8/8$ is $1$-dimensional, and so there are $q$ classes of $\Alt(6)$s inside the $D_5$ parabolic, where $G=E_6(q)$. We claim that all $q-1$ of those outside the $D_5$ Levi subgroup are permuted by the $T_1$ factor of the $D_5T_1$, at least in the adjoint form of $E_6$: to see this, any element of $T_1$ centralizes the maximal subgroup $D_5T_1$ of the $D_5$-parabolic, and so if it centralizes anything else in the parabolic, for example, another class of $\Alt(6)$, then it centralizes the whole parabolic so lies in $Z(G)$. Thus there is a unique class in the adjoint form of $E_6$ and either one or three in the simply connected form, depending on whether $|Z(G)|$ has order $1$ or $3$.

We can proceed in two different ways here, but both are similar. The first is to note that $L$ is contained in one of the two algebraic $A_2$s, say $X$, which act on $V_{\min}$ as $8/19$ (up to duality induced by the graph automorphism). Since the $19$-space is uniquely determined in the action of both $L$ and $H$, we must have that $H$ is contained in $\langle H,X\rangle$, which is just $X$ as $X$ is known to be maximal. (Notice that this inclusion was known in \cite{aschbacherE6Vun} but uniqueness of $H$ was not.)

Alternatively, if $\bar H$ is any other $\Alt(7)$ subgroup containing $L$, and with the same module structure on $V_{\min}$ as $H$, then both $H$ and $\bar H$ must fix the same $19$-space (with quotient $8$), since $L$ and $H$ both fix a unique $19$-space. This proves that either $H$ is uniquely determined up to conjugacy or that $H$ is not maximal. If $H$ is unique up to conjugacy, however, then it is contained in the algebraic $A_2$ above, and hence is not maximal either.

\medskip

If $3\cdot \Alt(7)$ embeds in $G$ with $Z(H)=Z(G)$, then there are three sets of composition factors for $L(G)$, two of which have at least three $1$s and at most one $13$, so fix a line. The remaining possibility has composition factors $21,6$ on $\Vmin$ and $35^2,8$ on $L(G)$. The $35$ is projective, so splits off, and thus $L(G)\downarrow_H$ is semisimple, with the action of $u=(1,2,3,4,5)$ on $L(G)$ being $5^{15},3$, so class $A_4+A_1$, which acts on $\Vmin$ as $5^5,2$. Thus $\Vmin\downarrow_H$ cannot be semisimple, but $\Ext^1(21,6)=0$ from the Brauer tree in Section \ref{sec:modulesalt}. Thus $H$ cannot embed with these factors, as needed.

\medskip

For $G=E_7$, the only set of composition factors embedding $\Alt(7)$ into $G$ with at least twice as many $13$s as $1$s for both $V_{\min}$ and $L(G)$ is
\[ (10,10^*)^2,8^2,\qquad (10,10^*),13^3,8^7,6^3,\]
and modules can be constructed so that $u=(1,2,3,4,5)$ acts on $V_{\min}$ as $5^{10},3^2$ and on $L(G)$ as $5^{26},3$. The first is semisimple, and the second is $10\oplus 10^*\oplus P(8)^{\oplus 3}\oplus 8$. Note that this exists as a copy of $H$ acting irreducibly on the natural module for the $A_7$ maximal-rank subgroup. There is also a possible set of composition factors for $2\cdot\Alt(7)$ embedding in $E_7$ with centres coinciding, that yields the same action of $u$ on the two modules: it acts as
\[ 20\oplus 4/14\oplus 14/4^*,\qquad 35^{\oplus 3}\oplus 10\oplus 10^*\oplus 8\]
on $V_{\min}$ and $L(G)$ respectively.

For $E_8$ and $p=5$ we again consult the tables of possible sets of composition factors from \cite{litterickmemoir}, and once we apply the statement that there needs to be twice as many $13$s as trivials, leaves us with five options, four with no trivial factors at all:
\[ [35^4,15^2,10,10^*],13^2,8^4,\qquad [35^4,15^3,10,10^*],13,8^3,6,\qquad [35^4,15^4,10,10^*],8^2,6^2,\]
\[[35,(10,10^*)^4],13,8^{12},6^4,\qquad 13^7,8^{14},6^7,1^3.\]
(We place brackets around the composition factors in blocks of defect zero, which simply become summands and $u=(1,2,3,4,5)$ acts projectively.)

In the first case we have at least $5^{46}$, and so there are few classes to which $u$ can belong, from Table \ref{tab:unipe8p5}. With just $13$ and $8$, the only indecomposable modules we can produce are $13$, $8$, and $13/8$ and $8/13$ on which $u$ acts as projective plus a block of size either $3$ or $1$. Thus $u$ acts only with blocks $5$, $3$ and $1$, so must come from class $A_4+A_2$ acting as $5^{46},3^5,1^3$; however, this cannot work either, so $H$ cannot embed with these factors.

In the second case, $u$ acts with at least $47$ blocks of size $5$, so must come from $A_4+A_2$ and act as $5^{48},4^2$. This means we can have no $13$ as a summand of $L(G)\downarrow_H$ (as it would contribute a $3$ to the action of $u$), thus it lies inside a self-dual summand, which is only $P(8)$, thus we get $P(8)\oplus 8$, another contradiction. Thus $H$ cannot embed with these factors either.

In the fourth case, the indecomposable modules containing an $8$ but no $1$s are (up to duality)
\[ 8,\quad 8/13,\quad 8/6,\quad 8/13,6,\quad P(8),\quad P(6).\]
With the exception of $8$, each of these has at least as many factors that are not of dimension $8$ as those that are, and since there are five factors that are not of dimension $8$ in the principal block for $L(G)\downarrow_H$, this means we have seven summands of $L(G)\downarrow_H$ of dimension $8$, thus $3^7$ contributing to the action of $u$ on $L(G)$. In addition, we have at least $5^{41}$ from the projective blocks in $L(G)\downarrow_{\langle u\rangle}$, and there is no such class in Table \ref{tab:unipe8p5}. Thus $H$ cannot embed with these factors.

In the fifth case, in order not to fix a line, $L(G)\downarrow_H$ must have three copies of $P(13)$, meaning that we have at least $5^{41}$, and leaves factors $13,8^{11},6^7$ to understand. In addition, applying the argument before yields $8^{\oplus 3}$ as a summand of $L(G)\downarrow_H$, so we have blocks in $u$ of $5^{41},3^3$; looking at Table \ref{tab:unipe8p5}, this means that $u$ lies in one of two classes, $A_4+2A_1$ acting as $5^{45},3^4,2^4,1^3$, or $A_4+A_2$ acting as $5^{46},3^5,1^3$. The four indecomposable modules contributing a $1$ to the action of $u$ are $1$, $6$, $8/13$ and $13/8$, and since we have three blocks of size $1$ (so we need a self-dual summand), this means we have a $6$ as a summand in $L(G)\downarrow_H$, thus another $8$ as a summand by the above argument. For the other $1^2$, we have either $6^{\oplus 2}$, which means we have another $8^{\oplus 2}$ and at least six blocks of size $3$ in $u$, not allowed, or we have $8/13\oplus 13/8$, not allowed because we only have a single $13$. This contradiction means that $H$ must fix a line if it embeds with these factors.

We are left with the third possibility, that $H$ embeds with factors projective plus $8^2,6^2$. Here there is a single possibility for the action of $u$, since we have at least forty-eight blocks in the action of $u$, namely $u$ acts as $5^{48},4^2$, which means that $H$ must act as
\[ 35^{\oplus 4}\oplus 15^{\oplus 4}\oplus 10\oplus 10^*\oplus 8/6\oplus 6/8.\]
This is the only module that satisfies the traces of elements and the action of the unipotent class, and so must be the reduction modulo $5$ of the embedding into $D_8$ given over $\C$.
\end{pf}

\begin{prop} Suppose that $p=3$.
\begin{enumerate}
\item If $G=F_4$ and $H\cong \Alt(7)$ then $H$ does not embed into $G$.
\item If $G=E_6$ and $H\cong \Alt(7)$ then $H$ fixes a line on $L(G)$.
\item If $G=E_7$ and $H\cong \Alt(7)$, or $H\cong 2\cdot \Alt(7)$ with $Z(H)=Z(G)$, then $H$ fixes a line on either $\Vmin$ or $L(G)$.
\item If $G=E_8$ and $H\cong \Alt(7)$ then either $H$ fixes a line on $L(G)$, and hence is not maximal, or acts on $L(G)$ with composition factors 
\[15^4,13^6,(10,10^*)^5,1^{10}.\]
\end{enumerate}
\end{prop}
\begin{pf} From Proposition \ref{prop:ext1simples}, the modules $10$ and $10^*$ each have $1$-dimensional $1$-cohomology, and $13$ has $2$-dimensional $1$-cohomology, with all other cohomology groups being zero. The projective covers $P(10)$ and $P(13)$ have structures
\[ 10/1/10^*,13/1/10\qquad \text{and}\qquad 13/1,1/13,10,10^*/1,1/13\]
respectively, and so if $M$ is a module whose socle consists solely of $10$s, $10^*$s and $13$s, with no trivial quotients, then it either has at most three socle layers or it contains a projective.

If $H$ embeds in $E_6$ then it does so with factors $15,13,(10,10^*)^2,6,1^3$ (remember that $\dim(L(G))=77$, not $78$, for $p=3$), and since $L(G)$ is self dual, if $13$ lies in the socle (or top) then it is a summand, and so does not need to be considered. The $6$ and $15$ lie in a non-principal block, and so form a summand of the form $6\oplus 15$ and will be ignored. We therefore may assume that the socle of $L(G)\downarrow_H$ consists of $10$s and $10^*$s, with possibly $13$ as a summand.

As $L(G)$ is self dual, if $P(10)$ is a summand of $L(G)\downarrow_H$ then so is $P(10^*)$, but this is not possible, as it uses too many $13$s. This also means that the $13$ is definitely a summand. However, we now do not have enough $10$-dimensional modules to cover three trivials, and so must fix a line on $L(G)$.

\medskip

When $G=E_7$ we have two possible sets of composition factors for embeddings of $H\cong\Alt(7)$, one with two trivial factors on $V_{\min}$ and no other modules from the principal block, hence centralizes a $2$-space, and one with nine trivial factors on $L(G)$ and a single $13$, with no $10$s, so centralizes at least a $7$-space. If $2\cdot\Alt(7)$ embeds in $E_7$ with centres coinciding, then there are twenty-two trivial composition factors, seven $13$s and $10,10^*$, so our module has negative pressure and fixes a line by Proposition \ref{prop:lowpressuremodules}.

\medskip

In $E_8$ there are four possible sets of composition factors for $L(G)\downarrow_H$, but only one case has negative pressure, hence fixes a line by Proposition \ref{prop:lowpressuremodules}. We examine the other three now. If the composition factors are $15^6,13^9,6^5,1^{11}$ then $H$ fixes a line, because without any $10$s we can only form $13/1,1/13$, and so need as many $13$s as $1$s not to fix a line.

With factors $15^2,13^{11},(10,10^*)^2,6^3,1^{17}$, we cannot have $P(10)$ as a factor because, together with its dual, we would need three $10$s, too many. Each $P(13)$ reduces the number of $13$s by three and $1$s by four, and we can have at most two of them, beyond which we can only have $13/1,1/13$, and this will not use up enough trivials. Hence we may assume that $L(G)\downarrow_H$, cut by the principal block, has three socle layers. Removing any simple summands, we have at most five $13$ and two $10$s in the socle, which support at most twelve trivials above, so we must fix a line again.

If the factors are $15^4,13^6,(10,10^*)^5,1^{10}$, however, then we cannot yet prove that $H$ fixes a line on $L(G)$, although it would be surprising if an embedding that has no fixed points does exist.
\end{pf}

Although in the last case we cannot prove that $H$ fixes a line on $L(G)$ or doesn't exist, we can prove a lot about $L(G)\downarrow_H$. Firstly, $x=(1,2,3)$ lies in class $2A_2$, and $y=(1,2,3)(4,5,6)$ lies in class $2A_2+2A_1$. Furthermore, $P(13)$ is not a summand of $L(G)\downarrow_H$, and the subgroups isomorphic to $\Alt(5)$ and $\PSL_2(7)$ both centralize a $2$-space on $L(G)$, with the $\Alt(5)$s having two trivial summands. Finally, $H$ fixes a line on $L(G)$ if and only if the subgroup $\Alt(6)$ does.

\begin{prop} Suppose that $p=2$.
\begin{enumerate}
\item If $G=F_4$ and $H\cong \Alt(7)$ then either $H$ or the image of $H$ under the graph automorphism fixes a line on $\Vmin$, and hence is not maximal.
\item If $G=E_6$ and $H\cong \Alt(7)$ then $H$ fixes a line on both $\Vmin$ and $L(G)$, and is hence not maximal.
\item If $G=E_7,E_8$ and $H\cong \Alt(7)$, or $G=E_6$ and $H\cong 3\cdot\Alt(7)$ with $Z(H)=Z(G)$, then $H$ fixes a line on $L(G)$, and hence is not maximal.
\end{enumerate}
\end{prop}
\begin{pf} From the table in Section \ref{sec:modulesalt}, there are three modules in the principal block: the projective covers have structure as follows:
\[ 1/14,20/1,1/20,14/1,\qquad 14/1/20/1,14/14,\qquad 20/1/14/1/20.\]

For $F_4$, $L(G)$ always has four trivial composition factors and at most one non-trivial factor in the principal block, so certainly fixes lines. Therefore if $\sigma$ denotes the graph automorphism then either $H$ or its image $H^\sigma$ under the graph automorphism fixes lines on $V_{\min}$. (This is because for $p=2$, $L(G)$ has composition factors $V_{\min}$ and $V_{\min}^\sigma$.)

For $E_6$, all embeddings of $H$ from \cite{litterickmemoir} have non-positive pressure on $V_{\min}$, hence fix lines or hyperplanes. For $3\cdot \Alt(7)$, the two possibilities for composition factors on $L(G)$ both have pressure $-2$, so again fix lines.

\medskip

For $G=E_7$ there are more possibilities, but in each case there are significantly more trivial factors than non-trivial factors in the principal block (even taking into account the fact that $L(G)$ has a trivial submodule in characteristic $2$), hence $H$ fixes lines on $L(G)$ again.

\medskip

For $G=E_8$, we need to be a little more precise: from the structure of the projectives, it is clear that we need at least three $14$s or $20$s for every two $1$s. This reduces the possibilities down to three, namely
\[ 20^6,14^8,1^8,4,4^*,\qquad 20^5,14^9,1^8,6,4,4^*,\qquad 20^4,14^{10},1^8,6^2,4,4^*.\]
The element $v=(1,2,3,4)(5,6)$ acts projectively on both $20$ and $4$, and acts with blocks $4,2$ and $4^3,2$ on $6$ and $14$ respectively. Thus, if the embeddings were semisimple, the action of $v$ on these three modules has blocks
\[ 4^{56},2^8,1^8,\quad 4^{55},2^{10},1^8,\quad 4^{54},2^{12},1^8.\]

Write $L(G)\downarrow_H=A\oplus B$, where $A$ is a sum of projectives and the $\{4,4^*,6\}$-radical, and $B$ is a module in the principal block with no projective summands. In particular, $B$ has three socle layers, and examining the structure of $P(14)$, we see that $B$ has at least as many $20$s as trivials. Since the same is true for $P(20)$, we need a copy of $P(14)$ in $L(G)\downarrow_H$ for every trivial composition factor above the number of $20$s, i.e., two, three and four respectively. Since there are three $14$s in $P(14)$, this means the third case is not possible, and the second case cannot work either as all $14$s are used up in the $P(14)^{\oplus 3}$, leaving none to cover the remaining two trivials. In the first case $P(14)^{\oplus 2}$ uses up $20^2,14^6,1^4$ and leaves $20^4,14^2,1^4,4,4^*$; this gives sixty blocks to the action of $v$ already, and so we cannot have any more $P(14)$s or $P(20)$s, but leads to a contradiction as, with three socle layers, $B$ must have at least twice as many non-trivial composition factors as trivial ones. Hence any embedding of $H$ into $E_8$ fixes a line.
\end{pf}

\section{$\Alt(8)$}

The cases to consider are firstly $p>7$, and then $p=7,5,3,2$ in that order. We continue our assumption that $E_6$ and $E_7$ are used to denote the simply connected forms, so that $Z(E_6)$ and $Z(E_7)$ have orders $3$ and $2$ respectively.

\begin{prop} Suppose that $p\neq 2,3,5,7$.
\begin{enumerate}
\item If $G=F_4,E_6$ and $H\cong \Alt(8)$ then $H$ does not embed into $G$.
\item If $G=E_7,E_8$ and $H\cong \Alt(8)$ then $H$ fixes a line on $L(G)$ and hence is not maximal.
\item There is no embedding of $H\cong 2\cdot \Alt(8)$ into $G=E_7$ with $Z(H)=Z(G)$.
\end{enumerate}
\end{prop}
\begin{pf} This follows immediately from the tables in \cite{litterickmemoir} together with Lemma \ref{lem:smallspacestabilizers}.
\end{pf}

\begin{prop} Suppose that $p=7$.
\begin{enumerate}
\item If $G=F_4,E_6$ and $H\cong \Alt(8)$ then $H$ does not embed into $G$.
\item If $G=E_7,E_8$ and $H\cong \Alt(8)$ then $H$ fixes a line on $L(G)$ and hence is not maximal.
\item There is no embedding of $H\cong 2\cdot \Alt(8)$ into $G=E_7$ with $Z(H)=Z(G)$,
\end{enumerate}
\end{prop}
\begin{pf} As before, using the Brauer tree in Section \ref{sec:modulesalt} and Lemma \ref{lem:covertrivialBrauer}, the only projective containing the trivial module other than $P(1)$ is $P(19)$, which has the form $19/1,45/19$. We then consult the tables in \cite{litterickmemoir}: for $E_7$ there is a unique embedding of $H$, which has two trivial factors and no $45$ on $L(G)$, so fixes lines, and for $E_8$ the three possible embeddings of $H$ each have at least three trivial factors and at most one $45$, so all fix lines.
\end{pf}

\begin{prop} Suppose that $p=5$.
\begin{enumerate}
\item If $G=F_4,E_6$ and $H\cong \Alt(8)$ then $H$ does not embed into $G$.
\item If $G=E_7,E_8$ and $H\cong \Alt(8)$ then $H$ fixes a line on $L(G)$ and hence is not maximal.
\item There is no embedding of $H\cong 2\cdot \Alt(8)$ into $G=E_7$ with $Z(H)=Z(G)$,
\end{enumerate}
\end{prop}
\begin{pf} Using the Brauer trees from Section \ref{sec:modulesalt} the projective indecomposable modules are
\[ 1/13/1,\quad 13/1,43/13,\quad 43/13,21_1/43,\quad 21_1/43/21_1,\quad 7/21_2/7,\quad 21_2/7,21_2/21_2.\]
Thus a module needs at least twice as many $13$s as $1$s in order not to fix a line. For $E_7$ there are three possible sets of composition factors, two of which have one trivial and at most one $13$, so fix lines on $L(G)$. The remaining set of composition factors acts with factors $21_2^2,7^2$ on $V_{\min}$ and acts as $20\oplus 43\oplus 70$ on $L(G)$. Clearly $u=(1,2,3,4,5)$ acts as $5^{26},3$ on $L(G)$, and so must come from class $A_4+A_2$, so act as $5^{10},3^2$ on $V_{\min}$. This can be realized by $21\oplus P(7)$.

However, this restricts to $\Alt(7)$ as 
\[ 13\oplus 8\oplus 1/13/1\oplus 6/8/6,\]
and so lies in the stabilizer of a unique line, either an $E_6$-parabolic or a $B_5$ type subgroup, the latter acting with composition factors $1,1,11,11,32$, clearly not possible. (To see this, the derived subgroups of the stabilizers of the other orbits from \cite[Lemma 4.3]{liebecksaxl1987} all lie inside an $E_6$-parabolic, hence stabilize more than one line.) Therefore this $\Alt(7)$ lies inside an $E_6$ parabolic, with composition factors $13,8,6$ on the minimal module $W$ for $E_6$. The $E_6$ parabolic acts as $1/W/W^*/1$, and so there must be a non-split extension between the $6$ and $8$ in $W$, but the $13$ must split off, so that up to duality we have that $\Alt(7)$ acts on $W$ as $8/6\oplus 13$. However, $u$ acts on this module as $5^4,4,3$, not listed in \cite[Table 5]{lawther1995}. Thus this embedding of $\Alt(8)$ into $E_7$ does not exist, and $H$ fixes a line, as needed.

For $E_8$ every possible set of composition factors from \cite{litterickmemoir} has at least three trivials and at most one $13$, so these all fix lines.
\end{pf}

\begin{prop} Suppose that $p=3$.
\begin{enumerate}
\item If $G=F_4,E_6$ and $H\cong \Alt(8)$ then $H$ does not embed into $G$.
\item If $G=E_7$ and $H\cong \Alt(8)$ then either $L(G)\downarrow_H$ has a trivial summand or $(1,2,3)$ lies in a generic unipotent class for $\Vmin$, and so there exists a positive-dimensional subgroup stabilizing the same subspaces of $\Vmin$ as $H$. In either case, $H$ is not maximal.
\item There is no embedding of $H\cong 2\cdot \Alt(8)$ into $G=E_7$ with $Z(H)=Z(G)$,
\item If $G=E_8$ and $H\cong \Alt(8)$ then $H$ fixes a line on $L(G)$.
\end{enumerate}
\end{prop}
\begin{pf} For $E_7$ there is a unique set of composition factors, both on $V_{\min}$ and on $L(G)$. On $V_{\min}$ we have that $H$ acts with factors $21^2,7^2$. Since the $21$ lies in a separate block to the $7$, and the $7$ has no self extensions, $V_{\min}$ must have $7^{\oplus 2}$ as a summand. The other summand is either $21^{\oplus 2}$ or $21/21$. The element $x=(1,2,3)$ acts on $7$ with blocks $3,1^4$ and on $21$ as $3^5,1^6$, so on $21/21$ as $3^{10},2^6$. This means that $x$ acts on $\Vmin$ as either $3^{12},1^{20}$ -- so lies in class $A_2$, which is generic, and so we are done by Lemma \ref{lem:unipotentreduction} -- or acts like $3^{12},2^6,1^8$ -- so lies in class $A_2+A_1$, and acts on $L(G)$ as $3^{37},2^8,1^6$.

The composition factors of $H$ on $L(G)$ are also uniquely determined as $35^2,21,13,7^4,1$. The $1$ can be covered only by $35/1,7/35$, on which both classes of unipotent element act as $3^{26}$. Even if the rest of the module is semisimple, the $7^4$ contributes four blocks of size $3$, the $21$ another five, and the $13$ gives three, taking the total to $38$, too many. Thus $35/1,7/35$ cannot be a subquotient of $L(G)\downarrow_H$, and so $L(G)\downarrow_H$ has a trivial summand and hence $H$ fixes a line of $L(G)$.

For $E_8$, there are two possible sets of composition factors for $L(G)\downarrow_H$, each with four trivial composition factors. As $H^1(H,M)=0$ unless $M$ is either $13$ or $35$, in which cases it is $1$-dimensional, the fact that there are exactly three such composition factors in $L(G)\downarrow_H$ implies that it has negative pressure and so $H$ fixes a line on $L(G)$ by Proposition \ref{prop:lowpressuremodules}, as needed.
\end{pf}

\begin{prop} Suppose that $p=2$.
\begin{enumerate}
\item If $G=F_4$ and $H\cong \Alt(8)$ then either $H$ or its image under the graph automorphism stabilizes a line on $\Vmin$, and hence is not maximal.
\item If $G=E_6$ and $H\cong \Alt(8)$ then $H$ has a trivial summand on $\Vmin$ and lies inside a conjugate of $F_4$ or the $D_5$ Levi subgroup. In particular, $H$ is not maximal in $G$.
\item If $G=E_7$ and $H\cong \Alt(8)$ then $H$ fixes a line on either $\Vmin$ or $L(G)$, and hence is not maximal.
\item If $G=E_8$ and $H\cong \Alt(8)$ then $H$ fixes a line on $L(G)$, and hence is not maximal.
\end{enumerate}
\end{prop}
\begin{pf}
We will use Proposition \ref{prop:alt8compfactors} firstly, together with the dimensions of $1$-cohomology from Table \ref{tab:exts}. If $G=F_4$ then we have two possible sets of composition factors, swapped by the graph automorphism, and the first one has pressure $-3$, so fixes a line on $\Vmin$.

\medskip

For $G=E_6$, there are again two sets of composition factors for $H$: the first has pressure $-4$ on $\Vmin$ and $-3$ on $L(G)$, so fixes a line on both. Moreover, since in $\Vmin\downarrow_H$ only one composition factor has $1$-cohomology, and there are five trivials, this means there are at least three trivial summands in $\Vmin\downarrow_H$.

If the composition factors of $\Vmin\downarrow_H$ are $14,6^2,1$, then we can check by computer that there are no modules of the form $6/1/6$ or $6/1/14$, so either $H$ fixes a line or hyperplane on $\Vmin$ or $\Vmin\downarrow_H$ is indecomposable, with socle $6$ up to duality. Since $\Ext^1(6,6)=0$, we get that $\Vmin\downarrow_H$ has shape $6/1,14/6$. There is a unique such module inside $P(6)$, but this has a trivial quotient, so $H$ always fixes a line or hyperplane on $\Vmin$. (Notice that the dimensions are incompatible with coming from the $D_5$ parabolic, as this acts with factors $16,10,1$, so in fact $H\leq F_4$ and has a trivial summand.)

\medskip

For $G=E_7$, Proposition \ref{prop:alt8compfactors} states that there are three possible sets of composition factors, with the first and third stabilizing lines on $L(G)$ (remember to remove one trivial since $\dim(L(G))=132$ for $p=2$) as they have pressure $-11$ and $-6$ respectively. The second case, of $14^2,6^4,1^4$, has pressure $2$, so assuming that $H$ does not fix a line on $\Vmin$, there are at most two composition factors in the socle by Proposition \ref{prop:lowpressuremodules}. The $\{1,6,14\}$-radicals of $P(6)$ and $P(14)$ have two and one trivial composition factor respectively, so the socle must therefore be $6{\oplus 2}$: however, we cannot have a submodule $14/6,6$ of $\Vmin\downarrow_H$, since this has pressure $3$. However, the $\{1,6,14\}$-radical of $P(6)$ is
\[ 6/14/1,6/1,14/6,\]
and so to support the trivial in the third socle layer we clearly need the $14$ in the second socle layer, so in fact $(14/6)^{\oplus 2}$ is a submodule of $\Vmin\downarrow_H$, a contradiction since it has pressure $4$.

\medskip

The remainder of this proof concerns $E_8$, and is very delicate and long. We will proceed in stages.

\medskip

\noindent \textbf{Step 1: Identifying the two difficult cases} Firstly, using Proposition \ref{prop:alt8compfactors}, we have exactly four possibilities for the set of composition factors on $L(G)$, namely
\[ 14,6^{10},(4,4^*)^{16},1^{46},\qquad 14^8,6^{17},(4,4^*)^2,1^{18},\quad 20^4,14^4,6^8,(4,4^*)^7,1^8,\qquad 20^4,14^{10},6^2,4,4^*,1^8.\]
Since the modules $6,14,20,20^*$ each have $1$-dimensional $1$-cohomology by Proposition \ref{prop:ext1simples}, the first case fails by Proposition \ref{prop:lowpressuremodules}.

In the second possibility, there are no $20$s at all, and so Lemma \ref{lem:alt8no20s} applies to $L(G)\downarrow_H$, and we get that $H$ fixes a line on $L(G)$ because it has fewer $6$s than trivials.

We now apply Lemma \ref{lem:frobrep} to $L(G)\downarrow_L$, where $L=\Alt(7)$. Since we proved in the previous section that $L$ fixes a line on $L(G)$, this yields a map from the permutation module of $H$ on the cosets of $L$ to $L(G)\downarrow_H$. Since this module is $1/6/1$, and we assume that $H$ does not fix a line on $L(G)$, we get that $1/6$ is a submodule of $L(G)\downarrow_H$.

\medskip

We examine the third set of composition factors of $H$ on $L(G)$ now, namely $20^4,14^4,6^8,(4,4^*)^7,1^8$. Let $v=(1,2,3,4)(5,6)$, an element of order $4$ that acts projectively on both $20$ and $4$, and acts with a single block of size $2$ and otherwise projectively on $6$ and $14$. We write $W$ for the subquotient of $L(G)\downarrow_H$ obtained by quotienting out by all submodules and taking the kernel of any quotient by factors of dimension $4$, i.e., take the $\{4,4^*\}$-residual of $L(G)$, and quotient out by the $\{4,4^*\}$-radical of that to get $W$. This does not alter the action of $v$, since it acts projectively on these factors. Thus $W$ is a submodule of copies of $P(6)$, $P(14)$ and $P(20)$.

\medskip

\noindent \textbf{Step 2: Finding $1/14/1$ subquotients} We show that there are four disjoint subquotients of $W$ of the form $1/14/1$. To see this, we first note that there does not exist any module $M$ with no $14$s as composition factors, with at least one trivial composition factor, and with no trivial submodule or quotient. This is proved easily: the largest submodules of $P(6)$ and $P(20)$ with no trivial quotients and no $14$s are
\[ 6/4,4^*/6/4,4^*/6\qquad \text{and}\qquad 4/20.\]
Hence we need a $14$ above or below every trivial in $W$. However, there are four $14$s and eight $1$s in $W$, so we must have $1/14/1$ four times.

\medskip

\noindent \textbf{Step 3: Connecting the $1/14/1$s to $6$s} The socle of $W$ has at most four $6$s and two modules of dimension $20$. We will prove that there is at most one $1/14/1$ lying above any $6$ in the socle.

Since the dimension of $P(6)$ is $320$, we cannot have $P(6)$ in $L(G)\downarrow_H$, so we may remove the top of $P(6)$, to leave its Jacobson radical. Since the top of $W$ consists of $6$s, $20$s and $20^*$s, we may take the $\{1,4,4^*,14\}$-residual of this module, a module of dimension $290$. This will contain any submodule of $L(G)\downarrow_H$ lying above a $6$ in the socle. We now take the $\{1\}'$-residual and then quotient out by the $\{1\}'$-radical to find the smallest subquotient containing all copies of $1/14/1$, and this leaves us with
\[ 1\oplus 1,1/14,20,20^*/1,1,1.\]
Thus we can have at most one $1/14/1$ lying above each $6$ in the socle of $W$.

Performing the same calculation with $P(20)$ yields
\[ 1/14/1/20,20^*/1/14/1,\]
which of course contains two $1/14/1$s, but needs $20\oplus 20^*$ inside the module to obtain them both.

\medskip

\noindent \textbf{Step 4: No $6^{\oplus 4}$ in the socle} Suppose that there are four $6$s in $\soc(W)$, so that there are no $6$s in $\rad(W)/\soc(W)$. We construct the largest submodule of $P(6)$ with exactly one $6$, yielding
\[ 20,20^*/1/1,4,4^*,14/6,\]
which has no $1/14/1$ inside it, so we cannot stack enough $1/14/1$s, and hence we cannot have $6^{\oplus 4}$ inside $\soc(W)$. In particular, this means that we have one of the following socles:
\[ 6^3,20,\qquad 6^2,20,\qquad 6^3,20^2,\qquad 6^2,20^2,\]
where here $20^2$ means two modules of dimension $20$, either $20$ or $20^*$.

\medskip

\noindent \textbf{Step 5: Each $1/14/1$s requires a $10$ or $20^*$ and a $6/4/6$ or $ 6^*/4/^*$ around it} We construct a submodule $M$ of $P(6)$ by the following process:
\begin{enumerate}
\item Start with the $\{4,4^*,6\}$-radical of $P(6)$ (this has dimension $40$);
\item Add to this all trivials possible (there is only one);
\item Add to this all $14$s possible (there are three);
\item Add to this all trivials possible (there are three);
\item Add to this all $4$s, $4^*$s and $6$s possible (this yields a module of dimension $112$);
\item Take the $\{6\}'$-residual of this module (this yields a module of dimension $110$).
\end{enumerate}
This process yields a module with five socle layers. If it were possible to construct a submodule of $P(6)$ with no $20$s or $20^*$ and with a $1/14/1$ subquotient but no trivial submodules or quotients, it would lie inside $M$. It has socle layers
\[ 6,6,6/4,4,4^*,4^*,14,14/1,6,6,6/1,4,4^*,14/6.\]
But this clearly does not have a $1/14/1$ inside it, and so we must need a $20$ or $20^*$ above any such submodule to prevent the $1/14/1$ floating to the top of $W$.

Since we need at least two $6$s in the socle of $W$, with $1/14/1$s above them, this means we need at least two $20$s above that, so in particular they cannot be used to string together the two $1/14/1$s in the subquotient of $P(20)$ above. We therefore see that we need at least four factors in $\soc(W)$, so either $6^2,20^2$ or $6^3,20$. However, with $6^3,20$, we would need three $20$s above the $6^3$, meaning the module cannot be self dual. To see this, placing a $20$ above a $1/14/1$ creates a module
\[ 20/1/14/1,\]
and since there are four $1/14/1$s and four $20$s or $20^*$s, for each $20$ placed above a $1/14/1$ in $W$ there must be a $20^*$ placed below a $1/14/1$, to maintain self duality of $W$. We therefore see that there are two factors of dimension $20$ in $\soc(W)$, and two $6$s.

We now construct the same module $M$ as in the process above, but instead of step (v), we place $20$ on top of it. (We could use $20^*$ as well, but this would yield the image under the outer automorphism of $H$.) Instead of (vi) we then take the $\{20\}'$-residual. This will contain any smallest submodule of $P(6)$ with no trivial quotient or submodule and a $1/14/1$ inside it.

This process produces the module $M$ with socle layers
\[ 20/1/4,14,20/1,6/1,4,4^*,14/6\]
and its dual $M^*$ has structure
\[ 6/4/1,6/14,14/1,1,4^*,4^*/20^*,20^*;\]
from this we easily see the submodule $6/4^*/6$ inside $M$ (there is a unique such uniserial module). Furthermore, since the only copy of $6$ not in the socle is at the top of this submodule, we can see what happens if we remove this submodule by requiring a unique $6$ in our module $M$: even before performing stage (vi) and removing quotients, the structure is
\[ 20/1/1,4,4^*,14/6,\]
and so does not contain a $1/14/1$. Hence for every $6$ supporting a $1/14/1$ above it, we need a $6/4/6$ or $6/4^*/6$ as a submodule.

\medskip

Bringing this to a conclusion, for each $1/14/1$ in our module, we need a $20$ or $20^*$ above it, and a $6/4/6$ or a $6/4^*/6$ below it, or vice versa. In particular, we will need four submodules or quotients of the form $6/4/6$ (up to duality) in $W$, in order for it not to have any trivial submodules or quotients.

\medskip

\noindent \textbf{Step 6: Contradiction for third set of composition factors} The element $v=(1,2,3,4)(5,6)$ acts projectively on this module, so since it already acts projectively on $1/14/1$, if there are $i$ disjoint subquotients of the form $6/4/6$ (or $6/4^*/6$) then we have that $v$ acts on $L(G)$ with blocks $4^{58+i},2^{8-2i}$. However, from Table \ref{tab:unipe8p4} we see that $i=2$. However, we have just shown that $i=4$, a contradiction. Thus $W$, and hence $L(G)\downarrow_H$, has a trivial submodule or quotient, as needed.

\medskip

\noindent \textbf{Step 7: Elimination of all possible socles for fourth set of composition factors} The last case to consider is $64^2,20,20^*,14^4,6^2,1^4$. The socle of $W$ is a submodule of $6\oplus 14^{\oplus 2}\oplus 20$, and we have already shown above that it contains $6$. If it is all of them then there must be a submodule
\[ 1,1,1,1/6,14,14,20,\]
and each of $6$, $14$ and $20^*$ must cover one of the trivials. However, there are two extensions of $6$ by this module, both of which lie above the $14$s and not the $1$s, so that this cannot be the socle.

On the other extreme, if the socle is just $6$ then we cannot work either: the largest submodule of $P(6)$ with a single $6$ has two trivial factors, so cannot work.

We do a similar thing if the socle has two or three factors: compute the largest submodule $M$ of the appropriate sum of projectives, with the proviso that the composition factors of $\top(M)/\soc(M)$ do not include $6$, and do not include $20$ if $20\in \Soc(W)$ and similarly if $14^{\oplus 2}\leq \soc(W)$, and remove any simple quotients from $M$ whose duals are not isomorphic to anything in the socle. 

If the socle has factors $14,6$ then we construct the module $M$, add all four $6$s onto $M$, then remove all quotients not of dimension $14$ or $6$ to leave a module of dimension $122$ with two trivial factors, so this doesn't work either.

We end with the socle being $6\oplus 14^{\oplus 2}$. Define $A_1$ to be the preimage in $P(14)$ of the $\{1,4,4^*\}$-radical of $P(14)/\soc(P(14))$. We can construct a single extension of this module by $20$, yielding a module $A_2$. In $P(14)\oplus P(14)$, adding a single $20$ to $A_1\oplus A_1$ always makes, by choosing diagonal submodules if necessary, a copy of $A_2\oplus A_1$, up to isomorphism. The preimage of the $\{1,4,4^*\}$-radical of $P(6)/\soc(P(6))$ is a module $A_3$, on which one cannot place a $20$ or $20^*$. For $P(14)$ we construct similar modules $A_1$ and $A_2$ as above, but this time by taking the $\{1,4,4^*\}$-radical. To the sum of $A_1$, $A_2$ and $A_3$, we then place as many $1$s, $4$s, $4^*$s and $20^*$s as we can, to make a module of dimension $108$ with six trivial factors. Placing on top of this as many $14$s as we can, then as many $6$s as we can, still yields three trivial quotients, so we can only cover three factors, and this is the final contradiction.

Hence there is no possibility for $\soc(W)$, so $H$ must always fix a line on $L(G)$, as needed.
\end{pf}

\section{$\Alt(9)$}

Characteristics other than $2$ and $3$ have already been solved by Litterick in \cite{litterickmemoir}, so we consider only $p=3$ and then $p=2$. We continue our assumption that $E_6$ and $E_7$ are used to denote the simply connected forms, so that $Z(E_6)$ and $Z(E_7)$ have orders $3$ and $2$ respectively.

\begin{prop} Suppose that $p=3$.
\begin{enumerate}
\item If $G=F_4,E_6$ and $H\cong \Alt(9)$, then $H$ does not embed into $G$.
\item If $G=E_7$ and $H\cong \Alt(9)$ then either $L(G)\downarrow_H$ has a trivial summand or $(1,2,3)$ lies in a generic unipotent class for $\Vmin$, and so there exists a positive-dimensional subgroup stabilizing the same subspaces of $\Vmin$ as $H$. In either case, $H$ is not maximal.
\item If $G=E_7$ and $H\cong 2\cdot \Alt(9)$ then there is no embedding of $H$ into $G$ with $Z(H)=Z(G)$.
\item If $G=E_8$ and $H\cong \Alt(9)$ then $H$ fixes a line on $L(G)$, so is not maximal.
\end{enumerate}
\end{prop}
\begin{pf} Recall the simple modules and dimensions of $\Ext^1$ from Section \ref{sec:modulesalt}. In characteristic $3$ we need to consider $G=E_7$ and $G=E_8$. In $E_7$, $H$ acts on $V_{\min}$ with factors $21^2,7^2$, and so the element $z=(1,2,3,4,5,6,7,8,9)$, which acts on $21$ as $9^2,3$ and on $7$ with a single block, must come from either class $A_6$, acting as $9^4,7^2,3^2$, or from class $E_6(a_1)$, acting as $9^6,1^2$.

In the previous section, we proved that $L=\Alt(8)$ either contained a unipotent element from the generic class $A_2$ -- and hence is contained in a positive-dimensional subgroup by Lemma \ref{lem:unipotentreduction} -- or has a trivial summand on $L(G)$ and $x=(1,2,3)$ comes from class $A_2+A_1$. In addition, $y=(1,2,3)(4,5,6)$ acts on $V_{\min}$ as $3^{18},1^2$, and comes from class $2A_2$ (as the trivial summand on $L(G)$ means it cannot come from class $2A_2+A_1$, which acts as $3^{42},2^2$).

If the socle of $V_{\min}\downarrow_H$ is $7\oplus 21$ then the module is one of $7/21\oplus 21/7$ or $7,21/7,21$, with both modules being uniquely determined. The element $y$ acts on the former module as $3^{12},1^{20}$, not right, so this is not the correct embedding. Thus we may assume that $H$ acts as either $7,21/7,21$ or as $7\oplus 7\oplus 21/21$. The element $z$ acts on the first module as $9^6,1^2$ and on the second as $9^4,7^2,3^2$, so acts on $L(G)$ as either $9^{11},7^4,3,1^3$ or $9^{14},7$.

The composition factors of $H$ on $L(G)$ are $35^2,27,21,7^2,1$, with the $27$ splitting off as a summand since it lies in a separate block. The action of $z$ on each composition factor is
\[ 9^3,7,1,\qquad 9^3,\qquad 9^2,3,\qquad 7,\qquad 1.\]

As $L$ fixes a line on $L(G)$, by Lemma \ref{lem:frobrep} there is a map from the permutation module of $H$ on $L$, which is $1/7/1$, to $L(G)\downarrow_H$. As we assume that $H$ fixes no line on $L(G)$, this means that $1/7$ is a submodule of $L(G)\downarrow_H$.

Since $z$ acts indecomposably on $1/7$, we cannot have that the factors of $z$ on $L(G)$ are $9^{11},7^4,3,1^3$, and so $z$ acts as $9^{14},7$ on $L(G)$. In particular, the only possible simple summands of $L(G)\downarrow_H$ are $27$ and $7$, and in particular $35$ cannot be. However, since $1/7$ is a submodule of $L(G)\downarrow_H$ (and $7/1$ is a quotient) we cannot have that $7$ is a summand. Hence the socle of $L(G)\downarrow_H$ (apart from the $27$) is either $7\oplus 35$ or $7$.

The largest submodule of $P(7)$ with a single $7$ and other factors $1$, $21$ and $35$ is $21,35/1,21,35/7$. Since this has only three socle layers, we cannot have that the $21$ has an extension with the $35$s (since as $L(G)$ is self dual we would need socle layers of $7$, $35$, $21$, $35$ and $7$). But then the two $35$s cannot lie above the $7$, a contradiction. Thus the socle of $L(G)\downarrow_H$ must contain $35$.

Removing the $35$s from the top and bottom of $L(G)\downarrow_H$, we are left with a self-dual module with factors $7^2,21,1$, and with $1/7$ as a submodule. Since there is no uniserial module $7/1/7$, this must be $M=7/1,21/7$.

We have that $\Ext^1(M,35)$ is $1$-dimensional, but when this module is constructed, the element $y$ has a block of size $2$ in its action. However, above we saw that $y$ acts as $3^{41},1^7$ on $L(G)$, and this block of size $2$ reduces the socle of $L(G)\downarrow_{\langle y\rangle}$ to at most $47$, not allowed. This final contradiction proves that $H$ fixes a line on $L(G)$, as needed for $G=E_7$.

\bigskip

We finally consider $G=E_8$, where the composition factors of $L(G)\downarrow_H$ are uniquely determined and are $35^2,27,21^5,7^6,1^4$. 

As $H^1(H,M)=0$ unless $\dim(M)=7,35,41$ (and of course we can ignore $41$), by taking the $\{1,21,27\}$-residual of $L(G)\downarrow_H$ and then quotienting out by the $\{1,21,27\}$-radical, we get a subquotient $W$ of $L(G)\downarrow_H$ whose top and socle consist of $7$s and $35$s. Suppose that $35$ lies in the socle: this means that $W$ is a submodule of $P(35)\oplus P(7)^{\oplus i}$ for some $i$, and removing the $35$ from the socle and top we have no more $35$s. The $\{1,7,21\}$-radical of $P(7)$ is $7/1,21/7$, and the preimage in $P(35)$ of the $\{1,7,21\}$-radical of $P(35)/\soc(P(35))$ has the form $1,7,21/35$. As we need to hide four trivials, we clearly need the socle to be $35\oplus 7^{\oplus 3}$, whence that $7$ lying above the $35$ cannot lie in $W$. However, $\Ext^1(35,1,21/35)$ has dimension $2$, but all of these extensions have a trivial quotient, whence we fix a line. Thus the socle of $W$ consists solely of $7$s.

The $\{1,7,21,35\}$-radical of $P(7)$ is
\[ 7/1,21,35/7,7,21,35/1,21,35/7,\]
and so since we need to conceal four trivial composition factors, the socle of $W$ needs at least two $7$s. However, in order to have both trivial composition factors, we need the whole of this module as a submodule of $W$, which contains three $35$s, a contradiction. Thus we need as many $7$s in the socle as $1$s above them, i.e., four, which is not possible. Thus $H$ always fixes a line on $L(G)$.
\end{pf}

For $p=2$, if $G\neq E_8$ then \cite[Theorem 1]{litterickmemoir} gives the result, but we add it for sake of completeness since the proofs are easy using Proposition \ref{prop:alt9compfactors}.

\begin{prop} Suppose that $p=2$.
\begin{enumerate}
\item If $G=F_4$ and $H\cong\Alt(9)$, then either $H$ or its image under the graph automorphism fixes a line on $\Vmin$, and so is not maximal.
\item If $G=E_6$ and $H\cong\Alt(9)$, then $H$ has a trivial summand on $\Vmin$, and hence is not maximal
\item If $G=E_7$ and $H\cong\Alt(9)$, then $H$ fixes a line on $\Vmin$, and hence is not maximal.
\item If $G=E_8$ and $H\cong\Alt(9)$, then $H$ fixes a line on $L(G)$, and hence is not maximal.
\end{enumerate}
\end{prop}
\begin{pf} The dimensions of $1$-cohomology for simple $H$-modules are given in Table \ref{tab:exts}, and the composition factors for $\Vmin\downarrow_H$ and $L(G)\downarrow_H$ are given in Proposition \ref{prop:alt9compfactors}. These together prove that for $G=F_4$ either $H$ or its image under the graph automorphism has pressure $-2$, and for $E_6$ either $\Vmin\downarrow_H$ has pressure $-3$ or has composition factors $26$ and $1$, hence fixes a line or hyperplane on $\Vmin$. However, if the factors are $8_1,8_2,8_3,1^3$ then the trivial split off and if the factors are $26,1$ then $H$ cannot lie in a $D_5$ parabolic, with factors $10,16,1$ on $\Vmin$, hence in $F_4$, which acts as $26\oplus 1$, proving the result. For $E_7$ we have factors $8_1^2,8_2^2,8_3^2,1^8$, hence having eight trivial summands, or $26^2,1^4$, of pressure $0$, hence $H$ fixes a line.

\medskip

We are left with the case of $G=E_8$. The proof of this is quite long, and we will break it into stages. As we saw in Section \ref{sec:modulesalt}, the simple modules here are $1$, $8_1$ and $8_2$, permuted by the outer automorphism, $8_3$, $20$, $20^*$, $26$, $48$ and $78$, plus modules that do not appear in our analysis. The modules $8_3$ and $48$ lie in a non-principal block of $H$.

\medskip

\noindent \textbf{Step 1: Eliminating all but one set of composition factors on $L(G)$.} By Proposition \ref{prop:alt9compfactors} we may assume that the composition factors of $H$ acting on $L(G)$ are one of
\[ 26,8_1^8,8_2^8,8_3^8,1^{30},\quad 26^8,8_1,8_2,8_3,1^{16},\quad 26^4,(20,20^*)^2,8_1^5,8_2^2,1^8\]
\[48^2,26^4,8_3^5,1^8,\quad  48,26^2,(20,20^*)^2,8_1^3,8_2^3,8_3^2,1^4.\]
As $H^1(H,M)$ is $2$-dimensional for $M=26$ and $1$-dimensional for $78$, $20$ and $20^*$, we have that these sets of composition factors yield modules of pressure $-28$, $0$, $4$, $0$ and $4$, respectively, so only the third and the fifth case might not fix a line on $L(G)$.

\medskip

Let $L\cong \Alt(8)\leq H$. Since $L$ fixes a line on $L(G)$, we may apply Lemma \ref{lem:frobrep} to find a map from the permutation module of $H$ on the cosets of $L$, which is $1\oplus 8_3$, to $L(G)\downarrow_H$. Thus if $8_3$ is not a submodule of $L(G)\downarrow_H$ then $H$ fixes a line on $L(G)$. This deals with the third possibility for the composition factors of $L(G)\downarrow_H$ above, so we may assume from now on that the composition factors of $L(G)\downarrow_H$ are
\[ 48,26^2,(20,20^*)^2,8_1^3,8_2^3,8_3^2,1^4,\]
and that the composition factors of $L(G)\downarrow_L$ are
\[ (20,20^*)^2,14^4,6^8,(4,4^*)^7,1^8.\]

\noindent \textbf{Step 2: $8_3/48/8_3$ is a summand of $L(G)\downarrow_H$.} We concentrate on the non-principal block, which has factors $48,8_3^2$. As $\Ext^1(48,8_3)$ is $1$-dimensional, and there are no self extensions of $8_3$ or $48$, we get a unique uniserial module of structure $8_3/48/8_3$, and so this summand of $L(G)\downarrow_H$ is one of
\[ 8_3/48/8_3,\qquad 8_3\oplus 8_3\oplus 48,\]
since it is self dual. We understand completely the restrictions of these modules to $L$, and they are
\[ 1/14/6,6/4,4^*/1,1,6/6,14/1,\qquad (1/6/1)^{\oplus 2}\oplus 14/6/4,4^*/6/14.\]
Notice that, as we may assume that $L(G)\downarrow_H$ has no trivial submodules, these two modules describe the summand of $L(G)\downarrow_L$ with trivial submodules.

Suppose that we are in the second case, i.e., $8_3$ is a summand of $L(G)\downarrow_H$. We now see that $14/1$ cannot be a submodule of $L(G)\downarrow_L$: to see this, let $\rho_1$ and $\rho_2$ denote the projection maps along the non-principal and principal block summands $M_1$ and $M_2$ of $L(G)\downarrow_H$. If $M=1/14/1$ is isomorphic to a submodule of $L(G)\downarrow_L$, then the sum of the images along $\rho_1$ and $\rho_2$ must equal all of $M$. However, as $M$ is uniserial, this means that one of $\rho_1$ and $\rho_2$ is an injection for $M$, hence $M$ is contained in $M_1$, but from the actions above we see that it is not. Similarly, any module that is an extension with quotient $14/1,6$ (indecomposable) and submodule a module with factors $4$ and $4^*$ cannot appear in $L(G)\downarrow_L$ either. With these facts in mind we try to understand \emph{all} $\Alt(8)$s inside $E_8$ in characteristic $2$, running through positive-dimensional subgroups and using Proposition \ref{prop:alt8compfactors} and its proof. We are only interested in those with factors as given above, so we only need concern ourselves with those in the $A_7$ parabolic with factors $4^2$ or $4,4^*$ on the natural, $D_7$ parabolic with factors $6,4,4^*$ on the natural, and $D_8$ with factors $(4,4^*)^2$ on the natural. Thus we may assume that $L$ lies in one of these

\bigskip

If $L$ lies inside the $A_7$ parabolic $X$ and its factors are $4,4^*$, then the action on the natural is either $4\oplus 4^*$ or $4/4^*$ (up to outer automorphism). The filtration of the action of $X$ on $L(G)$ is (up to duality)
\[ L(\lambda_7)/L(\lambda_2)/L(\lambda_3)/L(\lambda_1+\lambda_7)/L(\lambda_3)/L(\lambda_6)/L(\lambda_1),\]
and so we are interested in $L(\lambda_6)=\Lambda^2(L(\lambda_1))^*$. In the two cases this module is $1,6/14/1,6$ and $1/14/1\oplus 6^{\oplus 2}$, so has a submodule $14/1,6$ (indecomposable) or $14/1$. Of course, the layer below this has factors $4$ and $4^*$, so has no extensions with $1$ or $14$, so we have one of the two submodules of $L(G)\downarrow_L$ not allowed by the above arguments.

Alternatively, $L$ acts as $4\oplus 4$ on the natural, but then all of $L(\lambda_1)$, $L(\lambda_2)$, $L(\lambda_3)$ and their duals have no trivial composition factors. Hence the $\{1\}'$-residual, modulo its $\{1\}'$-radical, is a subquotient of $L(\lambda_1+\lambda_7)$, but this module is clearly
\[ (4\otimes 4^*)^{\oplus 4}=(1/14/1)^{\oplus 4},\]
so we can never have $1/6/1$ being a submodule of $L(G)\downarrow_L$. Thus $L$ does not lie in the $A_7$ parabolic.

\medskip

If $L$ lies inside $X=D_8$ with factors $(4,4^*)^2$ on the natural module $V$, then $\Lambda^2(V)$ is a summand of $L(G)\downarrow_X$, and so if $4\oplus 4^*$ is a submodule of $V\downarrow_L$ then $4\otimes 4^*=1/14/1$ is a submodule of $L(G)\downarrow_L$, not allowed. Hence $4\oplus 4$ is the socle of $V\downarrow_L$. (The $\{4,4^*\}$-radical of $P(4)$ is $4^*/4$.) This means that $L$ acts on $V$ as $(4^*/4)^{\oplus 2}$. The exterior square of this module has three summands, two of which are $\Lambda^2(4^*/4)$, which we earlier saw had $14/1,6$ as a submodule, so not allowed in $L(G)\downarrow_L$, and $L\not\leq X$ again.

\medskip

We are left with $X$ being the $D_7$ parabolic, with $L$ acting on the natural with factors $6,4,4^*$. Certainly $L$ stabilizes a $4$-space on the natural module, and acts as $\GL_4(2)$ on it, so the stabilizer in $X$ acts irreducibly on the $4$-space.

This means that the radical of the form on the $4$-space is either $0$ or the whole space, so that it is either non-singular or totally isotropic. These stabilizers are a maximal parabolic or a $D_2D_5$ subgroup, which does not contain $\GL_4(2)$, and so we lie in a parabolic subgroup of $X$, hence a different parabolic subgroup of $G$. However, the parabolic subgroups of $G$ contained in the $D_7$ parabolic are all (after $A_1$s and $A_2$s have been stripped out) in the $E_7$ or $A_7$ parabolics, and so have been dealt with.

This proves that we cannot have $8_3^{\oplus 2}\oplus 48$ as a summand of $L(G)\downarrow_H$. We note that all unipotent elements in $H$ act projectively on $8_3/48/8_3$.

\medskip

\noindent \textbf{Step 3: $26$ is not a submodule of $L(G)\downarrow_H$.} We come back to the principal block summand, which has composition factors $26^2,(20,20^*)^2,8_1^3,8_2^3,1^4$. The element $v=(1,2,3,4)(5,6)$ acts projectively on $8_1$, $8_2$, $20$ and $20^*$, and as $4^5,2^3$ on the $26$. Thus $v$ has at least fifty-eight blocks of size $4$ in its action on $L(G)$, so must act as $4^{60},2^4$ by Table \ref{tab:unipe8p4}.

Write $W$ for the $\{8_1,8_2\}$-residual modulo its $\{8_1,8_2\}$-radical, so $W$ is a submodule of a sum of $P(26)$, $P(20)$ and $P(20^*)$s. Remove also any simple summands from $W$.

Suppose that $26$ lies in the socle of $W$. The $\{26,78\}'$-radical of $P(26)/\soc(P(26))$, lifted back to $P(26)$, has the form
\[ 20,20^*/1,1,8_1,8_2/26.\]
The $\{26,78\}'$-radical of $P(20)$ is
\[ 8_2/20^*/1,8_1/20,\]
and so if $26$ lies in $\soc(W)$, then $\soc(W)=20^{\oplus 2}\oplus 26$ or $20\oplus 20^*\oplus 26$. In both cases, the action of $v$ on the submodule
\[ 1,1,1,1/26,20,20^\pm\]
is projective, plus $3,2^2,1^3$. Since $v$ acts projectively on all of $8_1$, $8_2$, $20$ and $20^*$, and acts on $26$ as $4^5,2^3$, extending this module (possibly with some $8_1$s and $8_2$s as well) by the $26\oplus 20\oplus 20^\mp$ on the top of $W$ requires us to, for the action of $v$, add $2^3$ onto $3,2^2,1^3$ to make $4^2,2^4$.

We claim that, for the cyclic group of order $4$, there is no submodule $N$ of $M=4^2,2^4$ of shape $3,2^2,1^3$, whose quotient is $2^3$. To see this, firstly note that the socles of $M$ and $N$ coincide, so we may work modulo $\soc(M)$. Now the socle of $M$ has dimension $6$ and the socle of $N$ has dimension $3$, so the extension of $N$ by $\soc(M/N)$ must be split; loosely speaking, the socle of $M/N$ is `contained in' the socle of $N$ (although of course this doesn't strictly make sense). We therefore have an extension of $1^3$ by $2,1^5$ to make $3^2,1^4$ (remember we are working modulo $\soc(M)$), but $1^3$ is a trivial module, so that the kernel of the map from $3^2,1^4$ to $1^3$ contains the Jacobson radical of $3^2,1^4$, namely $2^2$, a contradiction.

What we have therefore proved is that $26$ is not in $\soc(W)$, and hence the socle of $W$ is one of $20$, $20^{\oplus 2}$ and $20\oplus 20^*$.

\medskip

\noindent \textbf{Step 4: Eliminating the remaining socle possibilities.} If the socle is $20$ then $W$ is a submodule of $P(20)$, and we firstly take the $\{78\}'$-radical of $P(20)$, then take the $\{20^*\}'$-residual of that (since $\soc(W)=20$, $\top(W)=20^*$), to produce the self-dual module
\[ 20^*/1,8_2/20,26/1,1,8_1/20^*,26/1,8_2/20.\]
However, $v$ acts projectively on this, and since it acts projectively on the non-principal block part of $L(G)$ as well, this means that it would act projectively on all of $L(G)$, not allowed by Table \ref{tab:unipe8p4}, as $v$ is supposed to act as $4^{60},2^4$. Thus the socle isn't simple.

Thus the socle is $20\oplus 20^\pm$, and in particular all $20,20^*$s are taken up by the socle and top of $W$. We thus want to take the preimage in $P(20)/\soc(P(20))$ of the $\{1,8_1,8_2,26\}$-radical $M$, which is
\[ 1,8_1/26/1,8_2/20.\]
In theory this could work, but we will need to place $20$ or $20^*$ on top of this module and hide the trivial quotient. Firstly, $\Ext^1(20,M)=0$, so this will not work; $\Ext^1(20^*,M)$ is $1$-dimensional, but this module is
\[ 1,8_1/26,20^*/1,8_2/20,\]
and so the socle cannot consist of $20^\pm$s either. Thus $H$ fixes a line on $L(G)$, as needed.
\end{pf}

\section{$\Alt(10)$ and above}
\label{sec:alt10}
For $H\cong\Alt(n)$ for $n\geq 10$, Litterick in \cite{litterickmemoir} proved that $H$ lies inside a $\sigma$-stable positive-dimensional subgroup for all $G$, except possibly a single set of composition factors for $p=2$, $G=E_8$ and $H\cong\Alt(10)$, which was completed in the published version of his PhD thesis using a complicated argument. Using pressure and Frobenius reciprocity we can provide a shorter proof of this result.

The composition factors of $L(G)\downarrow_H$ are
\[ 48^2,26^4,8^5,1^8,\]
and since $H^1(H,48)=0$ and $H^1(H,26)=H^1(H,8)$ are $1$-dimensional we have that $L(G)\downarrow_H$ has pressure $1$. If we remove any $48$s in the socle and top of $L(G)\downarrow_H$ to get a module $W$, then since $W$ has pressure $1$ by Proposition \ref{prop:lowpressuremodules} we have that $\soc(W)$ is either $8$ or $26$, so that $W$ is a submodule of either $P(8)$ or $P(26)$.

Using Lemma \ref{lem:frobrep}, since $L(G)\downarrow_L$ has a trivial submodule (where $L=\Alt(9)$), there is a non-trivial map from $1_L\uparrow^H=1/8/1$ to $L(G)\downarrow_H$, and hence to $W$. Thus either $H$ fixes a line on $L(G)$, as needed, or $1/8$ is a submodule of $L(G)\downarrow_H$, and so $W\leq P(8)$.

Furthermore, since $W$ has pressure $1$, $W$ must have at least seventeen socle layers, since all $1$s, $8$s and $26$s must lie in different layers. The projective $P(8)$ has exactly nineteen socle layers, and is given in \cite{craven2015un3}. The seventeenth socle layer of $P(8)$ is $26\oplus 200$ and the eighteenth is $1$, which therefore means that if $\soc(W)=8$ then $W$ cannot be self dual. Hence $H$ fixes a line on $L(G)$, so is not maximal.

\bigskip
\bigskip

\noindent\textbf{Acknowledgements:} I would like to thank Kay Magaard and Chris Parker for helpful discussions throughout, and Cheryl Praeger and Michael Giudici for asking me to prove Theorem \ref{thm:alt5}. Thank you also to Darrin Frey and Alastair Litterick for letting me have preprint versions of their work. I would like to thank the Royal Society for supporting me financially with a University Research Fellowship.

\bibliographystyle{amsplain}
\bibliography{references}

\providecommand{\bysame}{\leavevmode\hbox to3em{\hrulefill}\thinspace}
\providecommand{\MR}{\relax\ifhmode\unskip\space\fi MR }
\providecommand{\MRhref}[2]{%
  \href{http://www.ams.org/mathscinet-getitem?mr=#1}{#2}
}
\providecommand{\href}[2]{#2}
\begin{thebibliography}{10}

\bibitem{alperin}
Jonathan Alperin, \emph{Local representation theory}, Cambridge Studies in
  Advanced Mathematics, vol.~11, Cambridge University Press, Cambridge, 1996.

\bibitem{aschbacherE6Vun}
Michael Aschbacher, \emph{The maximal subgroups of ${E}_6$}, preprint, 170pp.

\bibitem{craven2015un3}
David~A.\ Craven, \emph{The structure of projective indecomposable modules for
  {$A_n$}, $n\leq 12$}, preprint, arXiv:1605.04403.

\bibitem{cmp2015un}
David~A.\ Craven, Kay Magaard, and Christopher Parker, \emph{On generic
  subgroups of exceptional groups}, in preparation.

\bibitem{feit}
Walter Feit, \emph{The representation theory of finite groups}, North--Holland,
  Amsterdam--New York, 1982.

\bibitem{freyun}
Darrin Frey, \emph{Embeddings of ${A}lt_n$ and its perfect covers for $n\geq 6$
  in exceptional complex {L}ie groups}, to appear.

\bibitem{frey1998}
\bysame, \emph{Conjugacy of {$Alt_5$} and {$SL(2,5)$} subgroups of
  {$E_6(\mathbb{C})$}, {$F_4(\mathbb{C})$}, and a subgroup of
  {$E_8(\mathbb{C})$} of type {$A_2E_6$}}, J. Algebra \textbf{202} (1998),
  414--454.

\bibitem{hassanabadi1978}
A.~Mohammadi Hassanabadi, \emph{Automorphisms of permutational wreath
  products}, J.\ Austral.\ Math.\ Soc. \textbf{26} (1978), 198--208.

\bibitem{houghton1975}
Chris~H.\ Houghton, \emph{Wreath products of groupoids}, J.\ London Math.\ Soc.
  \textbf{10} (1975), 179--188.

\bibitem{lawther1995}
Ross Lawther, \emph{Jordan block sizes of unipotent elements in exceptional
  algebraic groups}, Comm.\ Algebra \textbf{23} (1995), 4125--4156.

\bibitem{lawther2009}
\bysame, \emph{Unipotent classes in maximal subgroups of exceptional algebraic
  groups}, J.\ Algebra \textbf{322} (2009), 270--293.

\bibitem{lawthertesterman1999}
Ross Lawther and Donna Testerman, \emph{${A}_1$ subgroups of exceptional
  algebraic groups}, Mem.\ Amer.\ Math.\ Soc. \textbf{141} (1999), no.~674,
  viii+131.

\bibitem{liebecksaxl1987}
Martin Liebeck and Jan Saxl, \emph{On the orders of maximal subgroups of the
  finite exceptional groups of {L}ie type}, Proc.\ London Math.\ Soc.
  \textbf{55} (1987), 299--330.

\bibitem{lst1996}
Martin Liebeck, Jan Saxl, and Donna Testerman, \emph{Simple subgroups of large
  rank in groups of {L}ie type}, Proc.\ London Math.\ Soc. \textbf{72} (1996),
  425--457.

\bibitem{liebeckseitz1999}
Martin~W. Liebeck and Gary~M. Seitz, \emph{On finite subgroups of exceptional
  algebraic groups}, J.\ reine angew. Math. \textbf{515} (1999), 25--72.

\bibitem{liebeckseitz2004}
\bysame, \emph{The maximal subgroups of positive dimension in exceptional
  algebraic groups}, Mem.\ Amer.\ Math.\ Soc. \textbf{169} (2004), no.~802,
  vi+227.

\bibitem{liebeckseitz2005}
\bysame, \emph{Maximal subgroups of large rank in exceptional groups of {L}ie
  type}, J.\ London.\ Math.\ Soc. \textbf{71} (2005), 345--361.

\bibitem{litterick}
Alastair Litterick, \emph{Finite simple subgroups of exceptional algebraic
  groups}, Ph.D. thesis, Imperial College, London, 2013.

\bibitem{litterickmemoir}
\bysame, \emph{Finite simple subgroups of exceptional algebraic groups}, Memo.\
  Amer.\ Math.\ Soc., to appear., 2015.

\bibitem{lusztig2002}
George Lusztig, \emph{Homomorphisms of the alternating groups ${A}_5$ into
  reductive groups}, J.\ Algebra \textbf{260} (2002), 298--322.

\bibitem{magaardphd}
Kay Magaard, \emph{The maximal subgroups of the {C}hevalley groups {${F}_4(F)$}
  where {${F}$} is a finite or algebraically closed field of characteristic not
  equal to 2 or 3}, Ph.D. thesis, California Institute of Technology, 1990.

\bibitem{seitz1991}
Gary~M. Seitz, \emph{Maximal subgroups of exceptional algebraic groups}, Mem.\
  Amer.\ Math.\ Soc. \textbf{90} (1991), no.~441, iv+197.

\bibitem{thomas2016}
Adam Thomas, \emph{Irreducible ${A}_1$ subgroups of exceptional algebraic
  groups}, J. Algebra \textbf{447} (2016), 240--296.

\bibitem{wilsonrob}
Robert~A.\ Wilson, \emph{The finite simple groups}, Graduate Texts in
  Mathematics, vol. 251, Springer-Verlag London Ltd., London, 2009.

\end{thebibliography}

\end{document}